\documentclass[11pt]{report} 
\usepackage{amsmath}
\usepackage{amsfonts}
\usepackage[latin1]{inputenc}
\usepackage[T1]{fontenc}
\usepackage[french]{babel}
\makeatletter
\def\@sect#1#2#3#4#5#6[#7]#8{%
  \ifnum #2>\c@secnumdepth
    \let\@svsec\@empty
  \else
    \refstepcounter{#1}%
    \protected@edef\@svsec{\@seccntformat{#1}\relax}%
  \fi
  \@tempskipa #5\relax
  \ifdim \@tempskipa>\z@
    \begingroup
      #6{%
        \@hangfrom{\hskip #3\relax\@svsec}%
          \interlinepenalty \@M #8\@@par}%
    \endgroup
    \csname #1mark\endcsname{#7}%
    \addcontentsline{toc}{#1}{%
      \ifnum #2>\c@secnumdepth \else
        \protect\numberline{\csname the#1\endcsname.}%
      \fi
      #7}%
  \else
    \def\@svsechd{%
      #6{\hskip #3\relax
      \@svsec #8}%
      \csname #1mark\endcsname{#7}%
      \addcontentsline{toc}{#1}{%
        \ifnum #2>\c@secnumdepth \else
          \protect\numberline{\csname the#1\endcsname.}%
        \fi
        #7}}%
  \fi
  \@xsect{#5}}
\def\@seccntformat#1{\csname the#1\endcsname.\quad}

\def\@begintheorem#1#2{\trivlist
   \item[\hskip \labelsep{\bfseries #1\ #2.}]\itshape}
\def\@opargbegintheorem#1#2#3{\trivlist
      \item[\hskip \labelsep{\bfseries #1\ #2\ (#3).}]\itshape}
\renewcommand\theequation{\thesection.\arabic{equation}}
\@addtoreset{equation}{chapter}
\makeatother

\newtheorem{theo}[equation]{Th\'eor\`eme}

\newtheorem{pr}[equation]{Proposition}
\newtheorem{fait}[equation]{Fait}
\newtheorem{lem}[equation]{Lemme}
\newtheorem{cor}[equation]{Corollaire}
\newtheorem{df}[equation]{D\'efinition}
\newtheorem{question}[equation]{Question}
\newcommand{\smat}[1]{\left(\begin{smallmatrix}#1\end{smallmatrix}\right)}
\newcommand{\carrenoir}{\rule{0.5em}{0.5em}}
\newenvironment{demo}[1]{\textbf{D\'emonstration #1 : }}
{ \hfill\carrenoir\nolinebreak\vspace{2mm}}
\newenvironment{exemple}{
\vspace{2mm}
\refstepcounter{equation}\noindent\textbf{Exemple \theequation.\ }}%
{\hfill\carrenoir\nolinebreak\vspace{2mm}}

\newenvironment{remarque}{
\refstepcounter{equation}\noindent\textbf{Remarque \thetheo.\ }}
{\nolinebreak\vspace*{1mm}}

\newcommand{\avec}{\mbox{ avec }}
\newcommand{\pour}{\mathrm{\ pour\ }}
\newcommand{\et}{\mbox{ et }}
\newcommand{\ou}{\mbox{ ou }}

\newcommand{\dep}[2]{\partial #1/\partial #2}
\newcommand{\Dep}[2]{\displaystyle\frac{\partial #1}{\partial #2}}
\newcommand{\oper}[2]{\newcommand{#1}{\mathop{\mathrm{#2}}\nolimits} }

\oper{\Ker}{Ker}
\oper{\Ima}{Im}
\oper{\Id}{Id}
\oper{\Ric}{Ric}
\oper{\Vol}{Vol}
\oper{\Hess}{Hess}
\oper{\Tr}{Tr}
\oper{\diam}{diam}
\oper{\injrad}{inj}
\oper{\dimension}{dim}
\oper{\codimension}{codim}
\oper{\Det}{Det}
\oper{\Aut}{Aut}
\newcommand{\de}{\mathrm{ d }}

\newcommand{\trans}[1]{\vphantom{#1}^t#1}
\oper{\LC}{D}
\DeclareFontFamily{OML}{eur}{\skewchar\font127}
\DeclareFontShape{OML}{eur}{m}{n}{<5> <6> <7> <8> <9> gen * eurm
<10><10.95><12><14.4><17.28><20.74><24.88>eurm10}{}
\DeclareSymbolFont{greek}{OML}{eur}{m}{n}
\DeclareMathSymbol{\codiff}{\mathord}{greek}{"0E}
\newcommand{\N}{\mathbb N}
\newcommand{\Z}{\mathbb Z}
\newcommand{\R}{\mathbb R}

\newcommand{\Hom}{\mathrm{Hom}}
\newcommand{\Aff}{\mathrm{Aff}}
\oper{\ad}{\mathrm{ad}}
\oper{\GL}{GL}
\oper{\SL}{SL}
\oper{\SO}{SO}
\def\vvdots{\mathinner{\mkern-.5mu\raise0pt\hbox{.}\mkern-6.6mu\raise4pt\hbox{.}
\mkern-6.5mu\raise8pt\hbox{.}}}

\def\ccdot{\mathinner{\mkern-.5mu\raise2pt\hbox{.}}}

\def\hhdots{\mathinner{\mkern-.5mu\raise2pt\hbox{.}\raise2pt\hbox{.}\raise2pt\hb
ox{.}}}
\newcommand{\vecdots}[2]{\begin{smallmatrix}#1\\\vvdots\\#2\\ \end{smallmatrix}}

\title{Sur le spectre du laplacien des fibrés en tores qui s'effondrent}
\date{}
\author{Pierre Jammes}
\begin{document}
\maketitle
\pagenumbering{Roman}
\renewcommand{\thepage}{\textsc{\roman{page}}}
\renewcommand{\theequation}{\arabic{equation}}
\chapter*{Introduction}

L'une des questions centrales de l'analyse sur les variétés est
l'estimation du spectre du laplacien agissant sur les fonctions 
d'une variété compacte en 
fonction d'invariants géométriques, et fait l'objet d'un vaste
littérature. Des résultats importants ont en particulier été
obtenus concernant la minoration de la première valeur propre du
laplacien. J.~Cheeger en a donné une estimation en fonction d'une
constante isopérimétrique:
\begin{df}
Soit $(M,g)$ une variété riemannienne compacte connexe de dimension $n$. 
On définit la constante de Cheeger $h$ de la variété $(M,g)$ par 
$$h(M,g)=\inf_S\frac{\Vol S}{\min(\Vol M_1,\Vol M_2)},$$ 
où $S$ parcourt l'ensemble des sous-variétés fermées de dimension
$n-1$ de $M$ qui partitionnent $M$ en deux sous-variétés $M_1$ et
$M_2$ dont le bord est $S$.
\end{df}
\begin{theo}[\cite{cr}]\label{intro:th1}
Si $(M,g)$ une variété compacte connexe, alors $\lambda_1(M,g)\geq
\frac{h(M,g)^2}4$, où $\lambda_1(M,g)$ est la première valeur propre 
du laplacien agissant sur les fonctions de $M$.
\end{theo}
Si on ajoute une hypothèse sur la courbure de Ricci, on peut encadrer
la première valeur propre en fonction de la constante de Cheeger :
\begin{theo}[\cite{buser}]
Soit $(M,g)$ une variété compacte connexe de dimension $n$. 
Si $\Ric(M,g)>-(n-1)ag$ alors $\lambda_1(M,g)\leq \tau(n)(\sqrt ah(M,g)+
h(M,g)^2)$, où $a$ est une constante positive et $\tau(n)$ une
constante ne dépendant que de la dimension.
\end{theo}

Une autre approche est de minorer cette valeur propre en fonction de
bornes sur la courbure et le diamètre de la variété. Un premier
résultat a été obtenu par M.~Gromov:
\begin{theo}[\cite{gr}]\label{intro:th2}
Soit $a$ et $d$ deux réels strictement positifs et $n$ un entier.
Il existe une constante $c(n,a,d)>0$ telle que si $(M,g)$ est
une variété riemannienne de dimension $n$ dont le diamètre et la courbure
de Ricci vérifient $\diam(M,g)\leq d$ et $\Ric(M,g)\geq -ag$, alors
$$\lambda_1(M,g)\geq c.$$
\end{theo}
 Il existe d'autres contributions à ce type de minoration (\cite{ly80}, 
\cite{bbg}), et l'estimation du théorème \ref{intro:th2} a été affinée.
En particulier, Cheng et Zhou donnent dans \cite{cz} :
$$c=\frac{\pi^2}{d^2}e^{-\frac12c_n\sqrt{ad^2}}, \avec c_n=
\max\{\sqrt{n-1},\sqrt2\}.$$

 Le spectre du laplacien de Hodge-de~Rham agissant sur les formes
différentielles a été moins étudié, et la question se pose 
naturellement de savoir si les théorèmes \ref{intro:th1} et
\ref{intro:th2} se généralisent aux formes. Rappelons que si $(M,g)$ est
une variété connexe compacte orientable, le laplacien est défini
sur l'ensemble $\Omega^p(M)$ des formes différentielles de
degré $p$ de $M$, $0\leq p\leq n$ par 
$$\Delta^p=\de\codiff+\codiff\de,$$ où $\de$ est
la différentielle extérieure et $\codiff$ son adjoint formel pour 
la norme $L^2$ sur $\Omega^p(M)$. La dimension du noyau de $\Delta^p$
est égale au $p$-ième nombre de Betti de $M$, et son spectre est discret.
On notera 
$$0=\lambda_{p,0}(M,g)<\lambda_{p,1}(M,g)\leq\lambda_{p,2}(M,g)
\leq\cdots\leq\lambda_{p,k}(M,g)\leq\cdots$$
 les valeurs propres de
$\Delta^p$, en répétant les valeurs propres non nulles s'il y a 
multiplicité. Si $p=0$, on retrouve le laplacien agissant sur les fonctions.

 L'étude du laplacien de Hodge-de~Rham montre que le comportement
du spectre pour $1\leq p\leq n-1$ est différent du spectre du laplacien
agissant sur les fonctions. En particulier, le théorème \ref{intro:th2}
ne se généralise pas aux formes différentielles, même avec 
l'hypothèse plus forte de courbure sectionnelle bornée.
En effet, C.~Colbois
et G.~Courtois ont donné dans \cite{cc1} des exemples de variétés
admettant des petites valeurs propres, c'est-à-dire qu'on peut faire
varier leur métrique de sorte que la ou les premières valeurs propres
non nulles tendent vers zéro, le diamètre et la courbure sectionnelles
restant bornés.

Le théorème \ref{intro:th2} précédent soulève la question
de savoir à quelles conditions une variété admet une petite valeur
propre à diamètre et courbure bornés. Un début de réponse a été donné
par B.~Colbois et G.~Courtois dans \cite{cc1} en montrant qu'on
peut obtenir pour les formes différentielles un résultat semblable au
théorème \ref{intro:th2} si on se donne une hypothèse supplémentaire
sur le rayon d'injectivité de la variété:
\begin{theo}\label{intro:th3}
Pour tous réels $a$, $d$, et $r$ strictement positifs et tout entier
$n$, il existe une constante $c'(n,a,d,r)>0$ telle que si
$(M,g)$ est une variété riemannienne compacte connexe de dimension $n$
vérifiant $\diam(M,g)\leq d$, $|K(M,g)|\leq a$ et $\injrad(M,g)\geq r$,
où $K(M,g)$ et $\injrad(M,g)$ désignent respectivement  la courbure 
sectionnelle et le rayon d'injectivité de $M$, alors 
$$\lambda_{p,1}(M,g)\geq c',$$
pour tout $p$.
\end{theo}
Ce résultat a été amélioré par S.~Chanillo et F.~Trèves, qui donnent
une valeur explicite de la constante $c'$ du théorème \ref{intro:th3},
et en mettant en évidence le rôle du rayon d'injectivité dans cette
constante (\cite{ct}):
\begin{equation}\label{intro:eq2}
\lambda_{p,1}(M,g)\geq c''(n,a,d)\cdot r^{4n^2+4n-2}.
\end{equation}
Une conséquence immédiate est que si $(M_i^n,g_i)$ est une suite
de variétés riemanniennes de dimension $n$ de diamètres et courbures
uniformément bornés telle que $\lambda_{p,1}(M_i,g_i)$ tende vers zéro
pour un certain $p$, $1\leq p\leq n$,
alors $\injrad(M_i,g_i)$ tend aussi vers zéro (remarque: si le diamètre
et la courbure sectionnelle sont bornés, minorer le rayon d'injectivité est
équivalent à minorer le volume). Plus simplement, si
on se donne une variété compacte $M$ et qu'on fait varier sa métrique
à diamètre et courbure bornée de sorte que $\lambda_{p,1}$ tende vers zéro,
alors son rayon d'injectivité tend vers zéro, c'est-à-dire qu'elle
s'effondre. Cette condition a des conséquences topologiques importantes.
En particulier, si $(M,g_i)$ tend pour la distance de Gromov-Haussdorf
vers une variété de dimension inférieure, K.~Fukaya a montré (entre
autres choses) dans \cite{fu2} et \cite{fu3} que $M$ possède une structure 
de fibré dont la base est la variété limite:
\begin{theo}\label{intro:th4}
Soit $(M_i,g_i)$ une suite de variétés compactes de dimension $n$ et
$(N,h)$ une variété riemannienne compacte de dimension $m<n$. Si pour 
tout $i$ la courbure sectionnelle de $M_i$ vérifie $|K(M,g_i)|\leq 1$,
et si $(M,g_i)$ converge 
vers $(N,h)$ pour la distance de Gromov-Hausdorff, alors pour tout
$i$ suffisamment grand il existe une fibration $\pi_i:M_i\rightarrow N$ 
dont la fibre est une infranilvariété.
\end{theo}

 Cependant, il n'est pas suffisant que le rayon d'injectivité tende
vers zéro pour pour que la première valeur propre non nulle
du laplacien tende vers zéro. On peut par exemple considérer le
produit riemannien d'une variété $(N,h)$ quelconque munie d'une
métrique fixée avec un tore plat dont on fait tendre le diamètre 
vers zéro: le rayon d'injectivité du produit tend vers zéro, mais
son spectre est la réunion des spectres de $(N,h)$ et du tore, et
donc ne contient pas de petites valeurs propres. Ce la motive la 
\begin{question}\label{intro:qu1}
À quelles conditions une variété qui s'effondre
admet-elle une ---~ou plusieurs~--- petites valeurs propres?
\end{question}

Par ailleurs, la minoration (\ref{intro:eq2}) n'est \emph{a priori}
pas optimale, en particulier en ce qui concerne l'exposant du rayon
d'injectivité. Il se pose donc aussi la
\begin{question}\label{intro:qu2}
Peut-on améliorer la minoration de la première valeur propre non nulle
du laplacien de Hodge-de~Rham?
\end{question}

 Des réponses précises aux questions \ref{intro:qu1} et \ref{intro:qu2}
ont été apportées dans des situations topologiques ou géométriques
particulières. C'est le cas des limites adiabatiques :
si $(M,g)$ est une variété 
riemannienne compacte, $A$ une distribution de sous-espaces de
$TM$ et $B$ la distribution orthogonale à $A$, on peut écrire
la métrique $g$ sous la forme $g=g_A\oplus g_B$ où $g_A$ et $g_B$ 
sont des métriques sur $A$ et $B$ respectivement, et
définir sur $M$ la famille de métriques $$g_t=t^2g_A\oplus g_B.$$ 
La limite de $(M,g_t)$ quand $t\rightarrow 0$ est appelée
<<~limite adiabatique~>>. Dans le cas où $M$ est un fibré riemannien, la
distribution $A$ est la distribution verticale $T^VM$, $B$
la distribution horizontale $T^HM$, et on fait
tendre le diamètre de la fibre vers zéro par des
homothéties de rapport $t$. Pour un tel effondrement,
Mazzeo et Melrose ont montré (\cite{mm}) que
le nombre de valeurs propres de l'ordre de $t^2$ quand $t$ tend
vers zéro peut se calculer à l'aide 
de la suite spectrale de Leray du fibré (notons que Forman, \'Alvarez
et Kordyukov ont généralisé ce résultat aux feuilletages riemanniens
dans \cite{fo} et \cite{alk}). Cependant, dans une
situation de limite adiabatique, la courbure n'est en général
pas bornée, en particulier si la fibre n'est pas un tore. D'autre
part, on peut en général effondrer le fibré sur sa base autrement 
que par homothétie de la fibre.

 Par ailleurs, B.~Colbois et G.~Courtois ont étudié dans \cite{cc2} le cas
où $M$ est une variété de dimension $n+1$ qui tend pour la distance de
Gromov-Hausdorff vers une variété $(N,h)$ de dimension $n$ et donnent
des estimations plus précises de la ou les premières valeurs
propres non nulles en fonction de la topologie de $M$ et de la
géométrie de l'effondrement. Selon le théorème \ref{intro:th4}, 
on sait que $M$ est nécessairement un fibré en cercle
sur $N$, et l'hypothèse d'orientabilité sur $M$ assure que ce 
fibré est principal. La topologie de ce fibré est caractérisée par sa
classe d'Euler $[e]\in H^2(N,\Z)$ (cf. \cite{bt}, p.~72). En notant
$e(M)$ le représentant harmonique de la classe de
cohomologie $[e]$ pour la métrique $h$ sur $N$, on a alors:
\begin{theo}\label{vol:th1}
Soit $a$ et $d$ deux réels strictement positifs et $(N,h)$ une 
variété riemannienne de dimension $n$. Il existe
des constantes $\varepsilon_0(n,a,d,(N,h))>0$ et $C_i(n,a,d,(N,h))>0$
pour $i=1,2,3$ telles que si $(M,g)$ est une variété riemanienne de
dimension $n+1$ vérifiant $\diam(M,g)\leq d$, $|K(M,g)|\leq a$
et $\varepsilon=d_{GH}((M,g),(N,h))\leq\varepsilon_0$, alors on a, pour
$1\leq p\leq n$,
\renewcommand{\theenumi}{\theequation.\arabic{enumi}}
\renewcommand{\labelenumi}{\theenumi.}
\begin{enumerate}
\item\label{vol:th1:1}$\displaystyle\lambda_{p,m_p+1}(M,g)\geq C_1$;
\item Si $e\neq0$, alors $\displaystyle C_2\|e(M)\|_2^2\varepsilon^2\leq
\lambda_{1,1}(M,g) \leq C_3\|e(M)\|_2^2\varepsilon^2$;
\item\label{vol:th1:3}Si $\dim H^2(N,\R)=1$, alors
$$ C_2\|e(M)\|_2^2\varepsilon^2\leq
\lambda_{p,k}(M,g) \leq C_3\|e(M)\|_2^2\varepsilon^2 \pour 1\leq k\leq m_p\ ;$$
\end{enumerate}
avec $m_p=b_p(N)+b_{p-1}(N)-b_p(M)$.
\end{theo}
\begin{remarque}
Si $e=0$, par exemple si le fibré est trivial, on a $m_p=0$ pour tout $p$. 
Le théorème se réduit alors à \ref{vol:th1:1}. D'autre part, si $e\neq0$ 
alors $m_1=1$. $\lambda_{1,1}(M,g)$ est donc la seule petite valeur propre
pour les $1$-formes.
\end{remarque}

\begin{remarque}
B.~Colbois et G.~Courtois montrent aussi dans \cite{cc2} que si 
$\dim H^2(N,\R)\geq1$, on ne peut pas trouver de constante $C_2$ vérifiant
\ref{vol:th1:3} qui soit indépendante de la topologie de $M$.
\end{remarque}

\begin{remarque}
Le théorème isole les rôles de la topologie (par l'intermédiaire de
$e(M)$) et du rayon d'injectivité (qui est de l'ordre de $\varepsilon$
quand $\varepsilon$ tend vers zéro) dans la minoration de la première
valeur propre non nulle du laplacien. On voit donc que dans la situation
d'un fibré en cercle qui s'effondre, si le laplacien admet une petite
valeur propre, elle se comporte asymptotiquement comme le carré
du rayon d'injectivité.
\end{remarque}

 Dans des situations topologiques et géométriques plus générales, on
ne connait pas de résultat semblable. Cependant, en ce qui concerne 
la question \ref{intro:qu1}, J.~Lott a apporté une contribution 
importante dans \cite{lo} et \cite{lo2} en définissant un opérateur
limite pour le laplacien, c'est-à-dire un opérateur dont le spectre
est la limite du spectre de laplacien de Hodge-de~Rham quand la 
variété s'effondre, le nombre de valeurs propres petites ou nulles étant 
alors donné par la multiplicité de la valeur propre nulle dans le spectre
de l'opérateur limite. Dans le chapitre 1, nous exposerons la construction
de cet opérateur limite, ainsi que les conditions nécessaires à
l'existence de petites valeurs propres  que Lott en déduit.

Les chapitres 2 et 3 seront consacrés à l'étude, pour des
fibrés munis d'une structure homogène ---~construits plus précisément 
comme quotients d'un groupe de Lie résoluble $G$ par un réseau cocompact
$\Gamma$~---, du spectre du laplacien $\Delta^p_{inv}$ restreint à l'espace
de dimension finie $\Omega^p(M)^G$ des formes différentielles invariantes
lors d'effondrements par des métriques homogènes, afin de mettre en
évidence comment varie précisément le nombre de petites valeurs propres
en fonction de la topologie du fibré et de la géométrie de l'effondrement.
En effet, les petites valeurs propres de $\Delta^p_{inv}$ sont aussi
petites valeurs propres de $\Delta^p$, et dans le cas où le groupe $G$ 
est nilpotent, J.~Lott a montré réciproquement que la recherche de petites 
valeurs propres de $\Delta$ se ramène à l'étude de $\Delta^p_{inv}$:
\begin{pr}[\cite{lo}]\label{pr:lott}
Il existe des constantes $a(n)$, $a'(n)$ et $c(n)$ strictement positives 
telles que si $M$ est une infranilvariété de dimension $n$ munie 
d'une métrique homogène pour laquelle 
$\|R\|_\infty\diam(M)^2\leq a'$, où $\|R\|_\infty$ est la norme
du tenseur de courbure, et si $\alpha$ est une forme
propre du laplacien sur $M$ dont la valeur propre $\lambda$ vérifie
$\lambda<a\diam(M)^{-2}-c\|R\|_\infty$, alors $\alpha$ est
invariante. 
\end{pr}
Cependant, ce résultat ne se généralise pas à tous les groupes
résolubles. On en verra un exemple au paragraphe \ref{cercle:ex2}.

Dans le chapitre 2, nous étudierons des fibrés de fibre $T^n$ sur
la base $S^1$. Le fait qu'une variété qui s'effondre sur un cercle 
admette une structure de solvariété est déjà connu (\cite{pe}, \cite{tu}).
Nous en ferons une construction explicite dans un cas simple:
les fibrés considérés seront définis comme suspension
d'un difféomorphisme linéaire de la fibre $T^n$ représenté
par un élément $A\in\SL_n(\Z)$. De plus, nous ferons 
l'hypothèse simplificatrice que $A$ admet un logarithme réel.
Les principales propriétés du spectre que
nous allons mettre en évidence sont données par le

\begin{theo}\label{cercle:th1}
Soit $n\geq2$, $A\in\SL_n(\Z)$ et $B\in\mathrm M_n(\R)$ tels que
$A=\exp(B)$, $d$ la dimension du sous-espace caractéristique 
associé à la valeur propre $0$ de $B$ et
$d'$ la dimension de son noyau.
Il existe un groupe $G=G(B)\subset\GL_{n+2}(\R)$ et un réseau $\Gamma
\subset G$ tel que $\Gamma\backslash G$ soit homéomorphe au
fibré $M$ de fibre $T^n$ de fibration $p:M\rightarrow S^1$ construit
par suspension du difféomorphisme $A$. Si de plus on suppose que les 
métriques sur $M$ sont homogènes et telles que $p$ soit une 
submersion riemannienne pour une métrique de volume $1$ sur $S^1$, alors:
\renewcommand{\theenumi}{\theequation.\arabic{enumi}}
\begin{enumerate}
\item\label{cercle:th1:1} $\dimension \Ker \Delta^1_{inv}=d'+1$ et 
$\Delta^1_{inv}$ admet $n-d'$ valeurs propres non nulles distinctes
ou non.
\item\label{cercle:th1:2} Pour tout $a>0$, il existe une constante $c(B,a)>0$ 
telle que pour toute métrique invariante $g$ sur $M$ dont la courbure 
sectionnelle vérifie $|K(M,g)|<a$, on a $\lambda^{inv}_{1,i}(M,g)<c$,
pour tout $i$.
\item\label{cercle:th1:3}
Si $d\neq n$, alors pour tout $a>0$, il existe une constante $c'(B,a)>0$ 
telle que pour toute métrique $g$ sur $M$ la courbure 
sectionnelle vérifie $|K(M,g)|<a$, on a $\lambda^{inv}_{1,d-d'+1}(M,g)>c'$. 

Si $d=n$ et $d'\neq n$, alors $G$ est nilpotent, et il existe une suite 
de métriques $g_\varepsilon$ sur $M$ telle que la courbure soit 
uniformément bornée et $\lambda^{inv}_{1,i}(M,g_\varepsilon)\rightarrow0$
quand $\varepsilon\rightarrow0$, pour tout $0<i\leq n-d'$.

Si $d=d'=n$, alors $G=\R^n$, $M$ est un tore et les formes harmoniques
sont exactement les formes invariantes.
\item\label{cercle:th1:4}
Pour tout $k\leq d-d'$, il existe une famille de  métriques 
$g^k_\varepsilon$ sur $M$ de courbure et diamètre uniformément 
bornés par rapport à $\varepsilon$ et une constante $c''(B,k)>0$
telle que $\lambda^{inv}_{1,i}(M,g^k_\varepsilon)\rightarrow 0$ pour 
$i\leq k$ quand $\varepsilon\rightarrow 0$, et
$\lambda^{inv}_{1,k+1}(M,g^k_\varepsilon)>c''$ si $k<n$.
\end{enumerate}
\end{theo}

\begin{remarque}
Le point \ref{cercle:th1:4} montre que le nombre de
petites valeurs propres peut, quand la topologie le permet, fortement
varier avec la géométrie de l'effondrement. Cette situation diffère
donc du cas où la fibre est de dimension $1$, étudié par B.~Colbois
et G.~Courtois dans \cite{cc2}.
\end{remarque}

\begin{remarque}
La démonstration de \ref{cercle:th1:4} met en évidence le fait 
que dans le cas d'un effondrement par homothéties de la fibre (situation
de limite adiabatique), il n'y a pas de petites valeurs propres.
\end{remarque}

\begin{remarque}\label{cercle:rq1}
Une conséquence immédiate du résultat~\ref{cercle:th1:3} est que si 
tout le spectre
de $\Delta_{inv}^1$ tend vers zéro, alors le groupe $G$ est nilpotent.
\end{remarque}

D'autre part, ce théorème permet de donner une condition
nécessaire et suffisante sur $B$ pour l'existence de petites valeurs 
propres pour les $1$-formes:
\begin{cor}\label{cercle:cor1}
Sous les hypothèses du théorème \ref{cercle:th1}, il existe une suite 
de métriques homogènes $g_\varepsilon$ sur $M$ telle
que  $\lambda^{inv}_{1,1}(M,g_\varepsilon)\rightarrow 0$ quand 
$\varepsilon\rightarrow 0$ si et seulement si
$d\neq d'$ (\emph{i.e.} si la réduite de Jordan de $B$ a un  bloc 
nilpotent non nul).
\end{cor}

\begin{remarque}
Les résultats \ref{cercle:th1} et \ref{cercle:cor1} montrent que le 
comportement asymptotique du spectre de $\Delta^1_{inv}$ est 
essentiellement lié à la nature, de la partie
nilpotente de la réduite de Jordan de $B$. En effet, le bloc nilpotent
est caractérisé par les invariants $d$ et $d'$ de $B$.
\end{remarque}

Dans le cas des $p$-formes, $p\geq2$, la situation est plus complexe,
et en particulier le corollaire \ref{cercle:cor1} n'est pas vrai pour $n$
et $p$ quelconque (on verra au paragraphe \ref{cercle:ex1}
qu'il peut exister une petite valeur propre pour les $2$-formes alors
que $B$ n'a pas de bloc nilpotent).
On peut cependant donner une condition nécessaire à l'existence de 
petites valeurs propres, à savoir que $B$ n'est pas semi-simple. 
En effet, on a le:

\begin{theo}\label{cercle:th2}
Sous les hypothèse du théorème \ref{cercle:th1},
si $B$ est semi-simple, alors il existe $c''(B,a)>0$ tel que pour toute 
métrique homogène $g$ sur $M$ dont la courbure sectionnelle vérifie
$|K(M,g)|<a$, on a $\lambda^{inv}_{p,1}(M,g)>c''$.
\end{theo}

Dans le chapitre 3, nous allons nous intéresser à des fibrés en tore
$T^k$ sur le tore $T^2$ muni d'une structure de fibré principal. 
Leur topologie est
assez simple. On sait en effet qu'un fibré principal en tore sur le tore
peut être muni d'une structure de nilvariété (\cite{pe}). On va ici
montrer un résultat plus précis en utilisant le fait que la base est de
dimension 2 : ce fibré est difféomorphe au produit
d'un tore et d'une nilvariété de dimension $3$. Puis on calculera le
spectre du laplacien restreint aux formes différentielles invariantes
en supposant le fibré muni d'une métrique homogène.

On obtient le r\'esultat suivant:

\begin{theo}\label{tore:th}
Soit $M$ un fibré principal non trivial en tore $T^n$ sur le tore $T^2$ . 
Alors
\renewcommand{\theenumi}{\theequation.\arabic{enumi}}
\begin{enumerate}
\item \label{tore:th:1} $M$ est une nilvariété et, si $n\geq2$, 
$M$ est homéomorphe à $N\times T^{n-1}$, où $N$ est une
nilvariété de dimension $3$.

\item \label{tore:th:2}  Il existe un vecteur $V$ vertical (tangent à
la fibre $T^n$) tel que si
$M$ est muni d'une métrique homogène, alors pour tout
$p\in[1,n+1]$,  $\Delta^p_{inv}$
admet une unique valeur propre non nulle $\lambda$. Sa  multiplicité
est $C^{p-1}_n$, et $\lambda=\Vol(B)^{-2}|V|^2$, où $\Vol(B)$ 
est le volume de la base du fibré pour la métrique induite.
\end{enumerate}
\end{theo}

\begin{remarque}
Le produit du \ref{tore:th:1} n'est pas
nécessairement riemannien pour les métriques considérées. Le spectre
ne peut donc pas se déduire directement de la formule de Künneth.
\end{remarque}

\begin{remarque}\label{tore:rem1}
On voit que contrairement à la situation du
théorème \ref{cercle:th1}, un effondrement par homothéties de la fibre
produit une petite valeur propre, et que $\lambda$ est alors
proportionnelle au carré du diamètre de la fibre, à topologie fixée. 
\end{remarque}

\begin{remarque}\label{tore:rem2}
Dans le cas o\`u $n=1$, la remarque pr\'ec\'edente rejoint les
r\'esultats de B.~Colbois et G.~Courtois
qui \'etudient dans \cite{cc2} le spectre des fibr\'es en cercles sur
des bases quelconques et sans restrictions sur la m\'etrique. Mais si
la dimension de la fibre est plus grande ($n\geq2$), un ph\'enom\`ene
nouveau appara\^\i t: il existe dans ce cas des effondrements du
fibr\'e tels que $\lambda$ ne tende pas vers z\'ero. Nous en
donnerons des exemples au paragraphe \ref{tore:eff}.
Si $n\geq2$, le nombre de petites valeurs propres ne d\'epend donc
pas uniquement de la topologie. Cependant, on a pas de libert\'e sur ce
nombre comme en \ref{cercle:th1:4}.
\end{remarque}

Dans une deuxième partie, du chapitre 4 au chapitre 7, nous allons nous
intéresser à la question \ref{intro:qu2}. 
Le théorème \ref{vol:th1} ne considère que des fibrés dont la fibre est 
de dimension fixée égale à $1$ ; on peut se demander dans quelle mesure 
il se généralise aux fibres de dimension plus grande. Nous allons
ici nous intéresser plus précisément aux situation de fibrations
principales s'effondrant sur leur base, la fibre étant alors un tore 
(ce sont les
seules infranilvariétés possédant une structure de groupe).
Notre but dans cette seconde partie est d'arriver au résultat suivant:
\begin{theo}\label{vol:th}
Soit $a$ et $d$ deux réels strictement positifs, $n$ un entier et 
$(N,h)$ une variété riemannienne de dimension strictement inférieure 
à $n$. Il existe des constantes $\varepsilon_0(n,a,d,(N,h))>0$ et 
$C(n,a,d,(N,h))>0$ telles que si $(M,g)$ est une variété riemannienne 
de dimension $n$ vérifiant $\diam(M,g)\leq d$, $|K(M,g)|\leq a$ et 
si $\pi:(M,g)\rightarrow(N,h)$ est une fibration principale de fibre 
$T^k$ qui soit une $\varepsilon$-approximation de Hausdorff avec 
$\varepsilon<\varepsilon_0$, alors
$$\lambda_{1,1}(M,g)\geq C\cdot\Vol^2(M,g).$$
\end{theo}
\begin{remarque}
On obtient une minoration en fonction du volume de la variété et
pas en fonction du rayon d'injectivité, mais avec un exposant égal à
2, indépendamment de la dimension de la variété.
\end{remarque}

Ce résultat soulève deux questions qui restent ouvertes:
\begin{question}
Peut-on obtenir une minoration en fonction du rayon d'injectivité
avec un exposant indépendant de la dimension?
\end{question}
\begin{question} Peut-on généraliser ce résultat aux $p$-formes
différentielles, pour tout $p$ ?
\end{question}

 Nous commencerons dans le chapitre 4 par discuter de la topologie des 
fibrés principaux en tore, dans le but de déterminer un invariant 
topologique généralisant la classe d'Euler des fibrés en cercle et qui 
pourra être utilisé pour contrôler le spectre.

Dans le chapitre 5, nous étudierons comment, dans le cas d'un fibré 
principal en tore $T^k$ muni d'une métrique invariante, on peut se
ramener à l'étude des petites valeurs propres du laplacien à celle
du spectre du laplacien resteint aux formes différentielles
invariantes par l'action de $T^k$.
Dans le cas des formes différentielles d'un fibré
en cercle, Colbois et Courtois donnent, avec certaines hypothèses
sur la métrique, une minoration du spectre sur l'orthogonal des formes
invariantes~:
\begin{theo}[\cite{cc2}]
Soit $(N^n,h)$ une variété compacte, $S^1\hookrightarrow M^{n+1}\rightarrow
N$ un $S^1$-fibré principal, et $g_\varepsilon$ une métrique $S^1$-invariante
sur $M$ telle que $|K(M,g_\varepsilon)|<a$, $\diam(M,g_\varepsilon)<d$ et
que la longueur des fibres soit inférieure à $\varepsilon$. Il existe 
des constantes $C=C(n,a,d,(N,h))$ et $\rho=\rho(n)$ telles
que toute forme propre de valeur propre $\lambda\leq\frac C
{\varepsilon^{\frac 1{2\rho}}}$ est $S^1$-invariante.
\end{theo}
Nous allons démontrer un résultat semblable à celui de
Colbois et Courtois, mais sans utiliser d'hypothèse sur le diamètre et
la courbure du fibré, et en donnant une minoration plus précise du
spectre du laplacien restreint à l'orthogonal des formes invariantes:

\begin{theo}\label{inv:th1}
Soit $(M,g)$ un fibré principal en cercle muni d'une métrique
$S^1$-invariante. On note $l_0$ le maximum des longueurs des fibres.
 
Soit $\lambda$ une valeur propre de $\Delta^p$. Si $\lambda<
\displaystyle\left(\frac{2\pi}{l_0}\right)^2$, alors les formes
propres associées sont $S^1$-invariantes.
  
\end{theo}
\begin{remarque}
La constante $\displaystyle\left(\frac{2\pi}{l_0}\right)^2$ du théorème
est optimale, en ce sens qu'il existe des fibrés pour lequelles
elle est égale à une valeur propre associée à une forme propre non
invariante. En effet, $\left(\frac{2\pi}{l}\right)^2$ est la première
valeur propre du cercle de longueur $l$. Dans le cas d'un fibré trivial
$M=N\times S^1$ muni d'une métrique produit, les formes différentielles
de la formes $\alpha\wedge f$ où $\alpha$ est une $p$-forme harmonique
de $N$ et $f$ une fonction propre de $S^1$ seront des formes propres
de $\Delta^p$, de même valeur propre que $f$. Comme la fonction $f$ n'est
pas constante, les formes propres associés à la valeur propre
$\left(\frac{2\pi}{l}\right)^2$ ne sont pas toutes invariantes.
\end{remarque}
Comme $\left(\frac{2\pi}l\right)^2$ est le $\lambda_{0,1}$ du cercle
de longueur $l$, la constante $\left(\frac{2\pi}{l_0}\right)^2$ 
du théorème \ref{inv:th1}
peut s'interpréter comme la borne inférieure sur l'ensemble des fibres 
de la première valeur propre du laplacien agissant sur les fonctions
de la fibre. Dans le cas des fibrés en tore, on peut montrer 
un résultat semblable, faisant intervenir la première
valeur propre du laplacien $\Delta^0$ restreint au tore~:

\begin{theo}\label{inv:th2}
Soit $k\in\N^*$, $T^k\hookrightarrow M\stackrel{\pi}{\rightarrow}N$ un 
fibré en tore $T^k$, $\bar g$ une métrique invariante sur $T^k$ et $f$ 
une fonction sur $N$ strictement positive. 
Supposons que $M$ est muni d'une métrique 
$T^k$-invariante $g$ telle que pour tout $x\in N$, la restriction 
$\bar g_x$ de $g$ à la fibre $\pi^{-1}(x)$ vérifie 
$\bar g_x\leq f(x)\cdot\bar g$.

Soit $\lambda$ une valeur propre du laplacien agissant sur les formes
différentielles de $M$. Si $\lambda<\displaystyle(\sup_{x\in N}f(x))^{-1}
\cdot\lambda_{0,1}(T^k,\bar g)$, alors les formes propres 
associées sont $T^k$-invariantes.
\end{theo}

\begin{remarque}
Comme pour le théorème \ref{inv:th1}, la constante est optimale dans
le cas d'un fibré trivial muni d'une métrique produit.
\end{remarque}
 
 On peut déduire du théorème \ref{inv:th2} une inégalité en fonction
du maximum des diamètres des fibres. Cependant, on doit ajouter 
une hypothèse sur la métrique $g$. On verra en effet que la donnée d'une 
borne sur le diamètre des fibres ne permet pas de majorer la fonction $f$ 
du théorème \ref{inv:th2}.
\begin{cor}\label{inv:cor}
Soit $k\in\N^*$, $T^k\hookrightarrow M\stackrel{\pi}{\rightarrow}N$ un 
fibré en tore $T^k$, et $\bar g$ une métrique invariante sur $T^k$.
Supposons que $M$ est muni d'une métrique $T^k$-invariante 
telle que pour tout $x\in N$, la restriction de la métrique
à la fibre $\pi^{-1}(x)$ soit proportionnelle à $\bar g$.

Soit $\lambda$ une valeur propre du laplacien.
Si $\lambda<\displaystyle\left(\frac\pi {d_0}\right)^2$, où $d_0$ est
le maximum des diamètres des fibres pour la métrique induite par $g$,
alors les formes propres associées sont $T^k$-invariantes.
\end{cor}
\begin{remarque}\label{inv:rem}
La démonstration des deux théorèmes met en évidence le fait que si la
multiplicité d'une valeur propre est impaire, alors le sous-espace
propre associé contient des formes invariantes.
\end{remarque}

Le chapitre 6 aura pour but de montrer que pour obtenir le théorème
\ref{vol:th}, on peut se ramener à une situation géométrique plus simple.
En particulier, on montrera qu'une métrique de courbure et diamètre
bornés sur le fibré est proche d'une métrique invariante pour laquelle
les fibres sont totalement géodésiques :
\begin{theo}\label{geom:th:intro}
 Soient $a$ et $d$ deux réels strictement positifs, $n$ un entier et 
$(N,h)$ une variété riemannienne de dimension strictement inférieure 
à $n$.  Il existe des constantes $\varepsilon_0(n,a,d,(N,h))>0$, 
$\tau(n,a,d,(N,h))>0$ et $\tau'(n,a,d,(N,h))>0$ telles que si 
$(M,g)$ est une variété riemannienne de dimension $n$ vérifiant 
$|K(N,h)|\leq a$, $|K(M,g)|\leq a$, $\diam(M,g)\leq d$ et
si $\pi:(M,g)\rightarrow(N,h)$ une fibration principale de fibre $T^k$
qui soit une $\varepsilon$-approximation de Hausdorff avec $\varepsilon
<\varepsilon_0$, alors il existe des métriques $\tilde g$ et $\tilde h$  
sur $M$ et $N$ respectivement et une fibration principale 
$\pi': (M,\tilde g)\rightarrow(N,\tilde h)$ telles que 
\begin{enumerate}
\item L'action de $T^k$ sur $(M,\tilde g)$ est isométrique;
\item Les fibres de la fibration $\pi'$ sont totalement géodésiques;
\item $\displaystyle\frac 1\tau g\leq\tilde g\leq \tau g$ et
$\displaystyle\frac 1\tau h\leq\tilde h\leq \tau h$;
\item La restriction de $\tilde g$ à la fibre est telle que 
$\diam(\pi'^{-1}(x))\leq\tau'\varepsilon$, pour tout $x\in N$.
\end{enumerate}
\end{theo}
\begin{remarque}
On verra aussi que si l'on suppose que la métrique $g$ sur $M$ est 
$T^k$-invariante, alors on peut remplacer dans le théorème 
\ref{geom:th:intro} l'hypothèse sur la courbure de $(M,g)$ par
l'hypothèse plus faible $K(M,g)\geq -a$.
\end{remarque}

 Enfin, dans le chapitre 7, nous démontrerons le théorème \ref{vol:th} en
utilisant les résultats des chapitres 4 à 6, 
et nous discuterons de la possibilité d'exprimer la constante $C$
du théorème \ref{vol:th} en fonction d'invariants géométriques de $(N,h)$.

Une grande partie des chapitres \ref{cercle} et \ref{pvp:tore} a
été publiée sous une forme légèrement différente dans \cite{ja03}, et 
l'article \cite{ja04} contient ---~entre autres choses~--- la
démonstration du théorème \ref{vol:th}.

\chapter{Existence de petites valeurs propres du laplacien}
\pagenumbering{arabic}
\renewcommand\theequation{\thechapter.\arabic{equation}}
\setcounter{page}{1}

\section{Un exemple : la nilvariété d'Heisenberg de dimension 3}
\label{heisenberg}
Nous allons commencer par présenter un exemple de variété compacte 
admettant une petite valeur propre. Cet exemple nous sera utile par 
la suite car il illustre un certain nombre de phénomènes.

On considère le groupe d'Heisenberg de dimension 3
\begin{equation}\label{intro:df1}
G=\left\{\left(\begin{array}{ccc}1&x&z\\0&1&y\\0&0&1\end{array}\right)
,\ x,y,z\in \R\right\}.
\end{equation}
C'est un groupe nilpotent, difféomorphe à $\R^3$, et son centre
vérifie $Z(G)=[G,G]=\{\smat{10z\\010\\001},z\in\R\}$. On note
$\Gamma$ le sous-groupe de $G$ formé des matrices à coefficients
entiers. C'est un sous-groupe discret cocompact de $G$, et on
définit la nilvariété d'Heisenberg de dimension 3 par $N=\Gamma
\backslash G$.
 Sa topologie peut être décrite de trois manières :
\begin{itemize}
\item par construction c'est une nilvariété, c'est-à-dire le
quotient d'un groupe nilpotent par un sous-groupe cocompact;
\item c'est aussi un fibré en tore sur le cercle dont les fibres 
sont, dans la paramétrisation donnée par (\ref{intro:df1}), les
sous-variétés d'équation $x=c^{te}$. Ces sous-variétés de $G$ sont 
difféomorphes à $\R^2$, et leurs quotients dans $N$ sont des tores;
\item enfin, c'est un fibré principal en cercle sur le tore, les fibres 
étant définies comme les orbites de l'action du centre $Z(G)$ sur $N$.
\end{itemize}

 Soit $X$, $Y$ et $Z$ les champs de vecteurs invariants à gauche sur 
$N$ engendrés en $(0,0,0)$ par $\dep{}x$, $\dep{}y$ et $\dep{}z$.
Ces champs passent au quotient sur $N$, et vérifient
$[X,Y]=Z$ et $[X,Z]=[Y,Z]=0$. On va, à l'aide de ces champs, construire
une famille de métrique pour laquelle le diamètre et la courbure
seront uniformément bornés: soit $\alpha$, $\beta$ et $\gamma$ trois
réels positifs fixés. Pour tout $\varepsilon\in]0,1]$, on définit sur $N$
la métrique $g_\varepsilon$ comme étant la métrique invariante à
gauche telle qu'en tout point, la base $(X_\varepsilon,Y_\varepsilon,
Z_\varepsilon)$ définie par $X_\varepsilon=\varepsilon^{-\alpha}X$,
$Y_\varepsilon=\varepsilon^{-\beta}Y$ et $Z_\varepsilon=\varepsilon^{-
\gamma}Z$ soit orthonormée. Les crochets de Lie entre les champs de 
vecteurs de cette base sont 
\begin{equation}\label{intro:eq1}
[X_\varepsilon,Y_\varepsilon]=Z_\varepsilon^{\gamma-\alpha-\beta}
\et [X_\varepsilon,Z_\varepsilon]=[Y_\varepsilon,Z_\varepsilon]=0.
\end{equation}

 Comme les paramètres $\alpha$, $\beta$ et $\gamma$ sont positifs,
les normes de $X$, $Y$ et $Z$ pour la métrique $g_\varepsilon$ sont
inférieures à 1, donc le diamètre de $N$ reste borné quand $\varepsilon$
varie. D'autre part, comme les champs $X_\varepsilon$, $Y_\varepsilon$
et $Z_\varepsilon$ sont invariants, le tenseur de courbure peut s'écrire 
en fonction des coefficients des crochets de Lie entre ces champs. On impose
donc à $\alpha$, $\beta$ et $\gamma$ de vérifier $\tau=\gamma-\alpha-\beta
\geq0$, de sorte que la courbure reste elle aussi bornée.
 
 Comme la métrique sur $N$ est invariante à gauche, l'espace des $1$-formes
différentielles invariantes est stable par le laplacien. On va calculer 
son spectre en restriction à cet espace. On peut déduire de 
(\ref{intro:eq1}) que 
\begin{equation}
\de X^\flat_\varepsilon=\de Y^\flat_\varepsilon=0\et \de Z^\flat_\varepsilon
=-\varepsilon^\tau X^\flat_\varepsilon\wedge Y^\flat_\varepsilon,
\end{equation}
où $U^\flat$ désigne la $1$-forme duale de $U$ pour la métrique
$g_\varepsilon$, et donc que
\begin{equation}
\codiff(X^\flat_\varepsilon\wedge Y^\flat_\varepsilon)=-\varepsilon^\tau
Z^\flat_\varepsilon \et \codiff(X^\flat_\varepsilon\wedge Z^\flat_\varepsilon)
=\codiff(Y^\flat_\varepsilon\wedge Z^\flat_\varepsilon)=0.
\end{equation}
En restriction aux $1$-formes invariantes, l'opérateur $\de\codiff$
est nul. On a donc finalement 
\begin{equation}
\Delta X^\flat_\varepsilon=\Delta Y^\flat_\varepsilon=0\et
\Delta Z^\flat_\varepsilon=\varepsilon^{2\tau}Z^\flat_\varepsilon.
\end{equation}
La forme différentielle $Z^\flat_\varepsilon$ est donc une forme 
propre de valeur propre $\lambda=\varepsilon^{2\tau}$. On voit
que si $\tau>0$, par exemple si $\alpha=1$, $\beta=1$ et $\gamma=3$, 
cette valeur propre tend vers 0. La variété $N$ admet donc
une petite valeur propre. On peut noter qu'en revanche, si
$\tau=0$, par exemple si $\alpha=\beta=1$ et $\gamma=2$, 
la valeur propre $\lambda$ ne tend pas vers zéro. La question
\ref{intro:qu1} doit donc être formulée de manière plus précise:
\begin{question}\label{pvp:qu1}
 Comment varie le nombre de petites valeurs propres avec la géométrie
de l'effondrement~?
\end{question}

\section{Opérateur limite}\label{lott}
\subsection{Le cas des fonctions}
 Avant d'aborder la construction faite par J.~Lott d'un opérateur 
limite pour le laplacien de Hodge-de~Rham, rappelons le résultat
obtenu par K.~Fukaya dans le cas des fonctions. Dans \cite{fu},
Fukaya montre l'existence ---~sous certaines conditions~--- 
d'un opérateur limite pour le laplacien agissant sur les fonctions.
 
 Soit $\mathcal M(n,d)$ l'ensemble des
variétés riemanniennes $(M^n,g)$ de dimension $n$ telles que
leur courbure sectionnelle et leur diamètre vérifient respectivement 
$|K(M,g)|\leq1$ et $\diam(M,g)\leq d$, et $\lambda_k(M,g)$ la $k$-ième
valeur propre du laplacien agissant sur les fonctions de $M$,
les valeurs propres étant répétées s'il y a multiplicité. Le
problème est d'étendre continument pour tout $k$ la fonction
$(M,g)\rightarrow\lambda_k(M,g)$ à l'adhérence de $\mathcal M(n,d)$
dans l'ensemble des espaces métriques compacts. C'est impossible si
on munit cet ensemble de la topologie de Gromov-Hausdorff, mais
Fukaya considère l'ensemble des espaces métriques mesurés ---~\emph{i.e.}
muni d'une mesure de probabilité~--- muni de la topologie de Hausdorff
mesurée définie comme suit: soit $(X_i,\mu_i)$ une suite d'espaces
métriques mesurés. Elle converge vers $(X,\mu)$ s'il existe des
applications mesurables $\Psi_i:X_i\rightarrow X$ et une suite
strictement positive $\varepsilon_i$ vérifiant
\begin{enumerate}
\item $\displaystyle\lim_{i\rightarrow\infty}\varepsilon_i=0$~;
\item L'$\varepsilon_i$-voisinage de $\Psi_i(X_i)$ est $X$~;
\item Pour tout $p,q$ dans $X_i$, on a
$$|d(\Psi_i(p),\Psi_i(q))-d(p,q)|<\varepsilon_i\ ;$$
\item Pour toute fonction $f$ sur $X$, on a
$$\lim_{i\rightarrow\infty}\int f\circ\Psi_i\de\mu_i=\int f\de\mu,$$
c'est-à-dire que la mesure $(\Psi_i)_*(\mu_i)$ converge vers $\mu$ pour la
topologie faible-$*$.
\end{enumerate}
Il obtient alors:
\begin{theo}\label{fu:1:1}
Soit $\overline{\mathcal M(n,d)}$ l'adhérence de $\mathcal M(n,d)$
dans l'ensemble des espaces métriques mesurés, en munissant chaque
variété $(M^n,g)$ de sa mesure riemannienne normalisée.
\begin{enumerate}
\item La fonction $\lambda_k$ s'étend continuement à l'ensemble
$\overline{\mathcal M(n,d)}\backslash\{(point,1)\}$.
\item Pour tout $(X,\mu)\in\overline{\mathcal M(n,d)}$, il existe un
opérateur auto-adjoint $P_{(X,\mu)}$ agissant sur $L^2(X,\mu)$ tel que 
$\lambda_k(X,\mu)$ soit égal à la $k$-ième valeur propre de $P_{(X,\mu)}$.
\item Supposons que $\lim_{i\rightarrow\infty}(M_i,g_i)=(X,\mu)$. Soit
$\varphi_{k,i}$ une fonction propre normalisée associée à la valeur
propre $\lambda_k(M_i,g_i)$ et $\Lambda_{k,i}=\{\varphi\circ\Psi_i,
\varphi\in L^2(X,\mu), P_{(X,\mu)}\varphi=\lambda_k(X,\mu)\varphi\}$.
Alors $\lim_{i\rightarrow\infty}d(\Lambda_{k,i},\varphi_{k,i})=0$.
\end{enumerate}
\end{theo}
\begin{remarque} Dans le cas où l'espace limite $(X,\mu)$ est
une variété riemannienne, la mesure limite $\mu$ et l'opérateur limite
$P_{(X,\mu)}$ ne sont pas nécessairement la mesure riemanienne ni le laplacien
de la variété, comme le montre l'exemple suivant :
\end{remarque}

\begin{exemple}
On considère le tore $T^2=\{(s,t),\ s,t\in S^1\}$, $c:S^1\rightarrow
\R$ une fonction positive $C^\infty$, et la famille de métriques 
$g_\varepsilon(c)$ définie sur $T^2$ par
$$g_\varepsilon(c)=\de s^2\oplus\varepsilon^2c(s)^2\de t^2.$$
Si $f\in C^\infty(T^2)$, on a alors
$$\Delta_{(T^2,g_\varepsilon(c))}f(s,t)=-c(s)^{-1}\Dep{}s\left(
c(s)\Dep{}sf(s,t)\right)-
\varepsilon^2c(s)^{-2}c(s)^{-2}\frac{\partial^2}{\partial t^2}f(s,t).$$
Si on note $\lambda_k(c)$ la $k$-ième valeur propre de
l'opérateur $P_c$ défini sur $S^1$ par 
$P_c(f)(s)=-c(s)^{-1}\frac\de{\de s}\left(c(s)\frac\de{\de s}f(s)\right)$,
alors
$$\lim_{\varepsilon\rightarrow0}\lambda_k(T^2,g_\varepsilon(c))=
\lambda_k(c),$$
et
$$\lim_{\varepsilon\rightarrow0}(T^2,g_\varepsilon(c))=(S^1,\mu)$$
avec $\mu=c\cdot\de s$.
\end{exemple}

\begin{remarque}\label{fu:1:2} Les fonctions contenues dans $\Lambda_{k,i}$
ont la propriété d'être constantes sur chacun des $\Psi_i^{-1}(x)$, $x\in X$.
Cela signifie, dans le cas où les $\Psi_i$ définissent des fibrations,
que les fonctions propres $\varphi_{k,i}$ sont
approximées par des fonctions constantes sur les fibres.
\end{remarque}

\subsection{Construction de l'opérateur}
Récemment, J.~Lott a généralisé le résultat de Fukaya aux formes
différentielles (\cite{lo}, \cite{lo2}). Nous allons ici présenter 
la construction de l'opérateur limite en nous restreignant par
soucis de clarté au cas
d'un fibré $F\hookrightarrow (M,g)\rightarrow (N,h)$ sur une variété 
riemannienne dont la fibre est une nilvariété $F=\Gamma\backslash G$,
et qui tend pour la distance de Gromov-Hausdorff vers sa base. Nous
noterons $\nabla^{aff}$ la connexion sur $F$ telle que les champs
invariants à gauche soient parallèles.

 La première difficulté est de déterminer un espace sur lequel
va agir l'opérateur limite. Dans la construction de Fukaya, l'opérateur
$P_{(X,\mu)}$ du théorème \ref{fu:1:1} agit sur $L^2(N)$, c'est-à-dire
sur un espace de fonctions à valeurs réelles sur la base. Dans le cas 
des formes différentielles, on utilise un espace de formes
différentielle sur la base, mais à valeur dans un espace plus
grand que $\R$. Plus précisement, on considère un fibré vectoriel
gradué $E=\bigoplus_{j=0}^mE^j$ sur la base dont chaque fibre est
munie d'un produit scalaire gradué ---~\emph{i.e.} tel que les $E^i$
soient orthogonaux entre eux~--- noté $h_E$, et on munit ce fibré d'une
superconnexion $A'$ de degré 1, c'est-à-dire d'un opérateur de la
forme $$A'=A_{[0]}'+A_{[1]}'+A_{[2]}'$$ où
\begin{itemize}
\item $A_{[0]}'\in C^\infty(N,\Hom(E^*,E^{*+1}))$~;
\item $A_{[1]}'$ est une connexion sur $E$ qui préserve la graduation~;
\item $A_{[2]}'\in\Omega^2(N,\Hom(E^*,E^{*-1}))$.
\end{itemize}
On peut étendre cette superconnexion par la règle de Leibniz à un
opérateur sur l'espace $\Omega(N,E)$ des formes différentielles
sur $N$ à valeur dans $E$. La métrique $h_N$ et le produit scalaire
$h_E$ permettent de construire un produit scalaire sur $\Omega(N,E)$,
et donc de définir un opérateur adjoint à $A'$ noté $(A')^*$,
et un laplacien sur $\Omega(N,E)$ par $\Delta_E=A'(A')^*+(A')^*A'$.
On notera $\Delta^p_E$ la restriction de $\Delta_E$ à
$\bigoplus_{a+b=p}\Omega^a(N,E^b)$.

Pour construire le fibré sur lequel agit l'opérateur limite, J.~Lott
se ramène d'abord en utilisant \cite{cfg} au cas où $(M,g)$ est un
fibré affine riemannien:
\begin{df}\label{lott:df} Un fibré riemannien $F\hookrightarrow
(M,g)\rightarrow(N,h)$ est une fibré affine riemannien si:
\begin{itemize}
\item le groupe de structure du fibré est contenu dans le groupe
$\Aff(F)$ des difféomorphismes de $F$ qui préservent $\nabla^{aff}$~;
\item $M$ est muni d'une distribution horizontale $T^HM$ dont
l'holonomie est dans $\Aff(F)$~;
\item chaque fibre est munie d'une métrique $g_{F_b}$ qui est parallèle
par rapport à la connexion affine $\nabla^{aff}$ sur la fibre~;
\item $N$ est muni d'une métrique $h_N$~;
\item la métrique $g$ sur $M$  s'écrit $g=h_N\oplus g_{F_b}$ relativement
à la distribution $T^HM$.
\end{itemize}
\end{df}
 En décomposant les formes différentielles de $M$ en leurs parties
verticale et horizontale, on peut écrire $\Omega^*(M)\simeq\Omega^*(N,W)$,
où $W$ est un fibré sur $N$ dont la fibre est isomorphe à $\Omega^*(F)$.
Dans le cas d'un fibré affine, le groupe de structure du fibré $M$ 
préserve $\nabla^{aff}$, et par conséquent l'action de ce groupe sur
$\Omega^*(F)$ préserve le sous-espace des formes invariantes. 
On peut donc définir le sous-fibré $E$ de $W$ de fibre 
$\Lambda^*(\mathfrak n^*)$, ainsi que l'espace $\Omega^*(N,E)$ sur lequel
va agir l'opérateur limite. On a de plus un plongement 
$\Omega^*(N,E)\hookrightarrow\Omega^*(M)$ qui permet d'identifier
chaque élément de $\Omega^*(N,E)$ à une forme différentielle sur $M$
parallèle pour la connexion $\nabla^{aff}$.

 Selon les résultats donnés dans \cite{bl} sur les superconnexions, 
et en utilisant le fait que sur un tel fibré, 
l'espace des formes parallèles le long de la fibre est stable par 
l'action de la différentielle extérieure, cette différentielle  induit 
par l'intermédiaire du plongement
$\Omega(N,E)\hookrightarrow\Omega(M)$ une superconnexion sur $E$ 
telle que $A_{[0]}'$ soit la différentielle sur $\Lambda^*(\mathfrak n^*)$,
$A_{[1]}'$ soit la connexion sur le fibré $E$ induite par $T^HM$ et 
$A_{[2]}'$ soit le produit intérieur $i_T$ par la forme de courbure $T$
de la distribution $T^HM$. D'autre part,
la métrique $g_{F_b}$ induit un produit scalaire $h_E$ sur les
fibres du fibré vectoriel $E$. 

La structure  de fibré affine riemannien
induit donc à la fois une superconnexion et un produit scalaire sur $E$,
et permet par conséquent de définir un laplacien $\Delta_E$ sur
$\Omega(N,E)$. En notant $\sigma(\Delta_E^p)$ le spectre du laplacien
$\Delta_E^p$, J.~Lott montre que dans cette situation, les petites 
valeurs propres du laplacien $\Delta_M$ 
sur $M$ sont celle de $\Delta_E^p$ (\cite{lo}, théorème 1): 
\begin{theo}\label{lott:th1}
Si $M$ est un fibré affine riemannien sur $N$ et si on note
$R_M$ et $R_F$ les tenseurs de courbures respectivement
sur $M$ et sur $F_b$ pour la métrique $g_{F_b}$, $\diam(F)$ la 
borne supérieure des diamètres des fibres et $\Pi$ la seconde forme
fondamentale des fibres, alors il existe des constantes 
$a$, $a'$ et $c$ qui ne dépendent que de 
$\dimension(M)$ telles que si $\|R_F\|_\infty\diam(F)^2\leq a'$ alors pour
tout $p\leq\dimension(M)$,
\begin{eqnarray*}
\sigma(\Delta_M^p)\cap[0,a\cdot\diam(F)^{-2}-c(\|R_M\|_\infty+\|\Pi\|^2_\infty
+\|T\|^2_\infty)[ & = &\\
\sigma(\Delta_E^p)\cap[0,a\cdot\diam(F)^{-2}-c(\|R_M\|_\infty+\|\Pi\|^2_\infty
+\|T\|^2_\infty)[. &  &
\end{eqnarray*}
\end{theo}

On est donc ramené à la recherche d'un opérateur limite sur
$\Omega(N,E)$. On ne peut considérer séparément la limite de la 
superconnexion $A'$ et de la métrique $h_E$. En effet, $A'$ ne dépend pas
de la métrique sur $M$ et donc elle est constante au cours de l'effondrement,
et $h_E$ dégénère donc sa limite ne permet pas de définir un opérateur
adjoint comme $(A')^*$. L'idée de Lott est de considérer l'ensemble
$(\mathcal S_E\times\mathcal H_E)$, où $\mathcal S_E$ est l'espace
des superconnexions de degré 1 sur $E$ et $\mathcal H_E$ l'espace des
produits scalaires euclidiens sur $E$, et de quotienter cet ensemble par
le groupe $\mathcal G_E$ des $\GL(E)$-transformations de jauge sur $E$ 
qui préservent la graduation. Il obtient alors le résultat 
de compacité suivant (\cite{lo}, théorème 3):

\begin{theo} Soit $(N,h_N)$ une variété riemannienne fixé, $n>
\dimension(N)$, $a>0$ et $\varepsilon>0$. Il existe une partie compacte 
$K(n,a,\varepsilon)\subset(\mathcal S_E\times\mathcal H_E)/\mathcal G_E$ 
telle que si $(M^n,g)$ est une variété riemannienne de dimension $n$ 
telle que $\|R_M\|_\infty\leq a$ 
et $d_{GH}(M,N)\leq\varepsilon$, alors la classe d'équivalence du couple
$(A',h_E)$ induit par $g$  appartient à $K$.
\end{theo}

 Étant donnée une suite de métriques $(g_i)_{i\in\N}$ qui effondre $M$ sur sa
base, on peut extraire de la suite $[(A'_i,h_{E,i})]_{i\in\N}$ d'éléments de
$(\mathcal S_E\times\mathcal H_E)/\mathcal G_E$ une sous-suite qui
converge vers une superconnexion limite qui permet de construire un
laplacien limite $\Delta_\infty$, dont le spectre est la limite
du spectre du laplacien sur $M$.

\begin{remarque}
La suite $[(A'_i,h_{E,i})]$ ne converge pas nécessairement. On peut par
exemple au paragraphe \ref{heisenberg} prendre $\alpha=0$, $\beta=1$ et 
choisir pour
$\gamma$ une fonction qui oscille entre $1$ et $2$. La première valeur
propre va osciller entre $1$ et $t$ sans converger. L'opérateur
$\Delta_E$ ne converge donc pas, ni la classe de $(A'_t,h_{E,t})$.
\end{remarque}

\begin{remarque}
Pour l'étude du spectre sur les fonctions, on peut se resteindre
au fibré $E^0$, qui est un fibré trivial en droite réelle sur $N$. 
Cependant, le produit scalaire $h_{E^0}$ n'est pas trivial. 
Il correspond à la mesure sur l'espace limite dans le travail de
Fukaya.
\end{remarque}

\subsection{Petites valeurs propres}
Pour déterminer la dimension du noyau de $\Delta_\infty$, on peut calculer
la cohomologie $H^*(A')$ pour l'action sur $\Omega(N,E)$ de la 
superconnexion limite $A'$. J.~Lott calcule une majoration de cette
dimension en utilisant la théorie des suites spectrales, et en
remarquant le fait suivant: le terme $A_{[0]}'$ de la superconnexion 
vérifie $(A_{[0]}')^2=0$, et définit donc un complexe 
différentiel sur les fibre de $E$, dont la cohomologie $H^*(A'_{[0]})$ 
est un fibré vectoriel gradué sur $N$. De plus, la connexion $A_{[1]}'$
sur le fibré $E$ passe au quotient sur $H^*(A'_{[0]})$ en une connexion
plate, c'est-à-dire telle que $(A_{[1]}')^2=0$, et définit donc aussi
un complexe différentiel dont la cohomologie est $H^*(N,H^*(A'_{[0]}))$.
Les premiers termes de la suite de Leray $(\mathcal E^{*,*}_r,d_r)$ sont
\begin{eqnarray*}
\mathcal E^{*,*}_0 &=& \Omega^*(N,E^*),\\
\mathcal E^{*,*}_1 &=& \Omega^*(N,H^*(A'_{[0]}))\et\\
\mathcal E^{*,*}_2 &=& H^*(N,H^*(A'_{[0]})),
\end{eqnarray*}
avec $d_0=A'_{[0]}$ et $d_1=A_{[1]}'$.
Lott en déduit:
\begin{pr}
$$\dimension\Ker\Delta_\infty^p\leq\sum_{a+b=p}\dimension
\left(H^a(N,H^b(A'_{[0]}))\right).$$
\end{pr}
 Cette formule peut se simplifier dans certains cas, en particulier
pour les $1$-formes:
\begin{cor}\label{lott:cor1}
$$\dimension\Ker\Delta_\infty^1\leq b_1(N)+\dimension(F).$$
\end{cor}
 On peut noter que cette majoration est très générale, en ce sens qu'on
ne fait pas d'hypothèse sur la structure du fibré, ni sur la géométrie
de l'effondrement. Cependant, le théorème énoncé en~\ref{cercle:th1} dans
l'introduction montre que la topologie impose des restrictions sur le 
nombre de petites valeurs propres possible. En particulier, dans la 
situation du théorème~\ref{cercle:th1}, le cas d'égalité de la majoration 
donnée par le corollaire~\ref{lott:cor1} n'est atteint que si $G$ est 
nilpotent (cf. remarque~\ref{cercle:rq1}).

Dans le cas d'un fibré en cercle, on obtient aussi une expression simple:
\begin{cor}
Si $M$ est un fibré en cercle sur $N$, alors 
$$\dimension\Delta_\infty^p\leq b_p(N)+b_{p-1}(N)$$.
\end{cor}
Cependant, on sait déjà (\cite{cc2}) qu'il y a 
nécessairement égalité dans l'inégalité ci-dessus.

 D'autre part, Lott montre qu'une petite valeur propre
ne peut être obtenue que selon trois mécanismes (\cite{lo}, Th. 5):
\begin{theo}\label{lott:th2}
Soit $g_i$ une suite de métriques qui effondre $M$ sur $N$. Supposons
que $\lim_{i\rightarrow\infty}\lambda_{1,p}(M,g_i)=0$. Alors au moins
l'une des trois conditions suivantes est vérifiée:
\begin{enumerate}
\item Il existe $q\in[0,p]$ tel que $b_q(F)<\dimension\Lambda^q
(\mathfrak n^*))$~;
\item Il existe $q\in[0,p]$ tel que l'holonomie du fibré de base $N$
et de fibre $H^q(F)$ n'est pas semi-simple~;
\item La suite spectrale de Leray qui calcule la cohomologie
$H^*(M,\R)$ du fibré $M$ ne dégénère pas au rang $2$.
\end{enumerate}
\end{theo}

 Les trois situations qui interviennent dans ce théorème sont illustrées
par les trois descriptions topologiques de la nilvariété d'Heisenberg
qu'on a données en \ref{heisenberg}

 Dans le premier cas, le terme $A_{[0],i}'$ de la superconnexion dégénère 
quand $i$ tend vers l'infini ---~c'est-à-dire que si on note 
$(A_{[0],\infty}',h_E)$ la limite de $(A_{[0],i}',h_{E,i})$ dans
$(\mathcal S_E\times\mathcal H_E)/\mathcal G_E$, la dimension
du noyau de $A_{[0],\infty}'$ est plus grande que celle du noyau de 
$A_{[0],i}'$~--- donc le laplacien
restreint à la fibre admet une petite valeur propre. Géométriquement,
cela signifie que la nilvariété $F$ n'est pas un tore. Un exemple 
simple est donné par $M=N\times F$ muni d'une métrique produit, où 
$F$ est la nilvariété d'Heisenberg de dimension $3$, et de considérer
sur $F$ la suite de métriques définie en \ref{heisenberg} avec 
$(\alpha,\beta,\gamma)=(1,1,3)$.

 La deuxième condition signifie que le terme $A_{[1],i}'$ de la
superconnexion dégénère quand $i$ tend vers l'infini. L'exemple
le plus simple est donné par la variété d'Heisenberg de dimension $3$
vue comme fibré en tore sur le cercle : comme la fibre est plate,
$A_{[0],i}'=0$, et comme la base est de dimension $1$, les $2$-formes
sur la base, et donc $A_{[2],i}'$, sont nulles. Cependant, si on
prend $\alpha=0$, $\beta=1$ et $\gamma=2$ en \ref{heisenberg}, on
a bien une petite valeur propre.

 La troisième condition est illustrée par les situations de fibré
principal. Si on considère la variété d'Heisenberg comme un fibré
principal en cercle sur le tore $T^2$, on a $A_{[0],i}'=0$ (la
fibre est plate) et l'holonomie de fibré $E$ est nécessairement
semi-simple car les fibres des fibrés $(E_i)_{i=0,1}$ sont de dimension $1$.
On verra d'autres exemples de fibrés principaux dans le chapitre 
\ref{pvp:tore}.

 Le théorème \ref{lott:th2} donne une condition nécessaire à l'existence
de petites valeurs propres, ce qui répond à la question \ref{intro:qu1},
mais pas à la question \ref{pvp:qu1}.

Dans le cas particulier d'une variété $M$, s'effondrant sur un cercle,
Lott donne le corollaire suivant (\cite{lo}, Cor. 4):
\begin{cor}\label{lott:cor2}
Soit $g_i$ une suite de métriques qui effondre $M$ sur $S^1$. Supposons
que $\lim_{i\rightarrow\infty}\lambda_{1,p}(M,g_i)=0$. Alors au moins
l'une des deux conditions suivantes est vérifiée:
\begin{enumerate}
\item Il existe $q\in[0,p]$ tel que $b_q(F)<\dimension\Lambda^q
(\mathfrak n^*))$~;
\item Il existe $q\in[0,p]$ tel que si on note $\Phi^*\in\Aut(H^*(Z))$ l'action
de l'holonomie sur le fibré $H^*(Z)$, alors la réduite de Jordan de 
$\Phi^q$ ou $\Phi^{q-1}$ contient un bloc unipotent non trivial.
\end{enumerate}
\end{cor}
Le chapitre suivant sera consacré à une étude des situations de fibrés
en tore sur le cercle, qui illustrent le point 2 du théorème \ref{lott:th2}.

\chapter{Effondrements homogènes de fibrés en tores sur le cercle}
\label{cercle}
\section{Structure homog\`ene}\label{cercle:st}
Nous commen\c{c}ons par d\'emontrer le d\'ebut du th\'eor\`eme \ref{cercle:th1}
en construisant le groupe $G$ et le r\'eseau $\Gamma$ qui nous
int\'eressent.
Consid\'erons un fibr\'e $M$ en tore $T^n$ sur le cercle qui est
la suspension d'un diff\'eomorphisme lin\'eaire $\varphi$ repr\'esent\'e par 
la matrice $A\in\SL_n(\Z)$. Un tel fibr\'e sera hom\'eomorphe \`a 
\begin{equation}\label{cercle:st:df}
M:=T^n\times [0,1]_{/(x,0)\sim(\varphi(x),1)},
\end{equation}
Pour construire $G$, on va munir $\R^{n+1}$ d'une structure de groupe
telle que $\Z^{n+1}\backslash\R^{n+1}=M$. Si on note $(x_1,\cdots,x_n,y)$ 
les \'el\'ements de $\R^{n+1}$, une telle structure devra v\'erifier
\begin{equation}\label{cercle:st:eq1}
(k_1,\cdots,k_n,0)\cdot(x_1,\cdots,x_n,y)=(x_1+k_1,\cdots,x_n+k_n,y)
\end{equation}
de sorte que les sous-espaces de $\R^{n+1}$ d'\'equation $y=c^{te}$
passent au quotient comme des tores $T^n$, et
\begin{equation}\label{cercle:st:eq2}
(0,\cdots,0,l)\cdot(x_1,\cdots,x_n,y)=(A^l\left(\vecdots{x_1}{x_n}\right),y+l)
\end{equation}
de sorte que la structure de fibr\'e soit bien celle d\'efinie par 
(\ref{cercle:st:df}). Cette structure est effectivement r\'ealis\'ee en 
d\'efinissant $G$ comme l'image du plongement
\begin{equation}\label{cercle:st:plongement}
(x_1,\cdots,x_n,y)\longmapsto\left(\begin{array}{c|cc}
A^y&\vecdots00&\vecdots{x_1}{x_n}\\\hline
0&\begin{smallmatrix}\\1\\0\end{smallmatrix}
&\begin{smallmatrix}\\y\\1\end{smallmatrix}
\end{array}\right) .
\end{equation}
Comme on se restreint aux matrices $A$ qui admettent un logarithme $B$,
l'expression $A^y$ est bien d\'efinie en posant $A^y=\exp(yB)$.
On peut facilement v\'erifier
que cette application est injective, que son image $G$ est bien un 
sous-groupe de $\GL_{n+2}(\R)$ et que sa structure est
bien celle d\'efinie par (\ref{cercle:st:eq1}) et  
(\ref{cercle:st:eq2}). Enfin,
l'image de $\Z^{n+1}$ par cette application est bien un sous-groupe
discret de $G$, qu'on notera $\Gamma$. La vari\'et\'e $M$ est 
donc hom\'eomorphe au quotient $\Gamma\backslash G$.

Remarque: on peut v\'erifier que si $A=\smat{1&1\\0&1}$, le groupe
$G$ obtenu est isomorphe au groupe d'Heisengerg de dimension $3$ tel
qu'il est pr\'esent\'e dans l'exemple du paragraphe \ref{heisenberg}.

\section{Laplacien}\label{cercle:laplacien}
Soit $X_i$ et $Y$ les champs invariants \`a gauche engendr\'es en
$I_{n+2}$ respectivement par 
\begin{equation}
\Dep{}{x_i}=\left(\begin{array}{c|c}
0&\begin{smallmatrix}
\vecdots0\ccdot&\vecdots00 \\
\ccdot&1\\
\vecdots{\vphantom{0}\ccdot}0&\vecdots00
\end{smallmatrix}\\\hline 0&0\end{array}\right) \et
\Dep{}{y}=\left(\begin{array}{c|c}
B&\begin{smallmatrix}\vecdots00&
\vecdots00\end{smallmatrix}\\\hline
0&\begin{smallmatrix}0&\:\:1\\
0&\:\:0\end{smallmatrix}\end{array}\right).
\end{equation}
Ces champs v\'erifient $[X_i,X_j]=0$ et $[Y,X_i]=\sum_{j=1}^nb_{ji}X_j$. 
On peut remarquer que l'application $X\mapsto [Y,X]$ est un endomorphisme
de  l'espace $\Gamma(T_VM)^{{}^G}$ des champs de vecteurs invariants 
verticaux, c'est-\`a-dire l'espace engendr\'e par les $X_i$, et dont la 
matrice est $B$. On notera $f$ cet endomorphisme.

On fixe une m\'etrique homog\`ene $g$ 
sur $M$ en se donnant une base $(V_i)_{i\in[1,n]}$ de l'espace 
$\Gamma(T_VM)^{{}^G}$, 
cette m\'etrique \'etant telle que $(V_1,\cdots,V_n,Y)$ soit orthonorm\'ee en
tout point. On notera $(V^\flat_1,\cdots,V^\flat_n,Y^\flat)$ sa base duale, 
et $C$ la matrice de $f$ dans la base $(V_1,\cdots,V_n)$.
On va d\'eterminer le spectre du laplacien $\Delta^1_{inv}$ restreint 
\`a  l'ensemble $\Omega^1(M)^G$ des $1$-formes invariantes \`a gauche 
en fonction des coefficients de $C$. Plus pr\'ecis\'ement, on a le
\begin{lem}
La matrice du laplacien $\Delta^1_{inv}$ dans la base 
$(V^\flat_1,\cdots,V^\flat_n,Y^\flat)$ est
$$\Delta^1_{inv}:\left(\begin{array}{cc}C\trans C&
\vecdots{\displaystyle0}{\displaystyle0}\\
0\cdots0&0\end{array}\right).$$
\end{lem}
\textbf{D\'emonstration:}
Les crochets de Lie entre les vecteurs de la nouvelle base sont
\begin{equation}
[V_i,V_j]=0 \et [Y,V_i]=\sum_{j=1}^nc_{ji}V_j.
\end{equation}
Soit $\alpha$ une $1$-forme diff\'erentielle invariante. Sa diff\'erentielle
ext\'erieure est déterminée par la relation $\de\alpha(U_1,U_2)=U_1\cdot
\alpha(U_2)-U_2\cdot\alpha(U_1)-\alpha([U_1,U_2])$, o\`u $U_1$ et $U_2$
sont des champs de vecteur. Si ces champs sont invariants \`a gauche, cette
relation devient: $\de\alpha(U_1,U_2)=-\alpha([U_1,U_2])$.
On en d\'eduit:
\begin{equation}\label{cercle:laplacien:eq1}
\de Y^\flat=0\et \de V_i^\flat=-\sum^n_{j=1} c_{ij}Y^\flat\wedge V_j^\flat.
\end{equation}
La matrice de la diff\'erentielle ext\'erieure 
$\de :\Omega^1(M)^G\rightarrow\Omega^2(M)^G$ sera, dans les bases
$(V^\flat_1,\cdots,V^\flat_n,Y^\flat)$ et $(Y^\flat\wedge V_i^\flat,V_i^\flat
\wedge 
V_j^\flat)$,
\begin{equation}
\de:\smat{\displaystyle-\trans C &\vecdots00\\\displaystyle 0&\vecdots00}.
\end{equation}
Les deux bases sont orthonorm\'ees, donc
la matrice dans ces bases de la divergence 
$\codiff:\Omega^2(M)^G\rightarrow\Omega^1(M)^G$ sera donc la transpos\'ee 
de la matrice ci-dessus.

Comme la diff\'erentielle restreinte \`a $\Omega^0(M)^G$ est nulle,
le laplacien $\Delta=\codiff\de+\de\codiff$ se r\'eduit sur $\Omega^1(M)^G$
\`a l'op\'erateur $\codiff\de$. On en d\'eduit la matrice du laplacien
$\Delta^{inv}$ restreint \`a $\Omega^1(M)^G$ est, dans la base 
$(V^\flat_1,\cdots,V^\flat_n,Y^\flat)$,

\begin{equation}
\left(\begin{array}{cc}
\displaystyle-C&\displaystyle 0\\0\cdots0&0\cdots0\end{array}\right)
\smat{\displaystyle-\trans C &\vecdots00\\\displaystyle 0&\vecdots00}
=\left(\begin{array}{cc}C\trans C&
\vecdots{\displaystyle0}{\displaystyle0}\\
0\cdots0&0\end{array}\right).
\end{equation}
\carrenoir

Remarque: On a fait ici le calcul pour un $Y$ fix\'e, c'est-\`a-dire
pour un certain choix de connexion du fibr\'e. Mais si on choisit $Y'$
tel que $Y'-Y\in\Gamma(T_VM)^{{}^G}$ et une m\'etrique telle que
$(V_1,\cdots,V_n,Y')$ soit orthonorm\'ee, le r\'esultat sera le m\^eme
car on aura toujours $[Y',V_i]=[Y,V_i]=\sum_{j=1}^nb_{ji}V_j$.

On peut noter que la métrique intervient par la réécriture de la 
matrice $B$ dans une base orthonormée. Ce travail de renormalisation
correspond dans le travail de J.~Lott au passage au quotient de 
$(\mathcal S_E\times\mathcal H_E)$ par le groupe de transformation
de jauge $\mathcal G_E$. Si deux métriques donnent la même matrice
$C$, cela signifie que les deux éléments correspondants dans 
$(\mathcal S_E\times\mathcal H_E)$ appartiennent à la même classe
dans $(\mathcal S_E\times\mathcal H_E)/\mathcal G_E$.

\section{Courbure}\label{cercle:courbure}
Nous allons d\'emontrer dans cette partie un lemme qui nous servira \`a
faire le lien entre le contr\^ole de la courbure et l'existence de petites
valeurs propres.
\begin{lem}\label{cercle:coubure:lem}
Soit $a$ la borne sup\'erieure de la valeur absolue de la courbure 
sectionnelle de $(M,g)$. Il existe des constantes $\tau(n)>0$ et 
$\kappa(B)$ telle que
$$\tau^{-1}a<\Tr(C\trans C)<\tau a+\kappa.$$
\end{lem}
\textbf{D\'emonstration:}
Rappelons tout d'abord l'expression suivante (dont
le lecteur pourra trouver la d\'emonstration dans \cite{ce})
de la courbure sectionnelle $K(U,V)$, o\`u $U$ et $V$ sont deux champs
invariants \`a gauche d'un groupe de Lie quelconque:
\begin{eqnarray}\label{cercle:courbure:eq}
K(U,V)&=&\frac14\|\ad_U^*V+\ad_V^*U\|^2-\langle\ad_U^*U,\ad_V^*V\rangle\\
&&-\frac34\left\|[U,V]\right\|^2-\frac12\langle[[U,V],V],U\rangle
-\frac12\langle[[V,U],U],V\rangle.
\nonumber\end{eqnarray}
Nous allons appliquer cette relation aux champs de la base $(V_i,Y)$. Pour
cela, remarquons d'abord que les matrices de $\ad_Y$ et $\ad_{V_i}$ 
sont, dans cette base
\begin{equation}
\ad_Y:\smat{\displaystyle C&\vecdots00\\0\cdots0&0}
\et \ad_{V_i}:\smat{\displaystyle 0&\vecdots{-c_{1i}}{-c_{ni}}\\0\cdots0&0}
\end{equation}
On en d\'eduit $\ad_{V_i}^*Y=0$, $\ad_Y^*Y=0$, 
$\ad_Y^*V_i=\sum_jc_{ij}V_j$ et $\ad_{V_i}^*V_j=-c_{ji}Y$,
et donc que
\begin{eqnarray}\label{cercle:courbure:eq1}
K(Y,V_i)&=&\frac14\|\ad_Y^*V_i\|^2-\frac34\left\|[Y,V_i]\right\|^2
-\frac12\langle[[V_i,Y],Y],V_i\rangle\nonumber\\
&=&\frac14\sum_j\left(c_{ij}^2-3c_{ji}^2-2c_{ij}c_{ji}\right)\nonumber\\
&=&-\sum_jc_{ji}^2+\frac14\sum_j(c_{ij}-c_{ji})^2
\end{eqnarray}
et
\begin{eqnarray}\label{cercle:courbure:eq2}
K(V_i,V_j)&=&\frac14\|\ad_{V_i}^*V_j+\ad_{V_j}^*V_i\|^2
-\langle\ad_{V_i}^*V_i,\ad_{V_J}^*V_j\rangle\nonumber\\
&=&\frac14(c_{ij}+c_{ji})^2-c_{ii}c_{jj}.
\end{eqnarray}
D'autre part, comme $C$ est la matrice de $f$, le terme de degr\'e $n-2$ du
polyn\^ome caract\'eristique est ind\'ependant de la m\'etrique choisie. Le
calcul montre que son coefficient est $\kappa=\sum_{ij}(c_{ii}c_{jj}-c_{ij}
c_{ji})$. On peut en d\'eduire que
\begin{equation}
\sum_{i,j=1}^nK(V_i,V_j)+\kappa=\sum_{i,j=1}^n\left(
\frac14(c_{ij}+c_{ji})^2-c_{ij}c_{ji}\right)=\frac14\sum_{i,j=1}^n
(c_{ij}-c_{ji})^2,
\end{equation}
et donc que
\begin{equation}
\sum_{i,j=1}^nc_{ji}^2=\sum_{i,j=1}^nK(V_i,V_j)-\sum_{i=1}^nK(Y,V_i)+\kappa
\leq(n^2+n)a+\kappa,
\end{equation}
ce qui montre l'une des deux in\'egalit\'es du lemme. La seconde d\'ecoule
du fait que la courbure sectionnelle s'\'ecrit comme un polyn\^ome homog\`ene
de degr\'e deux relativement aux $c_{ij}$.
\hfill\carrenoir\nolinebreak\vspace*{1mm}

\section{Petites valeurs propres}\label{cercle:pvp}

 Nous allons maintenant d\'emontrer les r\'esultats concernant le spectre de 
$\Delta^{inv}$.

\begin{demo}{de \ref{cercle:th1:1}}
 Si $U$ est un vecteur colonne tel que $C\trans CU=0$, alors
$\trans UC\trans CU=0$, et donc $\|\trans CU\|=0$. Par cons\'equent, 
$\dimension \Ker C\trans C=\dimension \Ker\trans C=\dimension \Ker C=d'$.
Comme $\dimension \Ker \Delta^1_{inv}=1+\dimension \Ker C\trans C$, 
on a bien $\dimension \Ker \Delta^1_{inv}=d'+1$.
\end{demo}

\begin{demo}{de \ref{cercle:th1:2}}
C'est une conséquence directe du lemme \ref{cercle:coubure:lem}: comme
la trace de $\Delta^1_{inv}$ est celle de $C\trans C$, cette trace
est majorée en fonction de $a$ et $B$. Comme les valeurs propres de
$\Delta^1_{inv}$ sont positives, chacune est majorée.
\end{demo}

\begin{demo}{de \ref{cercle:th1:3}}

Supposons que $d\neq n$.
Soit $E_0$ le sous-espace caract\'eristique de $f$ associ\'e \`a la
valeur propre $0$. On notera $E^\perp_0$ son orthogonal pour la dualit\'e
dans l'espace des $1$-formes invariantes verticales. Comme $d\neq n$,
l'espace $E^\perp_0$ est de dimension non nulle.
 On va montrer que le quotient de
Rayleigh est uniform\'ement minor\'e sur $E^\perp_0$, pour ensuite appliquer
le principe du minimax. 

Remarques: comme les formes et les m\'etriques 
consid\'er\'ees sont invariantes, la norme ponctuelle d'une forme
ne d\'ependra pas du point o\`u on la calcule, ce qui permet d'\'ecrire
que $R(\alpha)=\frac{\|d\alpha\|^2}{\|\alpha\|^2}=
\frac{|d\alpha|^2}{|\alpha|^2}$. D'autre part, il faut noter que la
notion d'orthogonalit\'e pour la dualit\'e est ind\'ependante
de la m\'etrique. En particulier, comme $E_0$ est d\'efini
ind\'ependamment de la m\'etrique, $E^\perp_0$ le sera aussi.

Soit $V^\flat\in E^\perp_0$ et $(V_i)$ une base orthonormée 
de $\Gamma(T_VM)^{{}^G}$ telle que $(V_1,\cdots,V_d)$ soit une 
base orthonorm\'ee de $E_0$ (si $d=0$ et donc  $E_0={0}$,
on choisit alors $(V_i)$ orthonorm\'ee quelconque, la suite 
de la d\'emonstration restant valide). L'espace $E_0$ est 
stable par $f$, donc $E^\perp_0$ est stable par $\trans f$, et 
la matrice de $(\trans f)_{|E^\perp_0}$ dans la base $(V_{d+1}^\flat,\cdots,
V_n^\flat)$ est $\trans D$, o\`u $D$ est une sous-matrice de $C$. Comme
la relation (\ref{cercle:laplacien:eq1}) peut s'\'ecrire 
$\de V^\flat=-Y^\flat\wedge (\trans f)(V^\flat)$, on a 
\begin{equation}
|\de V^\flat|^2=|(\trans f)(V^\flat)|^2\geq\lambda|V^\flat|^2,
\end{equation}
o\`u 
$\lambda$ est la plus petite valeur propre de $D\trans D$. D'une part,
le d\'eterminant de cette matrice v\'erifie 
\begin{equation}
\Det D\trans D=(\Det\trans D)^2=(\Det(\trans f)_{|E^\perp_0})^2,
\end{equation}
et donc $\Det D\trans D$
est ind\'ependant du choix de la base $(V_i)$. D'autre part, 
$\Det\trans(f_{|E^\perp_0})$ est non nul. En effet, si $\trans f_{|E^\perp_0}
(\alpha) =0$, alors $\alpha\circ f=0$, donc $\alpha$ 
est orthogonal \`a l'image de $f$, qui contient les sous-espaces
caract\'eristiques de $f$ autres que $E_0$, et par cons\'equent $\alpha$
est nul. On en d\'eduit que $\lambda$ est uniform\'ement minor\'ee:
s'il existe une suite de m\'etriques telle que $\lambda
\rightarrow0$, alors la plus grande valeur propre de $D\trans D$ tend
vers l'infini (car $\Det D\trans D$ est constant), ce qui est
impossible puisque  la courbure est born\'ee et que $\Tr C\trans C\geq\Tr 
D\trans D$ (car $D$ est une sous-matrice de $C$), et donc que la somme 
des valeurs propres de $D\trans D$ est born\'ee.
 
 On a montr\'e que le quotient de Rayleigh de $\alpha\in E^\perp_0$
est minor\'e par une constante $c(f,a)$ ind\'ependante de la m\'etrique 
et du choix de $\alpha$. Comme $\dim E^\perp_0=n-d$, le principe du 
minimax nous dit donc que les $n+1-d$ plus grandes valeurs propres de
$\Delta^1_{inv}$ sont minor\'ees par $c$. Comme $\dimension \Ker
\Delta^1_{inv}=d'+1$ et que $\dimension\Omega^1(M)^G=n+1$, on en 
d\'eduit que $\lambda^{inv}_{d-d'+1,1}>c$.

Si $d=n$, alors il existe $P\in\GL_n(\R)$ tel que $P^{-1}BP$ soit
triangulaire sup\'erieure avec des $0$ sur la diagonale, et comme
$P^{-1}AP=P^{-1}\exp(B)P=\exp(P^{-1}BP)$, la matrice $P^{-1}AP$ sera
triangulaire sup\'erieure avec des $1$ sur la diagonale.
On en d\'eduit, en posant $P'=\smat{P&0\\0&I}\in\GL_{n+2}(\R)$, que
le groupe ${P'}^{-1}GP'$, o\`u $G$ est le groupe construit au paragraphe 
\ref{cercle:st}, est constitu\'e de matrices triangulaires
sup\'erieures avec des 1 sur la diagonale. C'est donc un groupe
nilpotent. 

 L'existence d'un effondrement tel que toutes les valeurs propres 
de $\Delta^1_{inv}$ tendent vers z\'ero d\'ecoulera du 
\ref{cercle:th1:4} 
\end{demo}

\begin{demo}{de \ref{cercle:th1:4}}
 On vient de d\'emontrer que si $d$ (et donc $d'$) est nul, il n'y a
pas de petites valeurs propres.

 Supposons que $d>0$. Pour simplifier, nous allons montrer le r\'esultat 
dans le cas o\`u la partie nilpotente de la r\'eduite de Jordan de $B$ 
ne comporte qu'un seul bloc de Jordan, la construction  de $g^k_\varepsilon$
\'etant semblable dans le cas g\'en\'eral.

On construit une base $(V_1,\cdots,V_n)$ de $\Gamma(T_VM)^{{}^G}$
en choisissant une base de Jordan $(V_1,\cdots,V_d)$ de $E_0$ (en
particulier, $(V_1,\cdots,V_{d'})$ sera une base de $\Ker f$) que l'on
on la compl\`ete de mani\`ere quelconque en une base  
$(V_1,\cdots,V_{d'})$ de $\Gamma(T_VM)^{{ }^G}$. On notera $C$ la matrice
de $f$ dans cette base. La matrice $C$ n'est pas de Jordan,
mais sa restriction \`a $E_0$, c'est-\`a-dire le bloc carr\'e
sup\'erieur droit de taille $d$, l'est. Elle est de la forme
\renewcommand{\arraystretch}{-.6}
\addtolength{\arraycolsep}{-0.7mm}

\begin{equation}
C=\left(\begin{array}{c|c|c}
\begin{array}{cccc}
0&\cdots&\cdots\vphantom{\vdots}&0\\
\vdots&\ddots&&\vdots\\
&&\ddots&\vdots\\
\vphantom{\vdots}&&&0\\
&&&\vdots\\
&&&\vphantom{\vdots}\\
\vdots&&&\vdots\\
0&\cdots&\cdots&0\vphantom{\vdots}
\end{array}
&
\begin{array}{cccc}
0&\cdots&\cdots&0\vphantom{\vdots}\\
\vdots&&&\vdots\\
0&&&\vphantom{\vdots}\\
1&\ddots&&\vdots\\
0&\ddots&0&0\\
\vdots&\ddots&1&0\\
\vdots&&0&1\\
0&\cdots&0&0\vphantom{\vdots}
\end{array}
& C_1\\
\hline
\begin{array}{cccc}
\vdots&&&\vdots\\
0&\cdots&\cdots&0\vphantom{\vdots}
\end{array}
&
\begin{array}{cccc}
\vdots&&&\vdots\\
0&\cdots&\cdots&0\vphantom{\vdots}
\end{array}
& C_2\\
\end{array}\right)
\end{equation}
$$\hspace{5mm}\stackrel{\underbrace{\hspace{1.9cm}}}{d'\ \mathrm{colonnes}}
\hspace{4mm}\stackrel{\underbrace{\hspace{2.1cm}}}{d-d'\ \mathrm{colonnes}}$$

o\`u $C_2$ est un bloc carr\'e de taille $n-d$ et de d\'eterminant non nul.

Soit $k\leq d-d'$. On pose $V_i^\varepsilon=
\nu_i(\varepsilon)V_i$, avec $\nu_i(\varepsilon)=\varepsilon^{-1}$ pour 
$i\geq d'+k$, et $\nu_i(\varepsilon)=\varepsilon^{-(1+d'+k-i)}V_i$ 
pour $i<d'+k$.
 La matrice $C_\varepsilon$ de $f$ dans cette base v\'erifiera
\begin{equation}
c_{ij}^\varepsilon=\frac{\nu_j(\varepsilon)}{\nu_i(\varepsilon)}c_{ij},
\end{equation}
donc $c_{ij}^\varepsilon=c_{ij}$ pour $i\geq d'+k$ (en tenant compte
du fait que $c_{ij}=0$ pour $i\geq d$ et $j\leq d$), et $c_{ij}^\varepsilon
\rightarrow 0$ quand $\varepsilon\rightarrow 0$, pour $i<d'+k$.
 La matrice $C_\varepsilon$ tend donc vers une matrice $C_0$ de la forme

\begin{equation}
C_0=\left(\begin{array}{c|c|c}
\begin{array}{cccc}
0&\cdots&\cdots\vphantom{\vdots}&0\\
\vdots&\ddots&&\vdots\\
&&\ddots&\vdots\\
\end{array}
&
\begin{array}{cccc}
0&\cdots&\cdots&0\vphantom{\vdots}\\
\vdots&&&\vdots\\
0&&&\vphantom{\vdots}\\
\end{array}
&0\\ \cline{2-3}
\begin{array}{cccc}
\vphantom{\vdots}&&&0\\
&&&\vdots\\
&&&\vphantom{\vdots}\\
\vdots&&&\vdots\\
0&\cdots&\cdots&0\vphantom{\vdots}
\end{array}
&
\begin{array}{cccc}
1&\ddots&&\vdots\\
0&\ddots&0&0\\
\vdots&\ddots&1&0\\
\vdots&&0&1\\
0&\cdots&\cdots&0\vphantom{\vdots}
\end{array}
& C'_1\\
\hline
\begin{array}{cccc}
\vdots&&&\vdots\\
0&\cdots&\cdots&0\vphantom{\vdots}
\end{array}
&
\begin{array}{cccc}
\vdots&&&\vdots\\
0&\cdots&\cdots&0\vphantom{\vdots}
\end{array}
& C'_2\\
\end{array}\right)
\end{equation}
$$\hspace{14mm}\stackrel{\underbrace{\hspace{16mm}}}{d'+k\ \mathrm{colonnes}}
\hspace{3mm}\stackrel{\underbrace{\hspace{21mm}}\hspace{12mm}}
{d-d'-k\ \mathrm{colonnes}}$$

\addtolength{\arraycolsep}{0.7mm}
\renewcommand{\arraystretch}{1}

 Comme  les $\lambda^{inv}_{i,1}$ sont ces fonctions continues de $C$,
il suffit de calculer la dimension du noyau de $C_0\trans{C_0}$, qui
est \'egale \`a celle de $\Ker C_0$.
D'une part, les  $d'+k$ premi\`eres colonnes de  $C_0$ sont nulles, 
donc $\dimension\Ker C_0\geq d'+k$, et d'autre part, comme $\det C_2\neq0$,
la famille form\'ee par les lignes $d'+k$ \`a $d-1$ et les $n-d$ derni\`eres
lignes de $C_0$ est libre, donc $\dimension\Ker C_0\leq n-(n-d)-((d-1)-
(d'+k-1))=d'+k$. De m\^eme, $\dimension\Ker C=d'$, donc on a bien exactement
$k$ petites valeurs propres.

 Si la partie nilpotente de la r\'eduite de Jordan de $B$ contient
plusieurs blocs de Jordan, on obtient le r\'esultat en proc\'edant
de la m\^eme mani\`ere pour annuler le nombre souhait\'e de lignes dans $C$.

Remarques: la famille de matrice $C_\varepsilon$
est uniformement born\'ee par rapport \`a $\varepsilon$, et le lemme 
\ref{cercle:coubure:lem} donne la majoration $|K(M,g^k_\varepsilon)|\leq 
\tau\Tr(C_\varepsilon\trans{C_\varepsilon})$, pour tout $\varepsilon$. 
La courbure sectionnelle du fibr\'e est donc bien uniform\'ement born\'ee
au cours de l'effondrement. D'autre part, on voit que si on effondre
le fibr\'e par homoth\'etie de la fibre, par exemple en posant 
$\nu_i(\varepsilon)=\varepsilon^{-1}$ pour tout $i$, la matrice 
$C_\varepsilon$ est ind\'ependante de $\varepsilon$, et donc il
n'y a pas de petite valeur propre.
\end{demo}

\begin{demo}{du corollaire \ref{cercle:cor1}}
 Si $d=d'$ et $d\neq n$, alors $\lambda^{inv}_{1,1}$ est uniform\'ement 
minor\'e d'apr\`es \ref{cercle:th1:3}. Si $d=d'$ et $d=n$, 
alors $B=0$ et toutes les valeurs propres de $\Delta^1_{inv}$ sont nulles.

 Si $d\neq d'$, alors 
\ref{cercle:th1:4} garantit l'existence d'une petite valeur propre.
\end{demo}

\begin{demo}{du théorème \ref{cercle:th2}}
 Comme $B$ est semi-simple, son orbite par conjugaison est ferm\'ee 
(\cite{cmg}, p. 28).  Comme la courbure est born\'ee, la norme de 
$C$ reste born\'ee quand la base $(V^\flat_1,\cdots,V^\flat_n,Y^\flat)$
 ---~et donc la m\'etrique~--- varie. La matrice $C$ est par
construction dans l'orbite par conjugaison de $B$, donc elle prend 
finalement ses valeurs au cours de l'effondrement dans une partie 
compacte $K$ de cette orbite. 

 La base orthonorm\'ee $(V^\flat_1,\cdots,V^\flat_n,Y^\flat)$ de 
$\Omega^1(M)^G$
engendre, par produit ext\'erieur, une base orthonorm\'ee de $\Omega^*(M)^G$.
Les coefficients de la matrice de la diff\'erentielle ext\'erieure $\de$ dans
cette base, et donc ceux de la matrice de $\Delta=\de\codiff+\codiff\de$, 
sont des fonctions continues de $C\subset K$. Par cons\'equent, quand 
la m\'etrique varie, la matrice de $\Delta$ prend ses valeurs dans un 
compact image de $K$. 
S'il existe une famille de m\'etriques telle que $\lambda^{inv}_{1,p}$
tende vers z\'ero pour un $p\in[1,n]$, alors la matrice de $\Delta$ tend vers
une matrice de rang strictement inf\'erieur, ce qui est impossible puisque,
par compacit\'e, la matrice limite sera dans l'image de $K$, donc de
m\^eme rang que $\Delta$.

 Par cons\'equent, l'op\'erateur $\Delta$ restreint \`a $\Omega^*(M)^G$ 
n'admet pas de petite valeur propre.
\end{demo}

\section{Variétés de petites dimensions}

En petite dimension, on peut être plus précis que les
théorèmes \ref{cercle:th1} et \ref{cercle:th2}, et mettre en évidence
un lien simple entre l'existence de petites valeurs propres et la 
structure du groupe $G$:

\begin{cor}\label{cercle:cor2}
Supposons que $n=2 \ou 3$. S'il existe $p\in[1,n]$ et une suite de 
métriques homogènes sur $M$ telle que la courbure sectionnelle 
associée soit uniformément bornée et que $\lambda^{inv}_{1,p}$ 
tende vers $0$, alors $G$ est nilpotent.
\end{cor}

\begin{remarque}
C'est par exemple la situation exposée en \ref{heisenberg}, où on a 
$p=1$ et $n=2$.
\end{remarque}

\begin{demo}{du corollaire \ref{cercle:cor2}}
 Montrons d'abord que s'il existe $p$ tel que $\lambda^{inv}_{p,1}
\rightarrow0$, alors $d\neq d'$.

 Si $p=1$, cela d\'ecoule du corollaire \ref{cercle:cor1}. Si $p=n$, on 
est ramen\'e \`a la situation $p=1$ par dualit\'e de Hodge.  

Reste les cas $p=2$ et $n=3$. On a d\'ej\`a calcul\'e les matrices
de $\codiff:\Omega^2(M)^G\rightarrow\Omega^1(M)^G$ et
$\de:\Omega^1(M)^G\rightarrow\Omega^2(M)^G$. On en d\'eduit que la
matrice de $\codiff\de$, en restriction \`a $\Omega^2(M)^G$ est de la forme,
dans les bases introduites au paragraphe \ref{cercle:laplacien},
\begin{equation}
\de\codiff:\left(\begin{array}{cc}
\trans CC & 0\\ 0&0 \end{array}\right).
\end{equation}

Comme la vari\'et\'e est de dimension $4$, l'op\'erateur de Hodge $*$ est
une isom\'etrie de $\Omega^2(M)^G$. En restriction \`a $\Omega^2(M)^G$,
on aura $\codiff\de=*\de\codiff*$, et donc $\codiff\de$ et $\de\codiff$
ont m\^eme spectre. D'autre part, d'apr\`es la th\'eorie de Hodge, 
le spectre du laplacien est la r\'eunion des spectres de $\codiff\de$ et
$\de\codiff$, on d\'eduit de ce qui pr\'ec\`ede qu'une petite valeur propre 
non nulle de $\Delta^2_{inv}$ sera une petite valeur propre non nulle de 
$\de\codiff{|\Omega^2(M)^G}$, et donc une petite valeur propre de $\trans CC$.
Le raisonnement appliqu\'e \`a $C\trans C$  dans la d\'emonstration 
de \ref{cercle:th1:4} reste valable pour $\trans CC$. On peut
donc conclure que si $\lambda^{inv}_{2,1}$ tend vers z\'ero, alors 
$d\neq d'$.

 Supposons que $d\neq d'$. Alors le noyau de $B$ est non trivial,
par cons\'equent $d>d'>0$ et la multiplicit\'e de la valeur
propre $0$ de $B$ est au moins \'egale \`a deux. Si $n=3$ la troisi\`eme
valeur propre est \'egale \`a la trace de $B$ qui est nulle puisqu'elle 
est r\'eelle et que $\exp(\Tr B)=\det(\exp B)=\det A=1$. Donc $d=n$,
et $G$ est nilpotent, d'apr\`es \ref{cercle:th1:3}.
\end{demo}

\section{Homologie du fibré}\label{cercle:betti}
Nous allons ici montrer que l'on peut calculer le premier nombre 
de Betti du fibré $M$
indépendamment de la cohomologie, en utilisant le fait que le 
réseau $\Gamma$ est isomorphe au groupe fondamental de $M$.
\begin{theo}\label{cercle:th3}
Soit $M=\Gamma\backslash G$ un fibré en tore $T^n$ sur le cercle
construit selon \ref{cercle:th1}, définit par une matrice 
$A\in\SL_n(\Z)$. Alors le premier nombre de Betti de $M$
est $b_1(M)=1+\dimension\Ker(A-I)$.
\end{theo}
On en déduit:
\begin{cor}\label{cercle:cor3}
Si $1$ n'est pas valeur propre de $A$, alors les $1$-formes harmoniques
de $M$ sont $G$-invariantes.
\end{cor}
On verra au paragraphe \ref{cercle:ex2} un exemple qui montre ---~entre
autres choses~--- qu'on peut effectivement, dans certains cas, avoir
des formes harmoniques qui ne sont pas invariantes.

\begin{demo}{du théorème \ref{cercle:th3}}
Comme $G$ est simplement connexe, le réseau $\Gamma$ est isomorphe au
groupe fondamental du quotient $M=\Gamma\backslash G$. Pour déterminer
le premier nombre de Betti de $M$, on va calculer l'abélianisé de son
groupe fondamental.

Rappelons que le groupe $\Gamma$ est de la forme:
\begin{equation}
(x_1,\cdots,x_n,y)\longmapsto\left(\begin{array}{c|cc}
A^y&\vecdots00&\vecdots{x_1}{x_n}\\\hline
0&\begin{smallmatrix}\\1\\0\end{smallmatrix}
&\begin{smallmatrix}\\y\\1\end{smallmatrix}
\end{array}\right),
\end{equation}
avec $(x_1,\ldots,x_n,y)\in\Z^{n+1}$. Dans la suite de la démonstration,
nous noterons les éléments de $\Gamma$ indifféremment sous la forme
de vecteurs lignes ou de vecteurs colonnes.

Soit $g=(x_1,\ldots,x_n,y)$ et $g'=(x_1',\ldots,x_n',y')$ deux éléments
de $\Gamma$. Leur inverse respectives sont:
\begin{equation}
g^{-1}=\left(\begin{array}{c}A^{-y}\left(\vecdots{x_1}{x_n}\right)\\
-y\end{array}\right)\et 
g'^{-1}=\left(\begin{array}{c}A^{-y'}\left(\vecdots{x_1'}{x_n'}\right)\\
-y'\end{array}\right).
\end{equation}
Le calcul de leur commutateur $[g,g']=gg'g^{-1}g'^{-1}$ donne:
\begin{equation}
[g,g']=\left(\begin{array}{c}\left(\vecdots{x_1}{x_n}\right)
+A^y\left(\vecdots{x_1'}{x_n'}\right)-A^{y'}\left(\vecdots{x_1}{x_n}\right)
-\left(\vecdots{x_1'}{x_n'}\right)\\0\end{array}\right),
\end{equation}
soit $[g,g']=((A^y-I)\left(\vecdots{x_1'}{x_n'}\right)-(A^{y'}-I)
\left(\vecdots{x_1}{x_n}\right),0)$.

On voit que $[\Gamma,\Gamma]\subset((A-I)\Z^n,0)$. Réciproquement, si on
fixe $g=(0,\ldots0,1)$ et qu'on fait varier $(x_1',\ldots,x_n')$,
on obtient que $[\Gamma,\Gamma]\supset((A-I)\Z^n,0)$, et donc
$[\Gamma,\Gamma]$ est exactement $((A-I)\Z^n,0)$. 

L'image de $\Z^n$ par $(A-I)$ est un sous-réseau d'indice fini du
réseau des entiers de $\Ima(A-I)$, donc le quotient de $(\Z^n,0)$
par $((A-I)\Z^n,0)$ est de la forme $Z^k\times H$, où $k=\codimension\Ima(A-I)$
et $H$ est un groupe fini ---~éventuellement trivial~---, et finalement
l'abélianisé
$\Gamma'=\Gamma/[\Gamma,\Gamma]$ du groupe $\Gamma$ est donc de la forme 
$Z^{k+1}\times H$.
Le premier nombre
de Betti de $M$ est donc 
\begin{equation}
b_1(M)=1+\dimension\Ker(A-I).
\end{equation}
\end{demo}

\begin{demo}{du corollaire \ref{cercle:cor3}}
Si $1$ n'est pas valeur propre de $A$, alors $0$ n'est pas valeur
propre de $B$, et par conséquent, en vertu du point \ref{cercle:th1:1},
$b_1(M)=\dimension\Ker\Delta_{inv}^1=1$. Toutes les $1$-formes
harmoniques sont donc dans le noyau de $\Delta_{inv}^1$, et en particulier
sont $G$-invariantes.
\end{demo}
\section{Exemples}
\subsection{Petites valeurs propres pour les $2$-formes différentielles}
\label{cercle:ex1}
Nous allons donner ici un exemple de fibr\'e en tore sur le cercle 
pour lequel $d=d'=0$ et $\Delta^2_{inv}$ admet une petite valeur propre.
Cet exemple montre que le corollaire \ref{cercle:cor1} ne se g\'en\'eralise
pas \`a $n$ et $p$ quelconque.

On d\'efinit le fibr\'e consid\'er\'e par la matrice
\begin{equation}
A=\left(\begin{array}{c|c} A' & A''\\ \hline 0 & A' \end{array}\right),
\end{equation}
avec
\begin{equation}
A'=\left(\begin{array}{cc} 2 & 1\\ 1 & 1\end{array}\right) \et
A''=\left(\begin{array}{cc} 0 & 1\\ 0 & 0\end{array}\right).
\end{equation}

\begin{fait} La matrice $A$ est semblable \`a une matrice de la forme
\begin{equation}
\left(\begin{array}{cc|cc}
e^\lambda&1&0&0 \\
0&e^\lambda&0&0\\ \hline
0&0&e^{-\lambda}&1\\
0&0&0&e^{-\lambda}\end{array}\right),
\end{equation}
o\`u $\lambda$ est un r\'eel non nul.
\end{fait}

\textbf{Démonstration :}
La matrice $A'$ admet deux valeurs propres
r\'eelles positives, qui sont inverses l'une de l'autre car $\Det A'=1$. 
On notera $\lambda$ le r\'eel positif tel que ces deux valeurs propres
soient $e^\lambda$ et $e^{-\lambda}$. Elles sont aussi valeurs propres
de $A$ avec la multiplicit\'e deux. On peut v\'erifier que le polyn\^ome
caract\'eristique de $A$ est son polyn\^ome minimal. 
Les sous-espaces propres de $A$ sont donc tous les deux de dimension 1,
et par cons\'equent, les deux blocs de sa r\'eduite de Jordan sont 
$\smat{e^\lambda&1\\0&e^\lambda}$ et $\smat{e^{-\lambda}&1\\0&e^{-\lambda}}$.
\hfill\carrenoir\nolinebreak\vspace{2mm}

\begin{fait} Il existe une suite de m\'etrique $g_\varepsilon$ sur
$M=\Gamma\backslash G(B)$ et une suite de matrices $C_\varepsilon$ 
associ\'ees telles que
\begin{equation}
C_\varepsilon=\left(\begin{array}{cc|cc}
\lambda&\epsilon&0&0 \\
0&\lambda&0&0\\ \hline
0&0&-\lambda&\epsilon\\
0&0&0&-\lambda\end{array}\right).\end{equation}
\end{fait}

\textbf{Démonstration :}
Comme  on a 
$$\exp\left(\begin{array}{cc}\lambda&e^{-\lambda}\\
0&\lambda\end{array}\right)=
\left(\begin{array}{cc}e^\lambda&1\\0&e^\lambda\end{array}\right) \et
\exp\left(\begin{array}{cc}-\lambda&e^\lambda\\
0&-\lambda\end{array}\right)=
\left(\begin{array}{cc}e^{-\lambda}&1\\0&e^{-\lambda}\end{array}\right),$$
La matrice $A$ admet un logarithme semblable \`a 
\begin{equation}
C=\left(\begin{array}{cc|cc}
\lambda&e^{-\lambda}&0&0 \\
0&\lambda&0&0\\ \hline
0&0&-\lambda&e^\lambda\\
0&0&0&-\lambda\end{array}\right).
\end{equation}
Soit $(V_1,V_2,V_3,V_4)$ la base dans laquelle la matrice de l'endomorphisme
$f$ est \'egal \`a $C$. Si on pose $V^\varepsilon_i=\varepsilon^\alpha V_i$
pour $i=1,3$, $V^\varepsilon_2=\varepsilon^{\alpha+1}
e^\lambda V_2$ et
$V^\varepsilon_4=\varepsilon^{\alpha+1}e^{-\lambda}V_4$ où $\alpha$ est
un réel strictement positif, la matrice de $f$
dans cette base sera
\begin{equation}
C_\varepsilon=\left(\begin{array}{cc|cc}
\lambda&\epsilon&0&0 \\
0&\lambda&0&0\\ \hline
0&0&-\lambda&\epsilon\\
0&0&0&-\lambda\end{array}\right).
\end{equation}
Il suffit donc de d\'efinir $g_\varepsilon$ en posant que la
base $(V^\varepsilon_1,V^\varepsilon_2,V^\varepsilon_3,V^\varepsilon_4,Y)$ 
est orthonorm\'ee. Le fait que la courbure reste born\'ee quand 
$\varepsilon\rightarrow 0$ d\'ecoule du lemme \ref{cercle:coubure:lem}
\hfill\carrenoir\nolinebreak\vspace{2mm}

\begin{fait} 
La valeur propre $\lambda^{inv}_{2,1}(M,g_\varepsilon)$ tend vers z\'ero
quand $\varepsilon\rightarrow 0$.
\end{fait}
 
\textbf{Démonstration :}

On va calculer la matrice de $\de:\Omega^2(M)^G\rightarrow\Omega^3(M)^G$ dans
des bases de la forme $(V^\flat_i\wedge V^\flat_j,V^\flat_i\wedge Y^\flat)$ et 
$(V^\flat_i\wedge V^\flat_j\wedge Y^\flat,V^\flat_i\wedge V^\flat_j 
\wedge V^\flat_k)$.

En utilisant (\ref{cercle:laplacien:eq1}), on obtient que 
$\de V^\flat_i\wedge Y^\flat=0$ pour tout $i$, et que
\begin{eqnarray}
\de(V^\flat_1\wedge V^\flat_2)&=&(\lambda V^\flat_1+\varepsilon V^\flat_2)
\wedge Y^\flat\wedge V^\flat_2
-V^\flat_1\wedge\lambda V^\flat_2\wedge Y^\flat\nonumber\\
&=& -2\lambda V^\flat_1\wedge V^\flat_2\wedge Y^\flat,\\
\de(V^\flat_1\wedge V^\flat_3)&=&(\lambda V^\flat_1+\varepsilon V^\flat_2)
\wedge Y^\flat\wedge V^\flat_3
-V^\flat_1\wedge(-\lambda V^\flat_3+\varepsilon V^\flat_4)\wedge Y^\flat
\nonumber\\
&=& -\varepsilon V^\flat_2\wedge V^\flat_3\wedge Y^\flat
-\varepsilon V^\flat_1\wedge V^\flat_4 \wedge Y^\flat,\\
\de(V^\flat_1\wedge V^\flat_4)&=&(\lambda V^\flat_1+\varepsilon V^\flat_2)
\wedge Y^\flat\wedge V^\flat_4
-V^\flat_1\wedge(-\lambda V^\flat_4)\wedge Y^\flat\nonumber\\
&=&-\varepsilon V^\flat_2\wedge V^\flat_4 \wedge Y^\flat,\\
\de(V^\flat_2\wedge V^\flat_3)&=&\lambda V^\flat_2\wedge Y^\flat
\wedge V^\flat_3-V^\flat_2\wedge(-\lambda V^\flat_3+\varepsilon V^\flat_4)
\wedge Y^\flat\nonumber\\
&=&-\varepsilon V^\flat_2\wedge V^\flat_4\wedge Y^\flat,\\
\de(V^\flat_2\wedge V^\flat_4)&=&\lambda V^\flat_2\wedge Y^\flat
\wedge V^\flat_4-V^\flat_2\wedge(-\lambda V^\flat_4)\wedge Y^\flat=0,\\
\de(V^\flat_3\wedge V^\flat_4)&=&(-\lambda V^\flat_3+\varepsilon V^\flat_4)
\wedge Y^\flat\wedge V^\flat_4
-V^\flat_3\wedge(-\lambda V^\flat_4)\wedge Y^\flat\nonumber\\
&=&2\lambda V^\flat_3\wedge V^\flat_4\wedge Y^\flat.
\end{eqnarray}
La matrice de $\de$ dans les bases 
$$(V^\flat_1\wedge V^\flat_2,V^\flat_1\wedge V^\flat_3,V^\flat_1\wedge V^\flat_4, 
V^\flat_2\wedge V^\flat_3,V^\flat_2\wedge V^\flat_4,V^\flat_3\wedge V^\flat_4,
V^\flat_i\wedge Y^\flat)$$
et 
$$\begin{array}{l}
(V^\flat_1\wedge V^\flat_2\wedge Y^\flat,V^\flat_1\wedge V^\flat_3\wedge 
Y^\flat,
V^\flat_1\wedge V^\flat_4\wedge Y^\flat,V^\flat_2\wedge V^\flat_3\wedge 
Y^\flat,\\
V^\flat_2\wedge V^\flat_4\wedge Y^\flat,V^\flat_3\wedge V^\flat_4\wedge 
Y^\flat,
V^\flat_i\wedge V^\flat_j \wedge V^\flat_k)
\end{array}$$
est de la forme
\begin{equation}
\left(\begin{array}{c|c}
\begin{array}{rrrrrr}-2\lambda&0&0&0&0&0\\0&0&0&0&0&0\\
0&-\varepsilon&0&0&0&0\\0&-\varepsilon&0&0&0&0\\
0&0&-\varepsilon&-\varepsilon&0&0\\0&0&0&0&0&2\lambda
\end{array}&0\\ \hline
0&0\end{array}\right).
\end{equation}
On voit que quand $\varepsilon$ tend vers z\'ero, cette matrice tend
vers une matrice de rang strictement inf\'erieur. On peut en d\'eduire 
comme au paragraphe \ref{cercle:pvp} que l'op\'erateur $\codiff\de$ admet 
une petite valeur propre, qui sera aussi petite valeur propre de $\Delta$.
\hfill\carrenoir\nolinebreak\vspace{2mm}

\subsection{Structure homogène non abélienne sur le tore}\label{cercle:ex2}
Nous allons ici étudier plus en détail un exemple particulier de fibré
en tore $T^2$ sur le cercle, pour mettre en évidence plusieurs de
ses propriétés.

Ce fibré est construit par le théorème \ref{cercle:th1}, avec la donnée
de 
\begin{equation}
A=\left(\begin{array}{cc}1&0\\0&1\end{array}\right) \et
B=\left(\begin{array}{cc}0&2\pi\\-2\pi&0\end{array}\right).
\end{equation}
La matrice $\exp(xB)$ est de la forme 
\begin{equation}
\left(\begin{array}{cc}\cos2\pi x&\sin2\pi x\\
-\sin2\pi x&\cos2\pi x\end{array}\right).
\end{equation}
C'est donc la matrice d'une rotation d'angle
$2\pi x$ ; nous la noterons $R(2\pi x)$.

Le groupe $G$ s'écrit:
\renewcommand{\arraystretch}{1.2}
\begin{equation}\label{cercle:ex2:df}
\left(\begin{array}{c|c}
R(2\pi x)&\begin{smallmatrix}0&y\\0&z\end{smallmatrix}\\\hline
0&\begin{smallmatrix}1&x\\0&1\end{smallmatrix}
\end{array}\right),\ x,y,z\in\R.
\end{equation}

Un première remarque est de constater que la variété $M=\Gamma\backslash G$
est un tore:
\begin{fait}
$\Gamma$ est isomorphe à $\Z^3$ et $\Gamma\backslash G$ est
difféomorphe à $T^3$.
\end{fait}
En effet, le réseau $\Gamma$ s'écrit:
\begin{equation}
\left(\begin{array}{c|c}
\begin{smallmatrix}1&0\\0&1\end{smallmatrix}
&\begin{smallmatrix}0&y\\0&z\end{smallmatrix}\\\hline
0&\begin{smallmatrix}1&x\\0&1\end{smallmatrix}
\end{array}\right),\ x,y,z\in\Z.
\end{equation}
\renewcommand{\arraystretch}{1}
On peut vérifier que $\Gamma$ est bien abélien. Par ailleurs, comme
la topologie du fibré est entièrement déterminée par la matrice $A$, 
le fibré est bien trivial. On peut remarquer le fait que $\Gamma$
soit abélien est cohérent avec le fait que ce soit le groupe
fondamental d'un tore.

On a donc construit un groupe de Lie résoluble simplement connexe
dont un sous-groupe cocompact est commutatif.  En comparaison, on a
pour les groupes nilpotents le résultat suivant (\cite{ra}):
\begin{theo}
Soit $N_1$ et $N_2$ deux groupes de Lie nilpotents simplement connexes, 
et $\Gamma_1$, $\Gamma_2$ deux sous-groupes cocompacts de $N_1$ et
$N_2$ respectivement. Alors tout isomorphisme entre $\Gamma_1$ et
$\Gamma_2$ s'étend en un isomorphisme entre $N_1$ et $N_2$.
\end{theo}
En particulier, si un groupe nilpotent simplement connexe contient
un sous-groupe cocompact isomorphe à $\Z^n$, alors il est abélien.
Le groupe $G$ illustre donc le fait que ce théorème ne se généralise
pas aux groupes résolubles.

D'autre part, comme $M=\Gamma\backslash G$ et un tore, son premier
nombre de Betti est égal à sa dimension, donc $b_1(M)=3$.
Or, selon le théorème \ref{cercle:th1}, si on munit $G$ d'une métrique
invariante le noyau de $\Delta_{inv}^p$ est de dimension $1$. On en conclut :
\begin{fait}
Soit $g$ une métrique $G$-invariante à gauche sur $M$. Alors il existe
sur $M$ des $1$-formes harmoniques qui ne sont pas $G$-invariantes.
\end{fait}
On voit donc que la proposition \ref{pr:lott} et le corollaire 
\ref{cercle:cor3} ne se généralisent pas à toutes les solvariétés. 
Réciproquement, le groupe $G$ illustre la situation où la
multiplicité de la valeur propre $1$ dans $A$ est strictement supérieure
à la multiplicité de $0$ dans $B$.

Enfin, on peut remarquer que d'après les formules (\ref{cercle:courbure:eq1})
et (\ref{cercle:courbure:eq2}), pour la métrique invariante sur $M$ 
telle que la base $(\dep{}x,\dep{}y,\dep{}z)$ soit orthonormée à l'origine,
la courbure sur $M$ est nulle (on peut vérifier que c'est en général faux
pour une métrique invariante quelconque). On va reformuler ce résultat
et en donner une démonstration très simple qui ne fait pas appel aux 
formules du paragraphe \ref{cercle:courbure}:
\begin{fait}
Il existe sur $\R^3$ une métrique invariante pour la structure canonique
de groupe abélien, et invariante pour l'action à gauche du groupe $G$.
\end{fait}
\textbf{Démonstration :}
On considère sur $\R^3$ la métrique euclidienne canonique, et on note
$x$, $y$ et $z$ les cordonnées canoniques. Si $a$, $b$ et $c$ sont des
réels fixés, la paramétrisation de $G$ donnée par (\ref{cercle:ex2:df}) 
définit l'action à gauche de $(a,b,c)$ comme étant
\begin{equation}
\left(\begin{array}{c}
a\\b\\c\end{array}\right)\cdot
\left(\begin{array}{c}x\\y\\z\end{array}\right)=
\left(\begin{array}{c}
a+x\\
\left(\begin{array}{c}b\\c\end{array}\right)
+R(2\pi a)\left(\begin{array}{c}y\\z\end{array}\right)\end{array}\right).
\end{equation}
On voit que l'action de $(a,b,c)$ est la composée d'une rotation et
d'une translation. C'est donc une isométrie pour la norme euclidienne
canonique. Par conséquent, cette norme est invariante à gauche pour
l'action de $G$.
\hfill\carrenoir

Remarquons pour finir qu'on peut facilement généraliser cet exemple en
dimension supérieure en construisant une matrice $A$ contenant un bloc
de la forme $\smat{1&0\\0&1}$, par exemple :
\begin{equation}
A=\left(\begin{array}{cccc}1&0&0&0\\0&1&0&0\\0&0&1&1\\0&0&0&1
\end{array}\right)\et B=\left(\begin{array}{cccc}0&2\pi&0&0\\-2\pi&0&0&0\\
0&0&0&1\\0&0&0&0\end{array}\right).
\end{equation}
Le fibré obtenu dans ce cas sera une solvariété ayant une topologie
de nilvariété, mais dont les formes harmoniques ne sont pas toutes
invariantes.

\chapter{Effondrements homogènes de fibrés principaux en tores sur le tore}
\label{pvp:tore}

\section{Topologie et spectre du fibré}
Nous allons ici démontrer le théorème \ref{tore:th}.

\begin{demo}{de \ref{tore:th:1}}
 Soit $M$ un fibré principal de base $T^2$ et de fibre $F$.
La base du fibré peut s'écrire
$$[0,1]\times[0,1]_{/\sim},$$
o\`u $\sim$ est la relation d'équivalence engendrée par 
$(x,0)\sim(x,1)$ et $(0,y)\sim(1,y)$.  Le fibré $M$ peut
alors se définir par la donnée, pour tout point $p$ du bord $\partial K$
de $K=[0,1]\times[0,1]$ d'un difféomorphisme $\varphi_p$ de la fibre, et en 
posant
\begin{equation}
M=K\times F_{/(p,x)\sim (q,\varphi_q^{-1}\circ \varphi_p(x)), 
\forall x\in F, \forall p,q\in \partial K, p\sim q}.
\end{equation}
L'hypothèse de principalité se traduit ici par le fait que pour tout
$p,q\in\partial K$ tels que $p\sim q$ et pour tout $g,x\in F$, on a
\begin{equation}
(p,g\cdot x)\sim(q,g\cdot \varphi_q^{-1}\circ \varphi_p(x)),
\end{equation}
ce qui impose aux $\varphi_q^{-1}\circ \varphi_p$ d'\^etre des translations 
\`a droite sur la fibre. On peut donc, sans perte de g\'en\'eralit\'e
se restreindre, pour le choix des $\varphi_p$, au groupe des translations
de la fibre, qui est isomorphe \`a $F$. Le fibré est donc
déterminé par la donnée d'une application de $\partial K$ dans $F$.
Comme sa topologie ne dépend pas de la classe d'homotopie de cette
application, l'ensemble des fibrés principaux de fibre $F$ 
sur le tore $T^2$ est parametré par le groupe fondamental de $F$.
Il s'agit en fait d'un exemple de classe d'obstruction (\cite{steen}, \S~35) 
qui est, dans le cas général d'un $F$-fibré principal sur une variété
compacte $N$ un élément de $H_2(N,\pi_1(F))$ et qui mesure
l'obstruction du fibré à être trivial.

Considérons maintenant un fibré principal $M$ de fibre $T^n$, et
$(a_1,\cdots,a_n)\in\pi_1(T^n)=\Z^n$ sa classe d'obstruction.
Nous allons munir $\R^{n+2}$ d'une structure de groupe telle que la topologie 
du quotient à gauche par $\Z^{n+2}$ soit celle du fibré. Pour ce faire,
nous choisirons le repr\'esentant $\gamma:\partial K\rightarrow T^n$
de la classe d'obstruction de la manière suivante:
\begin{equation}\label{tore:topo1}
\begin{array}{c}
\gamma_{|\{0\}\times[0,1]}=\gamma_{|[0,1]\times\{0\}}
=\gamma_{|[0,1]\times\{1\}}=0,\\
\\
\gamma(1,t)=(ta_1,\cdots,ta_n), \forall t\in[0,1],\\
\end{array}
\end{equation}
de sorte qu'un élément $(x_1,\cdots,x_n)\in T^n$ de la fibre au dessus de
$(0,t)\in K$ sera identifié \`a l'élément $(x_1+ta_1,\cdots,x_n+ta_n)$
au dessus de $(1,t)$. Si on note $(x_1,\cdots,x_n,y_1,y_2)$ les éléments 
de $\R^{n+2}$, on veut donc d\'efinir sur cet ensemble un produit tel que
\begin{equation}\label{tore:topo2}
(k_1,\cdots,k_n,0,0)\cdot(x_1,\cdots,x_n,y_1,y_2)=
(x_1+k_1,\cdots,x_n+k_n,y_1,y_2)
\end{equation}
de sorte que d'une part les sous-espaces de $\R^{n+2}$ d'équation
$(y_1,y_2)=c^{te}$ passent au quotient comme des tores, et tel que 
\begin{multline}\label{tore:topo3}
(0,\cdots,0,l_1,l_2)\cdot(x_1,\cdots,x_n,y_1,y_2)=\\
(x_1+y_2a_1l_1,\cdots,x_n+y_2a_nl_1,y_1+l_1,y_2+l_2),
\end{multline}
de sorte que la structure de fibré en tore sera bien celle définie
 par (\ref{tore:topo1}).

On peut effectivement construire une telle structure de groupe en plongeant
$\R^{n+2}$ dans $M_{n+3}(\R)$ par l'application suivante:
\addtolength{\arraycolsep}{-1mm}
\begin{equation}
(x_1,\cdots,x_n,y_1,y_2)\longmapsto\left(\begin{array}{c|ccc}
 I_n&\vecdots00&
\vecdots{a_1y_1}{a_ny_1}&
\vecdots{x_1}{x_n}\\\hline
0&\begin{smallmatrix}\\1\\0\\0\end{smallmatrix}
&\begin{smallmatrix}\\0\\1\\0\end{smallmatrix}
&\begin{smallmatrix}\\y_1\\y_2\\1\end{smallmatrix}
\end{array}\right).
\end{equation}
\addtolength{\arraycolsep}{1mm}
Notons $G$ l'image de cette application. C'est un sous-groupe de 
$M_{n+3}(\R)$, et le quotient $\Gamma\backslash G$ où $\Gamma$
est le r\'eseau des entiers de $G$ est difféomorphe à la variété
$M$, qui est donc une nilvariété.

Supposons maitenant que $n\geq2$. On pose
$d=\mathrm{pgcd}(a_1,\cdots,a_n)$ et $a'_i=a_i/d$. Soit $P=(p_{ij})
\in M(n,\Z)$ une matrice telle que $p_{i1}=a'_i$ et que ses
vecteurs colonnes forment une base du réseau $\Z^n$. On a alors
\begin{equation}
P^{-1}\left(\vecdots{a_1}{a_n}\right)=\left(\begin{smallmatrix}d\\0\\
\vvdots\\0 \end{smallmatrix}\right) 
\end{equation}
et 
\begin{equation}
\left(\begin{array}{c|c}P^{-1}&0\\\hline0&I\end{array}\right)
\addtolength{\arraycolsep}{-1mm}
\left(\begin{array}{c|ccc}
 I_n&\vecdots00&
\vecdots{a_1y_1}{a_ny_1}&
\vecdots{x_1}{x_n}\\\hline
0&\begin{smallmatrix}\\1\\0\\0\end{smallmatrix}
&\begin{smallmatrix}\\0\\1\\0\end{smallmatrix}
&\begin{smallmatrix}\\y_1\\y_2\\1\end{smallmatrix}
\end{array}\right)
\addtolength{\arraycolsep}{1mm}
\left(\begin{array}{c|c}P&0\\\hline0&I\end{array}\right)
=\addtolength{\arraycolsep}{-1mm}
\left(\begin{array}{c|ccc}
 I_n&\vecdots00&
\begin{smallmatrix}dy_1\\0\\\vvdots\\0\end{smallmatrix}&
\vecdots{x'_1}{x'_n}\\\hline
0&\begin{smallmatrix}\\1\\0\\0\end{smallmatrix}
&\begin{smallmatrix}\\0\\1\\0\end{smallmatrix}
&\begin{smallmatrix}\\y_1\\y_2\\1\end{smallmatrix}
\end{array}\right)
\addtolength{\arraycolsep}{1mm}
\end{equation}
avec $\left(\vecdots{x'_1}{x'_n}\right)=P^{-1}
\left(\vecdots{x_1}{x_n}\right)$.
On peut voir que le groupe ${P'}^{-1}GP'$, avec $P'=\smat{P&0\\0&I}$, est
isomorphe \`a $\R^{n-1}\times G'$, o\`u $G'$ est le groupe
\begin{equation}
\left\{\left(\begin{array}{cccc}1&0&dy_1&x\\0&1&0&y_1\\0&0&1&y_2\\0&0&0&1
\end{array}\right), x,y_1,y_2\in\R\right\}.
\end{equation}
D'autre part, comme $\det P'=1$ et que $P'$ est à coefficients entiers,
le r\'eseau des matrices \`a coefficients entiers de ${P'}^{-1}GP'$ est 
exactement ${P'}^{-1}\Gamma P'$, o\`u $\Gamma$ est le réseau des entiers
de $G$.
La vari\'et\'e $M=\Gamma\backslash G$, qui est difféomorphe à
${P'}^{-1}\Gamma P'\backslash {P'}^{-1}GP'$ peut donc s'écrire
\begin{equation}
M\simeq (\Z^{n-1}\times\Gamma')\backslash(\R^{n-1}\times G')
\simeq T^{n-1}\times N,
\end{equation}
o\`u $N=\Gamma'\backslash G'$, en notant $\Gamma'=$ le réseau des 
entiers de $G'$.
\end{demo}

 Ce calcul montre qu'on peut se ramener au cas où les $a_i$, $i\geq2$
sont nuls. On supposera dans la suite que c'est le cas, et on posera
$a_1=a$.

\begin{demo}{de \ref{tore:th:2}}
  Soient $X_i$, $Y_1$ et $Y_2$ les champs de vecteurs invariants à gauche
engendr\'es en $I_{n+3}$ respectivement par
$$\Dep{}{x_i}=\left(\begin{array}{c|c}
0&\begin{smallmatrix}
\vecdots0\ccdot&\vecdots0\ccdot&\vecdots00 \\
\ccdot&\ccdot&1\\
\vecdots{\vphantom{0}\ccdot}0&\vecdots{\vphantom{0}\ccdot}0&\vecdots00
\end{smallmatrix}\\\hline 0&0\end{array}\right),\ 
\Dep{}{y_1}=\left(\begin{array}{c|c}
0&\begin{smallmatrix}\begin{smallmatrix}0\\ \\\vvdots\\0\end{smallmatrix}
\!&\begin{smallmatrix}a\\0\\\vvdots\\0\end{smallmatrix}
\!&\begin{smallmatrix}0\\ \\\vvdots\\0\end{smallmatrix}
\end{smallmatrix}\\\hline
0&\begin{smallmatrix}0\;\:&0&\;\:1\\
0\;\:&0&\;\:0\\ 0\;\:&0&\;\:0\end{smallmatrix}\end{array}\right) \et \Dep{}{y_2}
=\left(\begin{array}{c|c} 0&0\\\hline 
0&\begin{smallmatrix}0&0&0\\0&0&1\\0&0&0
\end{smallmatrix}\end{array}\right). $$

Ces champs vérifient $[X_i,X_j]=0,\ [X_i,Y_j]=0 \et 
[Y_1,Y_2]=aX_1$. On notera $V$ le vecteur $aX_1$, dont on 
peut remarquer qu'il est non nul (si $a$ est nul, le fibré est trivial).

 Soit $g$ une métrique homogène sur $M$, et $(V_i)$ une base
de $\Gamma(T_VM)^{{}^G}$ telle que $(V_1,\cdots,V_n,Y_1,Y_2)$ soit 
orthonormée en tout point et que $V_1$ soit colinéaire à $V$. Les
crochets de Lie entre les vecteurs de cette base sont:
\begin{equation}\label{tore:lie}
[V_i,V_j]=0,\ [V_i,Y_j]=0, \et [Y_1,Y_2]=\eta V_1,\ \eta\in \R^*.
\end{equation}

On en déduit:
\begin{equation}\label{tore:d}
\de V^\flat_1=-\eta Y^\flat_1\wedge Y^\flat_2\et \de V^\flat_i=
\de Y^\flat_j=0,\ i>1,
\end{equation}
o\`u $(V^\flat_1,\cdots,V^\flat_n,Y^\flat_1,Y^\flat_2)$ est la base duale de
$(V_1,\cdots,V_n,Y_1,Y_2)$. Les formes de cette base de $\Omega^1(M)^G$
engendrent, par produit extérieur, une base de $\Omega^*(M)^G$
composée de formes propres de $\Delta_{inv}$. En effet, il découle
de (\ref{tore:d}) qu'elles sont toutes fermées sauf celles de
la forme $V_1\wedge V_{i_1}\wedge\cdots\wedge V_{i_k}$ ($i_j\neq 1$), dont
la différentielle vaut:
\begin{equation}
\de(V_1\wedge V_{i_1}\wedge\cdots\wedge V_{i_k})=-\eta Y_1\wedge Y_2\wedge
V_{i_1}\wedge\cdots\wedge V_{i_k},
\end{equation}
et, puisque $\codiff=(-1)^{n(p+1)+1}*\de*$,  elles sont toutes cofermées 
sauf celles de 
la forme $Y_1\wedge Y_2\wedge V_1\wedge V_{i_1}\wedge\cdots\wedge V_{i_k}$ 
dont la codifférentielle vaut:
\begin{equation}
\codiff(Y_1\wedge Y_2\wedge V_1\wedge V_{i_1}\wedge\cdots\wedge V_{i_k})=
-\eta V_1\wedge V_{i_1}\wedge\cdots\wedge V_{i_k}.
\end{equation}
 En restriction à $\Omega^p(M)^G$, les formes de la base sont donc
harmoniques, sauf $C^{p-1}_{n-1}$ formes cofermées et $C^{p-2}_{n-1}$ 
formes fermées qui sont des formes propres de valeur propre 
égale à $\eta^2$. L'opérateur $\Delta^p_{inv}$ admet donc 
une unique valeur propre non nulle, de multiplicité $C^{p-1}_{n-1}
+C^{p-2}_{n-1}=C^{p-1}_n$ et égale à $\eta^2=|V|^2$.

Si on choisit une base orthonormée
de la forme $(V_1,\cdots,V_n,Y_1', Y_2')$, avec $Y_i'=Y_i+
\sum_{k=1}^n\xi_kV_k$, on aura toujours $[Y_1',Y_2']=[Y_1,Y_2]$.
Le résultat ne dépend donc pas du choix de la connexion sur le fibré.
Remarquons enfin que si l'on choisit une autre métrique sur la base 
(en se donnant deux
champs horizontaux  quelconques $Y'_1$ et $Y'_2$ et en les
supposant orthogonaux) on obtiendra le même résultat en rempla\c{c}ant
$V$ par $V'=[Y'_1,Y'_2]$, avec comme valeur propre $\eta^2=|V'|^2=\Vol(B)
^{-2}|V|^2$.
\end{demo}
\section{Effondrements des fibr\'es principaux sur le tore $T^2$}
\label{tore:eff}
Dans cette partie, nous allons montrer que  les fibr\'es 
construits dans le théorème \ref{tore:th} peuvent admettre, si $n\geq2$,
un effondrement \`a diam\`etre et courbure born\'es pour lequel 
$\lambda$ ne tend pas vers z\'ero.

Soit $M$ un tel fibr\'e. Nous allons d'abord montrer le lemme 
suivant qui nous permettra de contr\^oler la courbure:
\begin{lem}
Pour toute m\'etrique homog\`ene $g$ sur $M$, la courbure sectionnelle
de $M$ v\'erifie $|K(M,g)|\leq \frac34|V|^2$.
\end{lem}
\textbf{Démonstation :}
On se place dans la même base $(V_1,\cdots,V_n,Y_1,Y_2)$ que celle
utilisée dans la démonstration de \ref{tore:th:2}.
De (\ref{tore:lie}), on peut rapidement d\'eduire que $\ad_{V_i}=0$, et
que $\ad^*_{Y_i}Y_j=0$ car pour tout vecteur $U$,
$\langle\ad^*_{Y_i}Y_j,U\rangle=\langle Y_j,\ad_{Y_i}U\rangle=
\langle Y_j,[Y_i,U]\rangle=0$.
De plus, comme $\langle\ad^*_{Y_i}V_j,U\rangle=\langle V_j,[Y_i,U]\rangle$, 
on aura 
$\ad_{Y_i}^*V_j=0$ pour $j\neq 1$, $\ad_{Y_1}^*V_1=\mu Y_2$ et
$\ad_{Y_2}^*V_1=-\mu Y_1$.

La formule (\ref{cercle:courbure:eq}) donne donc:
\begin{equation}
K(V_i,V_j)=0,\ K(Y_1,Y_2)=-\frac34\|[Y_1,Y_2]\|^2,
\end{equation}
\begin{equation}
K(V_1,Y_i)=\frac{\mu^2}4, \et K(V_i,Y_j)=0 \mbox{ pour } i\neq1.
\end{equation}
Comme d'une part $\mu^2=\|V\|^2$, et d'autre part $V=[Y_1,Y_2]$, 
la majoration du lemme en d\'ecoule imm\'ediatement.
\hfill\carrenoir\nolinebreak\vspace{2mm}

Un effondrement \`a base fixe du fibr\'e $M$ par des m\'etriques homog\`enes
est d\'etermin\'e par une famille de bases $(V_1^\varepsilon,\cdots,
V_n^\varepsilon)$ de $\Gamma(T_VM)^{{}^G}$ (remarque: on ne suppose plus
ici que $V_1^\varepsilon$ est colinéraire à $V$). 
Nous allons pr\'esenter ici des
exemples d'effondrements associ\'es \`a des familles de bases de la forme
$(V_1^\varepsilon,\cdots,V_n^\varepsilon)=(\varepsilon^{-\alpha_1} V_1,
\cdots,\varepsilon^{-\alpha_n}V_n)$, o\`u $(V_1,\cdots,V_n)$ est une base
fix\'ee. Si $b_i$ et $b_i^\varepsilon$ sont les coefficients de $V$ dans 
les bases respectives $(V_1,\cdots,V_n)$ et $(V_1^\varepsilon,\cdots,
V_n^\varepsilon)$, on aura
\begin{equation}
b_i^\varepsilon=\varepsilon^{\alpha_i}b_i.
\end{equation}

\begin{exemple}\label{tore:ex1}
Si $\alpha_i>0$, pour tout $i$, alors le diam\`etre de la fibre tend
vers $0$, ainsi que les $b_i^\varepsilon$. On a donc un effondrement
\`a courbure born\'ee du fibr\'e sur la base $T^2$, et la valeur
propre $\lambda=\|V\|^2=\sum_{i=1}^n{(b_i^\varepsilon)}^2$ tend vers z\'ero.
\end{exemple}

On peut cependant construire des effondrements pour lesquels le
comportement du spectre est diff\'erent, et en particulier tels qu'il
n'y ait pas de petites valeurs propres :

\begin{exemple}\label{tore:ex2}
 Supposons que
$\alpha_i=0$, pour tout $i>1$, $\alpha_1>0$, et que les composantes
de $V_1$ dans la base $(X_1,\cdots,X_n)$ soient irrationnelles entre elles.
Une droite de la fibre de direction $V_1$ sera donc dense dans la fibre,
et par cons\'equent, il suffit que seul $\alpha_1$ soit non
nul pour que la fibre s'effondre sur un point. On aura alors
$b_i^\varepsilon=b_i$ pour $i>1$, et $b_1^\varepsilon\rightarrow0$.
La courbure reste donc born\'ee et $\lambda\rightarrow\sum_{i>1}b_i^2\neq 0$.
\end{exemple}

Ce dernier exemple justifie remarque \ref{tore:rem2} faite dans
l'introduction.

\section{Exemples de fibrés principaux sur des bases de dimension
strictement supérieures à 2}

Les exemples du théorème \ref{tore:th} sont relativement
simples, en ce sens que l'essentiel de la topologie est contenu
dans la structure de fibré en cercle de la nilvariété $N$. On peut
cependant facilement contruire des fibrés principaux en tore dont
le comportement du spectre est un peu plus riche:

\begin{exemple}
Considérons deux fibrés principaux $M_1$ et $M_2$, de fibre respective
$T^{k_1}$ et $T^{k_2}$ et de base $T^2$, muni d'une structure homogène
et d'une métrique invariante. Soient $\lambda_1$ et $\lambda_2$ leur
valeur propre associée définie par le théorème \ref{tore:th}.
La variété $M$ définie par le produit riemannien
$M=M_1\times M_2$ est un fibré principal de fibre
$T^{k_1}\times T^{k_2}=T^{k_1+k_2}$ et de base $T^2\times T^2=T^4$.
D'après la formule de Künneth, il admet $\lambda_1$ et $\lambda_2$
comme valeurs propres. En choissant sur $M_1$ et $M_2$ des suites
de métriques qui effondrent ces variétés sur $T^2$ la suite des
métriques produits effondre $M$ sur $T^4$.
\end{exemple}
 On voit que sur cet exemple, on peut avoir deux valeurs propres non
nulles distinctes en restriction aux formes invariantes. De plus, on
peut choisir des effondrements sur $M_1$ et $M_2$ tels que ces deux
valeurs propres tendent vers zéro à des vitesses différentes, ou
même, en choisissant pour $M_1$ et $M_2$ les effondrements
décrits dans les exemples \ref{tore:ex1} et \ref{tore:ex2} 
respectivement, tels que seule
l'une de ces valeurs propres tende vers zéro. Enfin, on peut remarquer
que considérer un effondrement de $M$ par homothétie de la fibre
revient à considérer des homothétie des fibres de $M_1$ et $M_2$,
et qu'on retrouve le fait remarqué en \ref{tore:rem1} que cet
effondrement produit des petites valeurs propres.

\chapter{Topologie des fibrés principaux en tores}\label{topo}
Nous allons dans cette partie nous attacher à décrire la topologie
des fibrés principaux en tore, et en particulier à construire
un invariant différentiel qui permettra, comme la classe d'Euler
dans la cas des fibrés en cercles, d'étudier le comportement du spectre
du laplacien lors d'un effondrement.

Soit $M$ un fibré principal en tore $T^k$ sur une base $N$. Le tore
$T^k=\R^k/\Z^k$ peut s'écrire comme le produit de $k$ cercles~:
$T^k=\prod_{i=1}^kS_{(i)}^1$. L'action de $T^k$ sur $M$ induit une action 
de chacun des $S_{(i)}^1$. On peut donc définir les variétés 
\begin{equation}\label{topo:eq1}
M_i=M/\prod_{j\neq i}S_{(j)}^1,
\end{equation}
chaque $M_i$ étant un fibré en cercle de
base $N$ sur lequel agit $S_{(i)}^1$. Réciproquement, la donnée des 
$k$ fibrés en cercles $(M_i)_{i\leq k}$ sur $N$ permet de construire 
un fibré en tore $T^k$ en prenant la somme de Whitney 
$\bigoplus_{i=1}^k M_i$ de ces fibrés en cercles, ce fibré en
tore étant difféomorphe au fibré $M$. Comme la structure
d'un fibré en cercle est déterminé par sa classe d'Euler, la topologie
de $M$ est déterminée par la donnée d'un $k$-uplet $(e_1,\cdots,
e_k)\in (H^2(N,\Z))^k$ de classes d'Euler. Cependant, la décomposition
de $T^k$ en produit de cercles n'est pas unique. En effet, pour chaque
base $(a_1,\cdots, a_k)$ du réseau $\Z^k$, on peut écrire $T^k$
comme le produit de la famille de cercles $(\R a_i/\Z a_i)_i$, auquel
correspond en général un $k$-uplet différent de classes d'Euler.

En homologie simpliciale, on peut définir la classe
d'obstruction $[c]$ d'un fibré $F\hookrightarrow M\rightarrow N$,
où $[c]$ est un élément de $H_2(N,\pi_1(F))$ qui est une mesure de
l'obstruction du fibré à admettre une section (voir \cite{steen},
\S 35). Si la fibre $F$ est un tore, les groupes
d'homotopies $\pi_n(F)$ sont triviaux pour $n\geq2$, et par conséquent
cette classe d'obstruction est nulle si et seulement si le fibré
admet une section (\cite{steen}, \S 29 et \S 35) ce qui, dans le cas
d'un fibré principal est équivalent à être trivial. Dans le cas d'une
fibre $T^k$, cette classe d'obstruction est un élément de $H_2(N,\Z^k)$.
On veut définir un objet semblable pour la cohomologie de Rham.

Dans le cas d'un fibré en cercle de classe d'Euler $[e]$, on a la 
propriété suivante (\cite{bt}, p.72)~: si $\omega$ est une $1$-forme 
verticale invariante
dont l'intégrale sur chaque fibre vaut $1$, alors $\de\omega$ est
une $2$-forme horizontale qui dépend du choix de la connexion sur le fibré, 
mais qui est, au signe près, le relevé d'une élément de $[e]$. Dans le
cas d'un fibré en tore, on va construire une invariant qui généralise
cette propriété. 

 Rappelons tout d'abord que si on note $\mathcal G$ l'algèbre de
Lie de $T^k$, l'action de $T^k$ sur le fibré $M$ induit un plongement de 
$\mathcal G$ dans l'espace des champs de vecteurs verticaux invariants
de $M$. Ce plongement s'étend naturellement aux tenseurs sur 
$\mathcal G$ en une application 
$$(\bigotimes^p\mathcal G)\otimes(\bigotimes^q\mathcal G^*)\rightarrow
\Gamma\left((\bigotimes^pT^VM)\otimes(\bigotimes^q{T^V}^*M)\right)=
\Gamma({T^V}^p_qM).$$
Le tenseur obtenu ne définit pas de manière canonique un élément de 
$\Gamma(T^p_qM)$ si $q\neq0$, mais si on se donne une connexion sur le
fibré $M$, la partie covariante du tenseur est bien définie sur $TM$ en
imposant à sa partie horizontale d'être nulle. Par abus de langage, 
si on se donne par exemple un élément de $\mathcal G^*$, on dira 
qu'il <<~induit~>> une $1$-forme verticale sur $M$, en précisant
la connexion utilisée s'il y a ambiguïté.

 On va montrer le résultat suivant, qui permet de définir une
généralisation de la classe d'Euler aux fibré principaux en tore:
\begin{pr}\label{topo:pr}
Soit $\bar\omega\in\mathcal G^*$, $\omega$ la $1$-forme différentielle
sur $M$ induite par $\bar\omega$ et $\alpha_\omega$ la $2$-forme
différentielle sur $N$ telle que $\de\omega=\pi^*(\alpha_\omega)$. Alors
l'application $e:\mathcal G^*\rightarrow H^2(N,\R)$
donnée par $\bar\omega\mapsto[\alpha_\omega]$ est bien définie (c.-à-d. que
la classe de cohomologie de $\alpha_\omega$ ne dépend pas du choix
de la connexion) et linéaire.
\end{pr}
\textbf{Démonstration:}
On pose $T^k=\R^k/\Z^k$ et on note $(\bar\omega_i)_{i\in[1,k]}$ les formes
coordonnées de $\R^k$ passées au quotient sur $T^k$. En utilisant la
décomposition $T^k=\prod_{i=1}^kS_{(i)}^1$, on définit la famille
de fibrés en cercle $(M_i)$ comme en (\ref{topo:eq1}). On a alors des
projections $M\stackrel{\pi_i}{\rightarrow}M_i\stackrel{\pi_i'}{\rightarrow}N$
qui vérifient $\pi_i'\circ\pi_i=\pi_j'\circ\pi_j=\pi$.
Chaque forme $\omega_i$ induite sur $M$ par $\bar\omega_i$ est
le relevé $\pi_i^*(\omega_i')$ de la forme de connexion du fibré
en cercle $M_i$. On peut écrire $\de\omega_i=\pi_i^*(\de\omega_i')
=\pi^*(e_i)$ où $e_i$ est une $2$-forme sur $N$. Or, on sait que
$e_i$ appartient à la classe d'Euler du fibré $M_i$, indépendamment
du choix de la connexion sur $M_i$ donc du choix de la connexion sur
$M$. Si on définit $e:\mathcal G^*\rightarrow H^2(N,\R)$ par
$e(\bar\omega_i)=e_i$ en l'étendant par linéarité, on aura bien,
par linéarité de la différentielle extérieure, $\de\omega=\pi^*(e(\bar
\omega))$ pour tout $\bar\omega\in\mathcal G^*$, en notant $\omega$
la forme induite sur $M$.
\hfill\carrenoir

\begin{remarque} 
Si $k=1$ et si $\bar\omega$ est la forme volume du cercle de longueur $1$,
alors $e(\bar\omega)$ est la classe d'Euler du fibré.
\end{remarque}

\begin{remarque}
La démonstration \ref{topo:pr} met en évidence le lien entre l'invariant
$e$ et la famille de classe d'Euler associée à une décomposition
particulière en somme de Whitney de $M$: à chaque décomposition possible
est associée une base de $\mathcal G^*$, et la famille de classe d'Euler
est l'image de cette base par $e$.
\end{remarque}

\begin{exemple}
Au chapitre \ref{pvp:tore}, on a considéré des fibrés principaux en tore
$T^k$ non triviaux dont la base est un tore $T^2$. Comme $H^2(T^2)$ est 
de dimension $1$, le noyau de $e$ est de dimension $k-1$, ce qui signifie
qu'on peut décomposer le fibré en une somme de Whitney de $k$ fibrés
en cercles dont $k-1$ sont triviaux. On retrouve donc le fait que
le fibré peut s'écrire comme le produit d'un fibré en cercle et d'un
tore de dimension $k-1$.
\end{exemple}

 Dans la suite, si $\omega$ est une forme induite par un élément
$\bar\omega$ de $\mathcal G^*$, on écrira parfois par abus de langage
<<~$e(\omega)$~>> au lieu de <<~$e(\bar\omega)$~>>. De plus, verra 
parfois $e$ comme une application de $\mathcal G^*$ dans l'espace
$\mathcal H^2(N,h)$ des $2$-formes harmoniques de $N$, en utilisant
le fait que $\mathcal H^2(N,h)$ est canoniquement isomorphe à 
$H^2(N,\R)$.

\chapter{Formes invariantes et petites valeurs propres}\label{2001}

Nous allons ici démontrer les résultats \ref{inv:th1}, \ref{inv:th2} et 
\ref{inv:cor} énoncés dans l'introduction. Ceux-ci s'appuient 
essentiellement sur le
  
\begin{lem}\label{inv:lem}
Soit $(M,g)$ une variété riemannienne, $\varphi_t$ un flot agissant par
isométrie sur $M$ et $X$ le champ de vecteur associé. On suppose de plus
que $X$ n'est pas uniformément nul.

 Soit $E$ un sous-espace de $\Omega^*(M)$ stable par le laplacien
et par $(\varphi^*_t)_{t\in\R}$, 
$\lambda$ une valeur propre du laplacien restreint à $E$, et $E_\lambda$ 
l'espace propre associé. 

 S'il existe $T$ tel que $\varphi^*_{t+T}=\varphi^*_t$ pour tout $t\in\R$ 
et que $\lambda<\displaystyle\left(\frac{2\pi}{T\|X\|_\infty}\right)^2$,
alors $(\varphi^*_t)_{t\in\R}$ agit trivialement sur $E_\lambda$.

\end{lem}

\begin{demo}{du lemme \ref{inv:lem}} Remarquons tout d'abord qu'on peut 
supposer que $T=2\pi$, le résultat général se déduisant par simple 
changement de variable. D'autre part, par théorie de Hodge, on peut se 
restreindre à l'étude des formes cofermées. On supposera donc que
$\lambda$ est une valeur propre de $\codiff d_{|\Ker\codiff}$ et $E_\lambda$
désignera le sous-espace propre associé dans $\Ker\codiff$.

Soit $\omega\in E_\lambda$. On sait que $\mathcal L_X\omega=i_X
\circ\de\omega+\de\circ i_X\omega$.

D'une part, on a $\mathcal L_X\omega=\displaystyle\lim_{t\rightarrow0}
\frac{\varphi_t^*\omega-\omega}{t}$. Puisque $\varphi_t$ est une isométrie,
les formes $\varphi_t^*\omega$ et $\frac{\varphi_t^*\omega-\omega}{t}$
sont dans $\Ker\codiff$, et donc $\mathcal L_X\omega$ aussi car
$\Ker\codiff$ est fermé. D'autre part, $\de\circ i_X\omega$ est dans $\Ima\de$. 

Comme $\Ker\codiff$ et $\Ima\de$ sont orthogonaux, les formes
$\mathcal L_X\omega$ et $\de\circ i_X\omega$ sont orthogonales. Le
théorème de Pythagore donne donc :
\begin{equation}
\|\mathcal L_X\omega\|^2_2+\|\de\circ i_X\omega\|^2_2=\|i_X\circ\de
\omega\|^2_2.
\end{equation}
On en déduit
\begin{equation}\label{inv:eq1}
\|\mathcal L_X\omega\|^2_2\leq\|i_X\circ\de\omega\|^2_2
\leq\|i_X\|^2\|\de\omega\|^2_2\leq\|i_X\|^2\lambda\|\omega\|^2_2.
\end{equation}
On va maintenant majorer la norme de $i_X$ d'une part, et évaluer
celle de $\mathcal L_X\omega$ d'autre part.

\subsubsection*{majoration de $\|i_X\|$}

Soit $\alpha\in\Omega^p(M)$.

Soit $m\in M$, et $(X_1,\cdots,X_n)$ une base orthonormée de $T_mM$
telle que $X=\mu X_1$. On pose $\alpha_m=\displaystyle
\sum_{i_1<\cdots<i_p}\alpha_{i_1,\cdots,i_p}X_{i_1}^\flat
\wedge\cdots\wedge X_{i_p}^\flat$. On a alors 
\begin{equation}
i_X(\alpha_m)=
\displaystyle\sum_{i_1<\cdots<i_p}\alpha_{i_1,\cdots,i_p}i_X
(X_{i_1}^\flat\wedge\cdots\wedge X_{i_p}^\flat).
\end{equation}

Si $i_1=1$, alors $i_X(X_{i_1}^\flat\wedge\cdots\wedge X_{i_p}^\flat)
=\mu X_{i_2}^\flat\wedge\cdots\wedge X_{i_p}^\flat$.

Si $i_1\neq1$, alors $i_X(X_{i_1}^\flat\wedge\cdots\wedge X_{i_p}^\flat)=0$.

Donc $i_X(\alpha_m)=\mu \sum_{i_2<\cdots<i_p}\alpha_{1,i_2,\cdots,i_p}
X_{i_2}^\flat\wedge\cdots\wedge X_{i_p}^\flat$, et
$$
|i_X(\alpha_m)|^2=\mu^2\sum_{i_2<\cdots<i_p}\alpha_{1,i_2,\cdots,i_p}^2
\leq\mu^2\sum_{i_1<\cdots<i_p}\alpha_{i_1,\cdots,i_p}^2=|X|^2|\alpha_m|^2.
$$
On en déduit
\begin{eqnarray*}
\|i_X(\alpha)\|_2^2&=&\int_M|i_X(\alpha_m)|^2\de v\leq
\int_M|X|^2|\alpha_m|^2\de v\leq\|X\|_\infty\int_M|\alpha_m|^2\de v\\
&\leq& \|X\|_\infty\|\alpha\|^2_2.
\end{eqnarray*}
et donc 
\begin{equation}\label{inv:eq2}
\|i_X\|\leq\|X\|_\infty
\end{equation}

\subsubsection*{calcul de $\|\mathcal L_x\omega\|$}
 Comme $\varphi_t$ est une isométrie, $\varphi^*_t$ agit par isométrie sur
$E_\lambda$. Si on suppose que $t\mapsto\varphi_t$ est $2\pi$-périodique,
$\varphi_t$ induit donc un morphisme $S^1\rightarrow\SO(E_\lambda)$ qu'on
peut décomposer en somme de représentation irréductibles.

 Les représentations irréductibles de $S^1$ sont :
\begin{itemize}
\item la représentation triviale $S^1\rightarrow\SO(1),\ t\mapsto \Id$;
\item les rotations du plan $S^1\rightarrow\SO(2),\ t\mapsto
R(kt)=\smat{\cos kt & -\sin kt\\\sin kt & \cos kt}, k\in Z^*$.
\end{itemize}
 Supposons maintenant que la décomposition fasse apparaître au
moins une rotation. On peut choisir une forme $\omega\neq0$ située 
dans le sous-espace stable associé. On a alors
\begin{equation}
\|\mathcal L_X\omega\|=\left\|\lim_{t\rightarrow0}\frac{\varphi^*_t
\omega-\omega}t\right\|=|k|\cdot\|\omega\|\geq\|\omega\|.
\end{equation}
Avec (\ref{inv:eq1}) et (\ref{inv:eq2}), on en déduit :
\begin{equation}
\|\omega\|^2\leq\lambda\|X\|^2_\infty\|\omega\|^2,
\end{equation}
et finalement
\begin{equation}
\lambda\geq\displaystyle\frac1{\|X\|^2_\infty}.
\end{equation}
La conclusion du lemme en découle immédiatement.
\end{demo}

\begin{demo}{du théorème \ref{inv:th1}}
C'est une application directe du lemme \ref{inv:lem} avec $E=
\Omega^p(M)$, le flot $\varphi_t$ étant induit par l'action de $S^1$. 
Le champ $X$ est alors un champ vertical $S^1$-invariant.

Si on paramètre $\varphi_t$ de manière à être $2\pi$-périodique, 
la longueur $l$ d'une fibre sera égale à $2\pi|X|$, la norme de $X$
ne dépendant pas du point choisi sur la fibre. On a donc
\begin{equation}
|X|=\frac{l}{2\pi}
\end{equation}
et par conséquent
\begin{equation}
\|X\|_\infty=\frac{l_0}{2\pi}.
\end{equation}
Si $\lambda<\left(\frac{2\pi}{l_0}\right)^2$, alors $\lambda<\frac1
{\|X\|^2_\infty}$ et donc les formes propres de $E_\lambda$ sont 
$S^1$-invariantes.
\end{demo}

\begin{demo}{du théorème \ref{inv:th2}}
Remarquons tout d'abord que dans la démonstration du théorème 
\ref{inv:th1}, la décomposition
de $\Omega^p(M)$ en éléments irréductibles est une décomposition
en série de Fourier par rapport à la fibre $S^1$, la représentation
triviale et les représentations $\theta\mapsto R(k\theta)$ 
correspondant respectivement aux fonctions constantes et aux
fonctions $\frac{2\pi}k$-périodiques du cercle. On va reprendre cette 
idée et l'appliquer au tore $T^k$ pour décomposer $\Omega^p(M)$ en somme
en sous-espaces de formes invariantes dans une direction et
périodiques dans une autre pour ensuite appliquer le lemme \ref{inv:lem}
à ces sous-espaces.

Plus précisément, si on pose $T^k=\R^k/\Z^k$, $\R^k$ étant muni d'une métrique
euclidienne (pas nécessairement la métrique euclidienne canonique), 
et si $\Gamma$ est le réseau dual de $\Z^k$, $\Gamma=
\{\gamma\in\R^k,\langle\gamma,\gamma'\rangle\in\Z,\forall\gamma'\in\Z^k\}$,
une base des fonctions propres de $T^k$ est donnée par $f_\gamma=\cos(2\pi
\langle\gamma,x\rangle)$ et $g_\gamma=\sin(2\pi\langle\gamma,x\rangle)$,
$\gamma\in\Gamma$ (voir par exemple \cite{ghl}, p.~200).
Si on note $\gamma^\perp$ l'orthogonal de $\gamma$ dans
$\R^n$, les fonctions $f_\gamma$ et $g_\gamma$ sont invariantes
sous l'action de $\gamma^\perp$. On peut remarquer
que, si $\Gamma$ dépend de
la métrique choisie sur $\R^k$, ce n'est pas le cas de l'ensemble
des $\gamma^\perp_{\gamma\neq0}$. En effet, c'est l'ensemble des
hyperplans vectoriels de $\R^k$ engendrés par des éléments de $\Z^k$.
Notons $A$ cet ensemble. En regroupant les
fonctions propres en fonction de leur direction invariante, on
obtient la décomposition  
\begin{equation}C^\infty(T^k)=\overline{\displaystyle
\bigoplus_{V\in A} C^\infty(T^k)^V\oplus\R}
\end{equation}
 où 
$\R$ représente les fonctions constantes, et où $C^\infty(T^k)^V$ 
est l'espace des fonctions $C^\infty$ d'intégrale nulle 
(c'est-à-dire orthogonales aux fonctions constantes) et invariantes 
dans la direction $V$. Par construction,
chaque $C^\infty(T^k)^V$ est invariant par $\Delta$ et
par l'action de $T^k$. 

De la même façon, $\Omega^p(M)$ peut s'écrire
\begin{equation}
\Omega^p(M)=\overline{\displaystyle\bigoplus_{V\in A}\Omega^p(M)^V\oplus
\Omega^p(M)^{T^k}},
\end{equation}
où $\Omega^p(M)^{T^k}$ est l'espace des $p$-formes $T^k$-invariantes
et $\Omega^p(M)^V$ l'espace des $p$-formes
différentielles invariantes par l'action de $V$ et 
orthogonales à $\Omega^p(M)^{T^k}$.
Comme $V$ ne dépend pas de la métrique sur chaque fibre, chacun des 
$\Omega^p(M)^V$ est bien défini, et sera de plus 
invariant par $\Delta$ et par l'action de $T^k$. On a ainsi partiellement
décomposé les formes différentielles de $M$ en série de Fourier 
par rapport à la fibre.

Soit $V\in A$, $\lambda_M^V$ une valeur propre du laplacien
restreint à $\Omega^p(M)^V$ et $\lambda_{T^k}^V$ la première valeur
propre de $\Omega^p(T^k)^V$ (~remarque : elle est non
nulle car les formes harmoniques du tore plat sont les formes 
invariantes~). On choisit sur $T^k$ un champ
invariant $X_{T^k}$ orthogonal à $V$ et on note $X_M$ le champ
vertical sur $M$ induit par $X_{T^k}$. L'action
sur $\Omega^p(M)^V$ du flot $\varphi_t$ associé à $X_M$ est périodique,
et si on note $T$ sa période, $\lambda_{T^k}^V$ est par construction
exactement $\displaystyle\left(\frac{2\pi}{T\|X_{T^k}\|_\infty}\right)^2$.
D'autre part, comme $\bar g_x\leq f(x)\cdot\bar g$, 
les normes de $\|X_{T^k}\|_\infty$ et $\|X_{M}\|_\infty$ 
sont liés~: 
\begin{equation}
\|X_{M}\|_\infty\leq(\sup_B f)^{1/2}\cdot\|X_{T^k}\|_\infty,
\end{equation}
donc si $\lambda_M^V<(\sup_B f)^{-1}\lambda_{T^k}^V$, alors
\begin{equation}
\displaystyle\lambda_M^V<
(\sup_B f)^{-1}\left(\frac{2\pi}{T\|X_{T^k}\|_\infty}\right)^2\leq
\left(\frac{2\pi}{T\|X_{M}\|_\infty}\right)^2,
\end{equation}
et le lemme \ref{inv:lem} s'applique 
avec $E=\Omega^p(M)^V$
et les formes propres associées à $\lambda_M^V$ sont 
$\varphi_t^*$-invariantes, donc $T^k$-invariantes.

 Si $\lambda_M$ est une valeur propre du laplacien agissant sur
$\Omega^p(M)$ et que $\lambda$ est strictement inférieure 
à la première valeur propre de $T^k$, alors elle sera
\emph{a fortiori} inférieure à tous les $\lambda_{T^k}^V$ et
donc les formes propres de $\lambda_M^V$ dans $\Omega^p(M)^V$
sont $T^k$-invariantes. Comme
les $\Omega^p(M)^V$ sont stables par le laplacien, l'espace
propre de $\lambda_M$ est la somme des espaces propres restreints
aux $\Omega^p(M)^V$. Par conséquent, toutes les formes propres
associées à $\lambda_M$ sont $T^k$-invariantes.
\end{demo}

\begin{demo}{du corollaire \ref{inv:cor}}
 L'hypothèse sur la métrique peut s'écrire $\bar g_x= f(x)\cdot\bar g$,
où $f$ est une fonction positive sur $N$. On peut alors appliquer
le théorème \ref{inv:th2}. Il reste à montrer que si $\displaystyle
\lambda<\left(\frac\pi{d_0}\right)^2$ alors $\displaystyle
\lambda<(\sup_{x\in B}f(x))^{-1}\cdot\lambda_{0,1}(T^k,\bar g)$.

Comme la métrique resteinte à la fibre $\pi^{-1}(x)$ est
$f(x)\cdot\bar g$, la première valeur propre de laplacien restreint
à cette fibre est $\frac{\lambda_{0,1}(T^k,\bar g)}{f(x)}$.
De plus, la première valeur propre d'un tore plat de diamètre $d$
est minorée par $(\frac\pi d)^2$, par conséquent 
\begin{equation}
\frac{\lambda_{0,1}(T^k,\bar g)}{f(x)}=\lambda_{0,1}(T^k,f(x)\bar g)
\geq\left(\frac\pi {d_x}\right)^2,
\end{equation}
 où $d_x$ est le diamètre
de la fibre $\pi^{-1}(x)$ pour la distance intrinsèque, et donc
\begin{equation}
\displaystyle\frac{\lambda_{0,1}(T^k,\bar g)}{\sup_B f}
\geq\left(\frac\pi {d_0}\right)^2,
\end{equation}
ce qui achève la démonstration.
\end{demo}

L'exemple suivant montre que si on ne suppose pas que les fibres
sont homothétiques entre elles, une majoration du diamètre des
fibres ne permet pas de majorer la fonction $f$ du théorème \ref{inv:th2}.

\begin{exemple}
On considère sur $\R^2$ muni de son système de coordonnés canonique 
la famille de métrique $g_t=(\de x+t\de y)^2+\de y^2$. Ces métriques
passent au quotient sur le tore $T^2=\R^2/\Z^2$. Quel que soit $t_0\in\R$,
il n'existe pas de constante $c>0$ telle que $g_t\leq c\cdot g_{t_0}$ pour
tout $t$. Cependant, le diamètre de $(T^2,g_t)$ reste borné quand $t$
varie. En effet, le difféomorphisme linéaire $\smat{1&-1\\0&1}$ est
une isométrie de $(T^2,g_t)$ dans $(T^2,g_{t+1})$, et par conséquent 
le diamètre de $(T^2,g_t)$ est une fonction périodique de $t$.
\end{exemple}

\begin{demo}{de la remarque \ref{inv:rem}}
 Il suffit de remarquer que dans le lemme \ref{inv:lem}, 
si la dimension du sous-espace propre $E_\lambda$ est impaire, 
la décomposition de cet espace en espace de représentations
irréductibles contient nécessairement une représentation triviale, et
donc $E_\lambda$ contient des formes invariantes par $\varphi_t$. 
Dans le cas du théorème \ref{inv:th1}, les espaces propres de dimension
impaire du laplacien agissant sur $\Omega^*(M)$ contiennent donc des
formes $S^1$-invariantes.

 Dans le cas du théorème \ref{inv:th2}, si un espace propre du laplacien
est de dimension impaire, l'un des éléments de la décomposition
de Fourier de cet espace sera aussi de dimension impaire, et la remarque
précédente s'applique.
\end{demo}

\chapter{Géométrie des fibrés principaux en tores}\label{geom}
\section{Métriques adaptées}
Nous allons ici montrer qu'on peut, dans le but d'obtenir le
théorème \ref{vol:th}, se ramener à une situation géométrique pour
laquelle l'étude du spectre d'un fibré en tore est plus simple. 
Cette situation est une généralisation de la notion de métrique 
adaptée définie dans le cas des fibrés en cercle par B.~Colbois
et G.~Courtois (\cite{cc2}) :
\begin{df}
On dit que le couple de métriques $(g,h)$ définies sur $M$ et $N$ 
respectivement est adapté à la fibration principale
$T^k\hookrightarrow M^n\stackrel{\pi}{\rightarrow} N$ si:
\renewcommand{\theenumi}{\textsc{\roman{enumi}}}
\renewcommand{\labelenumi}{\theenumi.}
\begin{enumerate}\label{geom:df}
\item $\pi:(M,g)\rightarrow (N,h)$ est une submersion riemannienne ;
\item L'action de $T^k$ sur $M$ est isométrique ;
\item Les fibres sont totalement géodésiques ;
\item Toute $1$-forme verticale $\omega$ induite par un élément de 
$\mathcal G^*$ vérifie $\de\omega=\pi^*(e(\omega))$.
\end{enumerate}
\end{df}

On veut montrer qu'une métrique de courbure bornée
sur un fibré principal en tore est proche d'une métrique 
adaptée:
\begin{theo}\label{geom:th}
 Soient $a$ et $d$ deux réels strictement positifs, et
$T^k\hookrightarrow (M^n,g)\stackrel{\pi}{\rightarrow} (N,h)$
un fibré principal en tore. Il existe 
des constantes $\varepsilon_0(n,a,d,(N,h))>0$, $\tau(n,a,d,(N,h))>0$,
$\tau'(n,a,d,(N,h))>0$ et $c(n,a,d,(N,h))>0$
telles que si $|K(N,h)|\leq a$, $|K(M,g)|\leq a$, $\diam(M,g)\leq d$ et
si $\pi$ est une $\varepsilon$-approximation de Hausdorff avec $\varepsilon
<\varepsilon_0$, alors il 
existe des métriques $\tilde g$ et $\tilde h$  sur $M$ et $N$ respectivement 
et une fibration $\pi': (M,\tilde g)\rightarrow(N,\tilde h)$  telles que 
\begin{enumerate}
\item Le couple $(\tilde g,\tilde h)$ est adapté à la fibration $\pi'$;
\item $\displaystyle\frac 1\tau g\leq\tilde g\leq \tau g$ et
$\displaystyle\frac 1\tau h\leq\tilde h\leq \tau h$;
\item La restriction de $\tilde g$ à la fibre est telle que 
$\diam(\pi'^{-1}(x))\leq\tau'\varepsilon$, pour tout $x\in N$;
\item La courbure sectionnelle de $(M,\tilde g)$ vérifie
$|K(X,Y)|\leq c$, pour toute paire de vecteurs horizontaux orthonormés
$(X,Y)$.
\end{enumerate}
\end{theo}
On pourra alors appliquer le résultat de J.~Dodziuk selon
lequel si deux métriques sont proches, alors les spectres du laplacien
pour ces deux métriques sont proches aussi:
\begin{theo}[\cite{dz}]
Soit $g$ et $\tilde g$ deux métriques riemanniennes sur une variété
$M$ de dimension $n$, et $\tau$ une constante positive. Si les deux
métriques vérifient $\frac1\tau g\leq\tilde g\leq\tau g$, alors
$$\frac1{\tau^{3n-1}}\lambda_{p,k}(M,g)\leq\lambda_{p,k}(M,\tilde g)
\leq\tau^{3n-1}\lambda_{p,k}(M,g),$$
pour tout entiers $k\geq0$ et $p\in[0,n]$.
\end{theo}
\begin{remarque} Le théorème \ref{geom:th} implique en particulier
le théorème \ref{geom:th:intro}. La restriction sur la géométrie 
imposée par le point 4 de la conclusion du théorème \ref{geom:th} permet 
de mieux contrôler le spectre du laplacien. 
\end{remarque}

\begin{remarque} En vertu d'un théorème de Hermann (\cite{her},
\cite{bes} p.~249),
le fait que les fibres soient totalement géodésiques implique 
qu'elles sont isométriques entre elles. On va voir dans la démonstration du
théorème \ref{geom:th} que réciproquement, sur les fibrés considérés, si
la métrique est invariante et que les fibres sont isométriques alors 
elles sont totalement géodésiques.
\end{remarque}
\section{Situation de métrique invariante}
 Nous allons dans un premier temps montrer que si on suppose qu'on 
a sur $M$ une métrique invariante, elle est proche d'une métrique
qui vérifie les points (\textsc{i}) à (\textsc{iii}) de la définition 
\ref{geom:df}. Plus précisément :
\begin{pr}\label{geom:pr}
Soit $T^k\hookrightarrow (M^n,g)\stackrel{\pi}{\rightarrow} (N,h)$
un fibré principal en tore muni d'une métrique invariante $g$ tel que 
$\pi$ soit une submersion riemannienne. Pour tout $a>0$ et $d>0$,
il existe des constantes $\tau(n,a,d)>0$ et $c(n,a)>0$ telles
que si $|K(N,h)|\leq a$, $K(M,g)\geq -a$ et $\diam(M,g)\leq d$, 
alors il existe une métrique invariante $\tilde g$ sur $M$ telle 
que la fibration $\pi:(M,\tilde g)\rightarrow(N,h)$ soit une 
submersion riemannienne à fibres totalement géodésiques et
$$\frac 1\tau g\leq\tilde g\leq \tau g.$$
\end{pr}
\begin{remarque} On peut noter qu'on utilise non pas une hypothèse
de courbure bornée sur $M$ mais seulement que la courbure sectionnelle
est minorée.
\end{remarque}

Pour montrer la proposition \ref{geom:pr}, on utilisera les deux lemmes 
suivants. Le premier est une application directe de la formule de O'Neill:
\begin{lem}\label{geom:oneill}
Soit $a>0$ et $T^k\hookrightarrow (M^n,g) 
\stackrel{\pi}{\rightarrow} (N,h)$ un fibré
principal en tore muni d'une métrique invariante $g$ tel que $p$
soit une submersion riemannienne, $|K(N,h)|\leq a$,
et $K_{(M,g)}(X,Y)\geq -a$ pour tout couple $(X,Y)$ de vecteurs 
horizontaux orthonormés.
Alors, pour toute $1$-forme différentielle 
verticale $\omega$ induite par l'action de $T^k$, on a :
\begin{enumerate}
\item $\displaystyle|\de\omega(X,Y)|_x^2\leq\frac{8a}3|\omega|_x^2$, pour tout 
$x\in M$ et tout couple de vecteurs horizontaux orthonormés $X$ et $Y$ ;
\item $\displaystyle\|\de\omega\|_\infty\leq\frac{4an(n-1)}3\|\omega\|_\infty$. 
\end{enumerate}
\end{lem}
\textbf{Démonstration:} Soit $x\in M$, $y=\pi(x)$, $\tilde X$ et $\tilde Y$ 
deux champs de $N$ orthonormés en $y$, et $X$ et $Y$ les relevés de
$\tilde X$ et $\tilde Y$ à $M$. La formule de O'Neill (\cite{ghl}
p.~127, \cite{bes} p.~241) donne
\begin{equation}
K_N(\tilde X,\tilde Y)=K_M(X,Y)+\frac34\left|[X,Y]^V\right|^2,
\end{equation}
où $[X,Y]^V$ désigne la composante verticale de $[X,Y]$.
D'autre part on a, en utilisant le fait que $\omega$ est verticale, 
\begin{eqnarray}
\de \omega(X,Y)&=&X\cdot\omega(Y)-Y\cdot\omega(X)-\omega([X,Y])\nonumber\\
&=&-\omega([X,Y]).
\end{eqnarray}
On en déduit :
\begin{eqnarray}
|\de \omega(X,Y)|_x^2&=&|\omega([X,Y])|_x^2\leq|\omega|_x^2|[X,Y]^V|_x^2
\nonumber\\
&\leq&\frac43|\omega|_x^2(K_y(\tilde X,\tilde Y)-K_x(X,Y)).
\end{eqnarray}
Comme chacun des couples $(\tilde X,\tilde Y)$ et $(X,Y)$ est orthonormé
en $x$ et $y$, on a les majorations $|K_y(\tilde X,\tilde Y)|\leq a$ 
et $K_x(X,Y)\geq -a$, et donc
\begin{equation}
|\de \omega(X,Y)|_x^2\leq\frac{8a}3|\omega|_x^2.
\end{equation}
 Et comme l'inégalité précédente est vraie quel que soit le choix
de $(\tilde X,\tilde Y)$, il en découle finalement
\begin{equation}
|\de\omega|_x^2\leq\frac{4an(n-1)}3|\omega|_x^2
\leq\frac{4an(n-1)}3\|\omega\|^2,
\end{equation}
ce qui achève la démonstration.
\hfill\carrenoir

Le second lemme montre, dans le cas d'un fibré en cercle, qu'à courbure bornée, 
la longueur des fibres varie peu d'une fibre à l'autre. 
\begin{lem}\label{geom:lem2}
Soit $S^1\hookrightarrow (M^n,g)\stackrel{\pi}{\rightarrow} (N,h)$
un fibré principal en cercle sur $N$, tel que $g$ soit invariante
et $\pi$ soit une submersion riemannienne. Pour tout $a>0$ et $d>0$, 
il existe $\tau(n,a,d)>0$ tel que si $|K(N,h)|\leq a$, $K(M,g)\geq -a$ 
et $\diam(M,g)\leq d$, alors pour tout $x,y\in N$, on a
$$\frac 1\tau l_y\leq l_x\leq \tau l_y,$$
où $l_x$ et $l_y$ désignent les longueurs des fibres au dessus de 
$\pi^{-1}(x)$ et $\pi^{-1}(y)$ respectivement.
\end{lem}
\textbf{Démonstration:} On choisit sur le fibré $M$ une $1$-forme 
verticale $\omega$ dont l'intégrale sur chaque fibre est égale à $1$.
Soit $U$ le champ vertical induit par l'action de $S^1$ qui vérifie
$\omega(U)=1$. La norme $|U|$ de ce champ est
constante sur chaque fibre, et s'écrit $|U|=\pi^*f$, où $f$ est une
fonction sur $N$. De plus, en tout point $x$ de $N$, la norme
$f(x)$ de $U$ est égale à la longueur de la fibre au dessus de $x$.
On va montrer que $f$ est bornée en fonction de $a$ 
et $d$. Remarque : $\omega$ n'est pas la forme duale de $U$ pour
la métrique. Sa norme ponctuelle sur la fibre $\pi^{-1}(x)$
est $|\omega|=\frac1f$, et on a $U^\flat=f^2\omega$.

Soit $x\in N$, et $\tilde X$ un vecteur unitaire tangent à $N$ en $x$.
Soit $\tilde X_i$ une base orthonormée de champs de vecteurs au voisinage 
de $x$, telle que $\LC_{\tilde X_1}\tilde X_1=0$ sur ce voisinage, et 
$\tilde X_{1|x}=\tilde X$.
On relève cette base à $T^HM$ en notant $X_i=\pi^*(\tilde X_i)$ et
$X=X_1$.  Ces champs vérifient 
\begin{equation}[X_i,U]=0.\end{equation}
En effet, ces crochets de Lie sont déterminés par $[X_i,U]=
\frac\de{\de t}(\Phi_t)_{*}X_i$, où $\Phi_t$ est le flot induit par
le champ $U$. Par définition de $U$, ce flot correspond à l'action de
$S^1$ sur $M$. Les crochets de Lie sont donc nuls, car les champs $X_i$ 
sont $S^1$-invariants. On notera par ailleurs $U'$ le champ de norme $1$
défini par $U'=U/|U|$.

On va calculer la courbure sectionnelle $K(X,U)$ en fonction de 
$f$ et de ses variations (Remarque: le champ $U$ n'est pas normé, 
mais ce calcul est plus simple que si on utilise le champ $U'$). 
Cette courbure s'écrit
\begin{eqnarray}
K(X,U)&=&\langle R(X,U)X,U\rangle\nonumber\\
&=&\langle\LC_U\LC_XX-\LC_X\LC_UX-\LC_{[X,U]}X,U\rangle,
\end{eqnarray}
où $\langle\cdot,\cdot\rangle$ désigne la métrique. Le vecteur $\LC_XX$ 
est horizontal et vaut $\pi^*(\LC_{\tilde X_1}\tilde X_1)$ 
(voir \cite{bes} p.~239), et par conséquent $\LC_XX=0$ au voisinage de $x$. 
Comme d'autre part $[X,U]=0$, on est ramené à calculer
\begin{equation}
K(X,U)=-\langle\LC_X\LC_UX,U\rangle.
\end{equation}
Pour ce faire, on utilisera la formule suivante, qui caractérise la
connexion de Levi-Civita:
\begin{eqnarray}\label{geom:lc}
2\langle\LC_{Z_1}Z_2,Z_3\rangle\!\!\!&=&\!\!\!Z_1\cdot\langle Z_2,Z_3\rangle
+Z_2\cdot\langle Z_3,Z_1\rangle-Z_3\cdot\langle Z_1,Z_2\rangle\nonumber\\
&&\!\!\!\!+\langle[Z_1,Z_2],Z_3\rangle-\langle[Z_1,Z_3],Z_2\rangle
-\langle[Z_2,Z_3],Z_1\rangle.\hspace{3mm}
\end{eqnarray}
En utilisant l'orthogonalité de $(X_1,\cdots,X_n,U)$ et le fait
que $[X_i,U]$=0, on obtient
$$2\langle\LC_UX,U\rangle=X\cdot\langle U,U\rangle=X\cdot f^2=2f\de f(X)$$
et $$2\langle\LC_UX,X_i\rangle=-\langle[X,X_i],U\rangle,$$
et donc
\begin{equation}\label{geom:DUX}
\LC_UX=\de f(X) U'-\frac12\sum_{i=1}^n\langle[X,X_i],U\rangle X_i.
\end{equation}
Le premier terme peut s'écrire $\de f(X) U'=\frac{\de f(X)}f U$, par
conséquent
\begin{eqnarray}
\LC_X(\de f(X) U')&=&\left(X\cdot\frac{\de f(X)}f\right)U
+\frac{\de f(X)}f\LC_XU\nonumber\\
&=&\left(\frac{(\LC_X\de f)(X)}f+\frac{\de f(\LC_XX)}f
-\frac{(\de f(X))^2}{f^2}\right)U\nonumber\\
&&+\frac{\de f(X)}f\LC_XU\\
&=&\frac{\Hess f(X,X)}fU-\frac{(\de f(X))^2}{f^2}U+\frac{\de f(X)}f\LC_XU,
\nonumber
\end{eqnarray}
en utilisant le fait que $\LC_XX=0$. De plus, la relation (\ref{geom:lc})
donne $$2\langle\LC_XU,U\rangle=X\cdot|U|^2=2f\de f(X),$$ et donc
\begin{equation}\label{geom:K1}
\langle\LC_X(\de f(X) U'),U\rangle=f\Hess f(X,X).
\end{equation}

La dérivation des termes suivants de (\ref{geom:DUX}) donne 
$$\LC_X(\langle[X,X_i],U\rangle X_i)=(X\cdot\langle[X,X_i],U\rangle)X_i
+\langle[X,X_i],U\rangle\LC_XX_i.$$
Quand on calcule le produit scalaire de cette expression avec $U$, le
premier terme s'annule, et comme la relation (\ref{geom:lc}) donne
$\langle\LC_XX_i,U\rangle=\langle[X,X_i],U\rangle$, il reste
\begin{equation}\label{geom:K2}
\langle\LC_X(\langle[X,X_i],U\rangle X_i),U\rangle=\langle[X,X_i],U\rangle^2.
\end{equation}
La somme des équations (\ref{geom:K1}) et (\ref{geom:K2}) donne
\begin{equation}
K(X,U)=-f\Hess f(X,X)+\frac12\sum_{i=1}^n\langle[X,X_i],U\rangle^2.
\end{equation}
En normalisant le vecteur $U$, on obtient
\begin{equation}
K(X,U')=-\frac{\Hess f(X,X)}f+\frac1{2f^2}\sum_{i=1}^n\langle[X,X_i],U
\rangle^2.
\end{equation}
Les derniers termes peuvent être majorés en fonction de la courbure. En
effet, on a
\begin{equation}
\langle[X,X_i],U\rangle=U^\flat([X,X_i])=f^2\omega([X,X_i])
=-f^2\de\omega(X,X_i),
\end{equation}
et donc, en vertu du lemme \ref{geom:oneill},
\begin{equation}
\langle[X,X_i],U\rangle^2\leq f^4\frac{8a}3|\omega|^2=\frac{8a}3 f^2.
\end{equation}
Par hypothèse, la courbure sectionnelle de $M$ est minorée par $-a$.
On a donc finalement :
\begin{equation}\label{geom:hess}
\frac{\Hess f(X,X)}f\leq \left(\frac{4n}3+1\right)a.
\end{equation}

Soient $x$ et $y$ deux points de $N$, et $\gamma$ une géodésique
minimisante joignant ces deux points. Notons $\mu$ la fonction
définie par 
\begin{equation}\mu(t)=\ln f\circ\gamma(t).\end{equation}
En dérivant $\mu$ par rapport à $t$, on obtient :
\begin{equation}
\mu'(t)=\frac{\de f(\gamma'(t))}{f\circ\gamma(t)}
\end{equation}
et
\begin{eqnarray}
\mu''(t)&=&-\left(\frac{\de f(\gamma'(t))}{f\circ\gamma(t)}\right)^2
+\frac{\de f(\LC_{\gamma'(t)}\gamma'(t))+\LC\de f(\gamma'(t),\gamma'(t))}
{f\circ\gamma(t)}\nonumber\\
&=&-\mu'(t)^2+\frac{\Hess f(\gamma'(t),\gamma'(t))}{f\circ\gamma(t)}.
\end{eqnarray}
On a donc, en vertu de la majoration (\ref{geom:hess}) :
\begin{equation}
\mu''(t)\leq\frac{\Hess f(\gamma'(t),\gamma'(t))}{f\circ\gamma(t)}
\leq\left(\frac{4n}3+1\right)a.
\end{equation}
Supposons que $x$ est un point où $f$, et donc $\mu$, atteint son 
minimum. On a alors $\mu'(0)=0$ et donc, en remarquant que $d$ majore
le diamètre de $N$,
\begin{equation}
\mu'(t)\leq\left(\frac{4n}3+1\right)at,
\end{equation}
et
\begin{equation}
\mu(t)\leq\frac{4n+3}6at^2\leq\frac{4n+3}6ad^2.
\end{equation}
Le rapport $\frac{f(y)}{f(x)}$ est donc majoré par une constante
$\tau=\exp(\frac{4n+3}6ad^2)$. Comme on a montré cette majoration 
en prenant pour $x$ un point où $f$ atteint son  minimum, elle
sera vraie \emph{a fortiori} pour un $x$ quelconque.

 Remarque : si les fibre sont isométriques, alors le champ $U$ est de
norme constante. Il est aisé de vérifier à l'aide de (\ref{geom:lc})
que $\LC_U U$ est alors nul, c'est-à-dire que les fibres sont totalement
géodésiques.
\hfill\carrenoir

\textbf{Démonstration de la proposition \ref{geom:pr}:}
 Le but est en fait de généraliser le lemme \ref{geom:lem2}
aux fibrés en tore pour montrer que $g$ est proche d'une métrique
pour laquelle toutes les fibres sont isométriques.

Soit $\bar U\in\mathcal G$ non nul et $U$ le champ vertical induit par 
$\bar U$ sur $M$. Soit $x_0\in N$. On choisit $\bar U$ de sorte que 
$|U|=1$ au dessus de $x_0$. De plus, on impose à $\bar U$ d'avoir un
coefficient directeur rationnel, c'est-à-dire que $\bar U$ est colinéaire
à un vecteur de $\Z^k\subset\mathcal G$. L'action du flot associé à
$U$ induit alors une fibration $S^1\hookrightarrow M
\stackrel{\pi}{\rightarrow}(M',g')$. On peut contrôler la courbure
de ce fibré. En effet, utilisant la formule de O'Neill, on peut
écrire
\begin{eqnarray}
K_{M'}(\tilde X,\tilde Y)&=&K_M(X,Y)-\frac34|[X,Y]^V|^2\nonumber\\
&=&K_M(X,Y)-\frac34\frac{|\de\omega(X,Y)|^2}{|\omega|^2},
\end{eqnarray}
où $\tilde X$ et $\tilde Y$ sont deux vecteurs orthonormés de $M'$,
$X$ et $Y$ leurs relevés respectifs sur $M$, et $\omega$ la $1$-forme
induite par l'action de $T^k$ telle que 
$\omega(U)=1$. Le lemme \ref{geom:oneill} permet
de contrôler le dernier terme en fonction de la courbure de $M$, et donc
\begin{equation}
K_{M'}(\tilde X,\tilde Y)\geq -a-2a=-3a.
\end{equation}
Le lemme \ref{geom:lem2} assure alors qu'il existe une constante 
$\tau(n,a,d)$ telle que 
\begin{equation}\label{geom:eq1}
\frac1\tau\leq|U|\leq\tau.
\end{equation}
Notons $\tilde g$ la métrique invariante sur $M$ obtenue en modifiant 
$g$ dans la direction verticale de sorte que les fibres soient isométriques
à $\pi^{-1}(x_0)$, et en conservant la distribution horizontale et la
métrique horizontale associées à $g$. Pour cette nouvelle métrique, la
norme de $U$ est uniformément égale à 1.
La relation (\ref{geom:eq1}) peut s'écrire
\begin{equation}\label{geom:eq2}
\frac1{\tau^2}\tilde g(U,U)\leq g(U,U)\leq\tau^2\tilde g(U,U),
\end{equation}
Par continuité, (\ref{geom:eq2}) s'étend à n'importe quel vecteur vertical.
Comme $g$ et $\tilde g$ sont identiques sur la direction horizontale, 
on aura finalement
\begin{equation}
\frac1{\tau^2}\tilde g\leq g\leq\tau^2\tilde g.
\end{equation}

Pour conclure, remarquons que si la métrique sur $M$ est telle que les 
fibres soient isométriques, la fibration $S^1\hookrightarrow M
\stackrel{\pi}{\rightarrow}(M',g')$ induite par le champ $U$ est à fibre 
totalement géodésique. Par continuité, les fibres du fibré $T^k
\hookrightarrow M \rightarrow N$ sont elles aussi totalement géodésiques.
\hfill\carrenoir

\section{Cas général}

 On va maintenant démontrer le théorème \ref{geom:th}. Pour ce faire, 
on va s'inspirer d'une démonstration d'un théorème de Lott 
(\cite{lo}, théorème 2), qui utilise les résultats de \cite{cfg}.

 Soit $g$ une métrique sur $M$ vérifiant les hypothèse du théorème 
\ref{geom:th}. Tout d'abord, en utilisant un résultat de régularisation
d'Abresch (\cite{cfg}, théorème 1.12), on construit une métrique
$g_1$ sur $M$ telle que $\frac1{\tau_1}g\leq g_1\leq\tau_1g$,
$|K(M,g_1)|\leq a$ et 
$\|\LC^iR\|\leq A_i(n,a,\tau_1)$, où $\tau_1>1$ est un réel fixé, et
$\LC$ et $R$ désignent respectivement la dérivée covariante et
le tenseur de courbure pour la métrique $g_1$.

 On applique ensuite le théorème 2.6 de \cite{cfg}, qui assure l'existence
de constantes $\epsilon_0(n,(N,h))$, $\kappa(n,A)$, $\kappa'(n,A,(N,h))$
et $\kappa_i(n,A,(N,h))$ et d'une fibration $\pi':(M,g_1)\rightarrow
(N,h)$ tels que si $\pi$ est une $\varepsilon$-approximation de 
Hausdorff avec $\varepsilon<\epsilon_0(n,(N,h))$, alors :

\begin{itemize}
\item pour tout $x\in N$, le diamètre de $\pi'^{-1}(x)$ pour 
la métrique $g'$ est inférieur à $\kappa\varepsilon$;
\item la seconde forme fondamentale de la fibre vérifie 
$\|I\!I_{\pi^{-1}(x)}\|_\infty\leq\kappa'$ pour tout $x\in N$;
\item la submersion $\pi'$ est $\kappa_i$-régulière, c'est-à-dire que
$\|\LC^i\pi'\|_\infty\leq \kappa_i$, pour tout $i\in\N$.
\end{itemize}
Enfin, pour une telle fibration $\pi'$, les parties 3 et 4 de 
\cite{cfg} donnent la construction d'un métrique $g_2$ sur $M$ qui
est $T^k$-invariante et telle que $|\LC^i(g_2-g_1)|\leq c(n,A,(N,h),i)$,
pour tout $i\in\N$. Cette dernière égalité assure l'existence d'une
constante $\tau_2(n,A,d,(N,h))$ telle que 
$\frac1{\tau_2}g_1\leq g_2\leq\tau_2g_1$, et permet aussi de contrôler la
courbure pour la métrique $g_2$.

On peut alors appliquer la proposition \ref{geom:pr}, qui nous
donne une métrique $g_3$ qui vérifie les points (\textsc{i}) à (\textsc{iii})
de la définition \ref{geom:df}.

Pour obtenir la métrique $\tilde g$ du théorème \ref{geom:th}, il reste
à modifier la distribution horizontale de manière à ce que 
$\tilde g$ vérifie le point (\textsc{iv}) de la définition \ref{geom:df}.
Remarquons tout d'abord que comme l'application $e:\mathcal G^*\rightarrow
\mathcal H^2(N)$ est linéaire, il suffit de montrer l'égalité
$\de\omega=\pi'^*(e(\omega))$ pour une base de $\mathcal G^*$.
Soit $(\omega_i)$ une base de $\mathcal G^*$ orthonormée pour la
métrique $g_3$. Pour chaque $i$, $\de\omega_i$ s'écrit
\begin{equation}
\de\omega_i=\pi'^*(\alpha_i+\de\beta_i),
\end{equation}
où $\alpha_i$ est une forme harmonique et $\beta_i$ une forme cofermée.
On définit une nouvelle forme verticale $\omega_i'=\omega_i-\pi'^*(\beta_i)$.
Cette forme vérifie 
\begin{equation}\label{geom:eq3}
\de\omega_i'=\de\omega_i-\pi'^*(\de\beta_i)=\pi'^*(\alpha_i)\in\mathcal H^2(N).
\end{equation}
L'intersection des noyaux des formes $\omega_i'$ définit une nouvelle
distribution horizontale. On définit $\tilde g$ comme étant la métrique
sur $M$ telle que $\pi':(M,\tilde g)\rightarrow(N,h)$ soit une
submersion riemannienne et $\tilde g=g_3$ sur l'espace vertical.
Cette métrique vérifie le point (iv) de la définition \ref{geom:df}
du fait de (\ref{geom:eq3}).

On doit encore vérifier $\tilde g$ est proche de $g_3$. Remarquons
que 
\begin{equation}
\tilde g-g_3=\sum_i(\omega'^2-\omega^2)=\sum_i(2\pi'^*(\beta_i)
\otimes\omega_i+\pi'^*(\beta_i)^2).
\end{equation}
Or, B.~Colbois et G.~Courtois ont montré dans \cite{cc2} (lemme
A.32) qu'il existe
une constante $\kappa(n,a,d,(N,h))>0$ telle que les formes $\beta_i$ telles
qu'on les a définies vérifient $\|\beta_i\|_\infty\leq\kappa$, ce qui permet de
conclure. Remarque: le lemme A.32 de \cite{cc2} utilise le fait
que pour la métrique $g_3$, la norme de la seconde forme fondamentale 
est contrôlé et que la submersion $\pi'$ est $\kappa_i$-régulière. Il n'est
donc pas évident qu'on puisse obtenir le théorème \ref{geom:th} en supposant
que la métrique initiale $g$ est invariante et en se passant des 
résultats de \cite{cfg}.

Enfin, il reste à montrer que la courbure de $(M,\tilde g)$ reste
bornée dans la direction horizontale. 
Soit $x\in N$, $\tilde X$ et $\tilde Y$ deux vecteurs orthonormés
tangents à $N$ en $x$, $\bar\omega$ une 
$1$-forme invariante de $T^k$, $\omega$ la $1$-forme induite sur $M$
pour la distribution horizontale associée à $g$ et $\omega'$
la $1$-forme induite pour la distribution associée à $\tilde g$.
Ces deux formes vérifient $\de\omega'=\pi'^*(\alpha)$ et
$\de\omega=\pi'^*(\alpha+\de\beta)$, où $\alpha$ est une $2$-forme
harmonique de $N$ et $\beta$ une $1$-forme de $N$.

D'après la formule de O'Neill, il suffit pour contrôler la courbure
sectionnelle $K_{(M,\tilde g)}(\pi'^*(X),\pi'^*(Y))$ de majorer
la norme de $[\pi'^*(X),\pi'^*(Y)]^V$. Or, on peut écrire d'une part,
\begin{equation}
\omega'([\pi'^*(X),\pi'^*(Y)]^V)=\de\omega'(\pi'^*(X),\pi'^*(Y))=\alpha(X,Y).
\end{equation}
D'autre part, on a
\begin{equation}\label{geom:li}
\|\alpha\|_\infty\leq \tau'(n,a,d)\|\alpha\|_2,
\end{equation}
d'après \cite{li}, car $\alpha$ est harmonique, et
\begin{equation}
\|\alpha\|_2\leq\|\alpha+\de\beta\|_2=\|\de\omega\|_2\leq
\|\de\omega\|_\infty,
\end{equation}
en utilisant le fait qu'une forme harmonique est le plus petit
élément de sa classe de cohomologie pour la norme $L^2$.
Enfin, le lemme \ref{geom:oneill} permet de contrôler
la norme de $\de\omega$ en fonction de $a$ et $\|\omega\|_\infty$,
et la norme de $\omega$ est contrôlé en fonction de $\|\omega'\|_\infty$.
Comme la majoration
de $\omega'([\pi'^*(X),\pi'^*(Y)]^V)$ obtenue est indépendante du choix
de $\bar\omega$, on a bien une majoration de $|[\pi'^*(X),\pi'^*(Y)]^V|$
en fonction de $n$, $a$ et $d$.
\hfill\carrenoir

\chapter{Petites valeurs propres des fibrés principaux en tores}
\label{adap}
\section{Minoration du spectre des $1$-formes par le volume du fibré}
 Les résultats des chapitres précédents nous permettent maintenant
de démontrer le théorème \ref{vol:th}. On a vu qu'on pouvait
se ramener au cas d'un fibré muni d'une métrique adaptée. On va donc
montrer le résultat du théorème \ref{vol:th} pour un fibré vérifiant
les conclusions du théorème \ref{geom:th}:
\begin{theo}\label{adap:th}
Soient $a>0$, $d>0$ deux réels, $n$, $k$ et $m$ trois entiers
tels que $n=k+m$, et $(N^m,h)$ une variété riemannienne. 
Il existe des constantes 
$c(n,a,d,(N,h))$ 
et $\varepsilon(n,a,d,(N,h))$ strictement
positives telles que si $\bar g$ est une métrique sur le tore 
$T^k$  telle que $\diam(T^k)<\varepsilon$ et si 
$T^k\hookrightarrow M^n\rightarrow N$ est un fibré principal
muni d'un couple de métriques $(g,h)$ adapté au fibré et tel que 
$g=\bar g$ en restriction à la fibre, 
$\diam(M,g)<d$ et $|K_M(X,Y)|\leq a$ pour toute
paire $(X,Y)$ de vecteurs horizontaux orthonormés, alors on a
$$\lambda_{1,1}(M,g)\geq c\cdot \Vol^2(T^k).$$
\end{theo}
On s'est ici donné comme hypothèse que la métrique sur $N$ est fixée.
En effet, une hypothèse sur la courbure ne nous sera pas suffisante.
On verra au paragraphe suivant dans quelle mesure on peut espérer
obtenir le même résultat avec des hypothèses plus faibles.

\textbf{Démonstration :}

Dans un premier temps, nous allons démontrer le théorème dans
le cas où le fibré $M$ ne contient pas de sous-fibré trivial,
c'est-à-dire quand  l'application $e:\mathcal G^*\rightarrow\mathcal H^2
(N,h)$ est injective. Nous généraliserons ensuite le résultat 
à un fibré principal quelconque. D'autre part, on se restreindra 
aux formes $T^k$-invariantes, en vertu des résultats du chapitre 
\ref{2001} (corollaire \ref{inv:cor}). En effet, le spectre des formes 
orthogonales aux formes invariantes sera minoré en fonction de 
la constante $\varepsilon$ du théorème, et on pourra toujours
choisir cette constante suffisamment petite de sorte que le
spectre des formes orthogonales aux formes invariantes soit
plus grand que le terme $c\cdot \Vol^2(T^k)$.

 Supposons donc $e$ injective. La démonstration se déroule en
deux étapes. D'abord, on se ramène à l'étude des valeurs
propres de l'opérateur $e^*e$, l'adjoint étant défini en munissant
$\mathcal H^2(N)$ de sa norme $L^2$ :

\begin{fait}\label{adap:ft1}
Il existe $\varepsilon(n,a,\lambda_{0,1}(N,h),\lambda_{1,1}(N,h))>0$ et 
$c(n,a,\lambda_{0,1}(N,h))>0$ tel que
pour toute $1$-forme $\varphi$ sur $M$ $T^k$-invariante et
orthogonale à $\Ker\Delta^1(M,g)$,
si le quotient de  Rayleigh de $\varphi$ vérifie 
$R(\varphi)<\varepsilon$, alors il existe une forme $\omega$
induite par un élément de $\mathcal G^*$ telle que 
$\|e(\omega)\|^2\leq c\cdot\varepsilon\|\omega\|^2$.
\end{fait}

\textbf{Démonstration :}
Soit $\varphi$ une $1$-forme différentielle $T^k$-invariante 
de $M$. On peut écrire 
\begin{equation}\label{adap:phi}
\varphi=\pi^*(\alpha)+\sum_{i=1}^k \pi^*(a_i)\cdot\omega_i,
\end{equation}
où $\alpha$ est une $1$-forme de $N$, $a_i$ des fonctions de $N$ et
$\omega_i$ les $1$-formes verticales induites par une base
orthonormée de $\mathcal G^*$. On a alors~:
\begin{equation}\label{adap:dphi}
\de\varphi=\pi^*(\de\alpha+\sum_{i=1}^k a_i\cdot e_i)+
\sum_{i=1}^k\de \pi^*(a_i)\wedge\omega_i,
\end{equation}
où $e_i$ désigne l'image de $\omega_i$ par l'application
$e:\mathcal G^*\rightarrow\mathcal H^2(N)$.
De plus, pour tout $i$ on a 
$$\|\codiff(\pi^*(a_i)\omega_i)\|^2=
(\codiff(\pi^*(a_i)\omega_i),\codiff(\pi^*(a_i)\omega_i))=
(\pi^*(a_i)\omega_i,\de\codiff(\pi^*(a_i)\omega_i)),$$
où $(\cdot,\cdot)$ désigne le produit scalaire $L^2$.
Comme $\pi^*(a_i)\omega_i$ est une forme $T^k$-invariante, 
$\codiff(\pi^*(a_i)\omega_i)$ est une fonction invariante, c'est-à-dire 
que c'est le relevé d'une fonction sur $N$. Par conséquent, 
$\de\codiff(\pi^*(a_i)\omega_i)$ est le relevé d'une $1$-forme
sur $N$, et est donc orthogonale à $\pi^*(a_i)\omega_i$. Finalement, on a:
\begin{equation}\label{adap:deltaphi}
\codiff\varphi=\pi^*(\codiff\alpha).
\end{equation}
On peut calculer précisément quelles sont les $1$-formes harmoniques
de $M$. On sait déjà que les formes harmoniques sont invariantes, donc de 
la forme donnée en (\ref{adap:phi}). Si $\varphi$ est harmonique, on a
de plus $\de\varphi=0$ et $\codiff\varphi=0$, donc 
\begin{equation}
\codiff\alpha=0,\ \de a_i=0\text{ pour tout $i$, et } 
\de\alpha+\sum_{i=1}^k a_i\cdot e_i=0.
\end{equation}
Comme les fonctions $a_i$ sont constantes, $\sum_{i=1}^k a_i\cdot e_i$ est
une $2$-forme harmonique de $N$, donc orthogonale à la forme exacte 
$\de\alpha$. On a donc $\Delta\alpha=0$ et $\sum_{i=1}^k a_i\cdot e_i=0$.
Comme $e$ est injective, les $e_i$ forment une famille libre, et donc
$a_i=0$ pour tout $i$. On obtient finalement que les formes harmoniques
de $M$ sont les relevés des formes harmoniques de $N$.

Supposons que $\varphi$ est de norme $1$, c'est-à-dire que 
$\|\alpha\|^2+\sum_{i=1}^k\|a_i\|^2=1$. Le quotient de 
Rayleigh de $\varphi$ s'écrit alors
\begin{equation}\label{adap:eq1}
R(\varphi)=\|\codiff\alpha\|^2+\|\de\alpha+\sum_{i=1}^k a_i\cdot e_i\|^2
+\sum_{i=1}^k\|\de \pi^*(a_i)\|^2.
\end{equation}
 Supposons que $R(\varphi)<\varepsilon$ pour un $\varepsilon>0$ donné.
Comme $\|\alpha\|^2+\sum_{i=1}^k\|a_i\|^2$ est égal à $1$, l'un
des termes de la somme est plus grande que $\frac1{1+k}$. Nous
allons distinguer les cas $\|\alpha\|^2>\frac1{1+k}$, et 
$\|a_p\|^2>\frac1{1+k}$ pour un $p\in[1,k]$.

Supposons $\|\alpha\|^2\geq\frac1{1+k}$ :

La forme $\varphi$ est orthogonale à $\Ker\Delta^1$, c'est-à-dire à
l'ensemble des relevés de formes harmoniques de $N$. Par conséquent,
la forme $\alpha$ est elle-même orthogonales aux formes harmoniques
de $N$, et donc $\|\de\alpha\|^2+\|\codiff\alpha\|^2\geq\lambda_{1,1}(N,h)$.

Si $\|\codiff\alpha\|^2\geq\|\de \alpha\|^2$, alors $\frac{\|\codiff\alpha\|^2}
{\|\alpha\|^2}\geq\frac{\lambda_{1,1}(N,h)}2$. Or, on a les
deux inégalités
\begin{equation}
\|\alpha\|^2\geq\frac1{k+1}
\end{equation}
et 
\begin{equation}
\|\codiff\alpha\|^2\leq R(\varphi)<\varepsilon,
\end{equation}
dont on peut déduire la majoration
\begin{equation}
\frac{\lambda_{1,1}(N,h)}{2(k+1)}<\varepsilon.
\end{equation}
On peut choisir $\varepsilon$ suffisamment petit pour éliminer ce cas.

Si $\|\de \alpha\|^2\geq\|\codiff\alpha\|^2$, alors $\frac{\|\de
\alpha\|^2}{\|\alpha\|^2}\geq\frac{\lambda_{1,1}(N,h)}2$. On
peut déduire de (\ref{adap:eq1}) que
\begin{equation}\label{adap:eq2}
\|\de\alpha+\sum_{i=1}^k a_i e_i\|^2<\varepsilon.
\end{equation}
On veut montrer à partir de (\ref{adap:eq2}) que puisque 
$\de\alpha$ est minoré, les $a_i$, et donc les $\de a_i$ le sont
aussi, ce qui contredira le fait que $R(\varphi)<\varepsilon$. 
Il découle de (\ref{adap:eq2}):
\begin{equation}
\|\de\alpha\|^2+\|\sum_{i=1}^k a_i e_i\|^2+2(\de\alpha,
\sum_{i=1}^k a_i e_i)<\varepsilon
\end{equation}
et donc
\begin{equation}
\frac{\lambda_{1,1}(N,h)}{2(k+1)}-\varepsilon<-2
(\de\alpha,\sum_{i=1}^k a_i e_i),
\end{equation}
De plus, en utilisant le fait que $e_i$ est harmonique,
donc cofermée et qu'en restriction aux $1$-formes on a 
$\codiff=*\de*$, on obtient $(\de\alpha,a_i e_i)=
(\alpha,\codiff(a_i e_i))=(\alpha,*(\de a_i\wedge(*e_i)))$.
En appliquant ponctuellement l'inégalité de Schwarz, on en déduit
\begin{eqnarray}
(\de\alpha,a_i e_i)&\leq&
\int_N|\alpha||\de a_i\wedge(*e_i)|\de v\leq
\int_N|\alpha||\de a_i||e_i|\de v\nonumber\\
&\leq&\|e_i\|_\infty\int_N|\alpha||\de a_i|\de v = \|e_i\|_\infty
(|\alpha||\de a_i|)\nonumber\\
&\leq&\|e_i\|_\infty\|\alpha\|\|\de a_i\|\leq
\|e_i\|_\infty\|\de a_i\|
\end{eqnarray}
et finalement
\begin{equation}
\frac{\lambda_{1,1}(N,h)}4-\varepsilon<2\sum_{i=1}^k\|e_i\|_\infty
\|\de a_i\|.
\end{equation}
Comme selon le lemme \ref{geom:oneill}, $\|e_i\|_\infty$ est majoré en
fonction de $a$ et $n$, on obtient si $\varepsilon$ est
suffisamment petit une minoration de 
$\sum_{i=1}^k\|\de a_i\|$, ce qui contredit le fait que
$\sum_{i=1}^k\|\de a_i\|^2<\varepsilon$.

En choisissant $\varepsilon$ suffisamment petit en fonction
de $\lambda_{1,1}(N,h)$, $a$ et $n$, on peut donc écarter le
cas $\|\alpha\|^2\geq\frac1{1+k}$.

Supposons $\|a_p\|^2\geq\frac1{1+k}$ :

Si on note $\bar a_i$ la valeur moyenne de la fonction $a_i$, on 
peut choisir la base $(\omega_i)$ de sorte que $\bar a_i=0$ pour
$i\geq2$ (il suffit de choisir comme nouveau $\omega_1$ la forme
différentielle $\sum_{i=1}^k\bar a_i\omega_i$). 
On a alors 
\begin{equation}
\|\de a_i\|^2\geq\lambda_{0,1}(N,h)\|a_i\|^2 
\end{equation}
pour $i\geq2$ et 
\begin{equation}
\|\de a_1\|^2\geq\lambda_{0,1}(N,h)\|a_1-\bar a_1\|^2.
\end{equation}
Comme $\|\de a_i\|^2< \varepsilon$ pour tout $i$, on a donc
\begin{equation}
\|a_i\|_2^2<\frac\varepsilon{\lambda_{0,1}(N,h)}
\end{equation}
pour $i\geq2$ et 
\begin{equation}
\|a_1-\bar a_1\|_2^2<\frac\varepsilon{\lambda_{0,1}(N,h)}.
\end{equation}
En particulier, si $\varepsilon$ est suffisamment petit, on a $p=1$.

On peut alors écrire 
\begin{eqnarray}
\|\de\alpha+\sum_{i=1}^ka_ie_i\|^2&=&\|\bar a_1e_1\|^2
+\|\de\alpha+(a_1-\bar a_1)e_1+\sum_{i=2}^ka_ie_i\|^2\nonumber\\
&&+2(\bar a_1e_1,\de\alpha+(a_1-\bar a_1)e_1+\sum_{i=2}^ka_ie_i).
\end{eqnarray}
On a d'une part
$$(\bar a_1e_1,\de\alpha)=\bar a_1(\codiff e_1,\alpha)=0,$$
et d'autre part
\begin{eqnarray}
(\bar a_1e_1,a_ie_i)&\leq&
\|\bar a_1 e_1\|_2\cdot\|a_ie_i\|_2\leq
\|\bar a_1 e_1\|_2\|e_i\|_\infty\|a_i\|_2\nonumber\\
&<&\|\bar a_1 e_1\|_2
\frac{\sqrt\varepsilon\|e_i\|_\infty}{\lambda_{0,1}(N,h)},
\end{eqnarray}
cette dernière inégalité restant vraie pour $i=1$ en remplaçant
$a_i$ par $(\bar a_1-a_1)$.

Comme $\|\de\alpha+\sum_{i=1}^ka_ie_i\|^2<\varepsilon$, on a finalement
$$\|\bar a_1e_1\|^2<2|
(\bar a_1e_1,\de\alpha+(a_1-\bar a_1)e_1+\sum_{i=2}^ka_ie_i)|
+\varepsilon$$
et donc 
\begin{equation}
\|\bar a_1e_1\|^2-\frac{2n\sqrt\varepsilon\|e\|}{\lambda_{0,1}(N,h)}
\|\bar a_1e_1\|-\varepsilon<0.
\end{equation}
On en déduit que $\|\bar a_1e_1\|$ est encadré par les racines du
polynôme $X^2-\frac{2n\sqrt\varepsilon\|e\|}{\lambda_{0,1}(N,h)}X-
\varepsilon$, et en particulier majoré par sa plus grande racine.
Comme $\|e\|$ est majoré en fonction de $n$ et de la borne
sur la courbure (lemme \ref{geom:oneill}), 
la plus grande racine du polynôme s'écrit comme
une constante dépendant de $n$, $a$ et $\lambda_{0,1}(N,h)$, multipliée par
$\sqrt\varepsilon$. On obtient donc le résultat souhaité en
prenant comme $\omega$ la forme $\bar a_1\omega_1$.
\hfill\carrenoir

 On va maintenant minorer le spectre de $e^*e$ en fonction du
volume de la fibre $T^k$.

\begin{fait}\label{adap:ft2}
Il existe une constante $c(n,a,(N,h))>0$ telle que
la première valeur propre de $e^*e$ soit minorée par 
$c\cdot\Vol(T^k)^2$.
\end{fait}

\textbf{Démonstration :}
 Soient $\lambda_1,\cdots,\lambda_k$ les valeurs propres de $e^*e$
classées dans l'ordre croissant. Comme $e$ est injective, ces
valeurs propres sont non nulles et la premi\`ere vérifie
\begin{equation}
\lambda_1=\frac{\prod_i\lambda_i}{\prod_{i\neq1}\lambda_i}
=\frac{\Det(e^*e)}{\prod_{i\neq1}\lambda_i}.
\end{equation}
Par ailleurs, les valeurs propres de $e^*e$ v\'erifient
$\lambda_i\leq\|e^*e\|$, donc
\begin{equation}
\lambda_1\geq\frac{\Det(e^*e)}{\|e^*e\|^{k-1}}\geq\frac{\Det(e^*e)}
{\|e\|^{2k-2}}.
\end{equation}

L'image de $e$ est un sous-espace de $\Ker \Delta^2(N)$ de
dimension $k$, engendré par un sous-réseau du réseau des
formes harmoniques entières de $N$. Si on restreint $e$ et
$e^*$ à ce sous-espace, on peut écrire 
$\Det(e^*e)=(\Det e)^2$, où $\Det e$ est le déterminant d'une
matrice de $e$ écrite dans des bases orthonormées de 
$\mathcal G^*$ et $\Ima e$, ce qui donne 
\begin{equation}
\lambda_1\geq\frac{(\Det e)^2}{\|e\|^{2k-2}}.
\end{equation}

 Le lemme \ref{geom:oneill} donne une majoration de $\|e\|$ en fonction
de $n$ et de la borne $a$ sur la courbure de $M$, il ne reste 
donc qu'\`a minorer $\Det e$. Notons $\Det' e$ le déterminant de la
matrice de $e$ dans la base canonique de $\mathcal G^*={\R^k}^*$ 
et une base orthonormée de $\Ima e$. On a alors
$\Det e=(\Det' e)(\Vol T^k)$. Comme les images dans $\Ker\Delta^2(N)$
des éléments de la base canonique de $\mathcal G^*$ sont des formes 
entières, le déterminant $\Det' e$, qui est aussi le volume de 
$e([0,1]^k)$, est un multiple du volume d'un domaine fondamental 
du réseau des formes entières dans $\Ima e$. Comme par ailleurs
$\Det' e$ est non nul, il sera donc minoré par le volume
de ce domaine fondamental. Si on note $\rho$ le minimum des normes
des $2$-formes harmoniques entières non nulles, ce volume 
est minoré par le volume d'une boule de rayon $\frac\rho2$ dans 
$\Ima e$, et donc minoré par une constante ne dépendant que de 
$n$ et de la métrique $h$ de $N$. On peut donc bien écrire
$$\lambda_1\geq c(n,a,(N,h))\cdot\Vol(T^k)^2.$$

\hfill\carrenoir

 Nous allons maintenant supposer que $e$ n'est pas injective. Notons
$l$ la dimension de son noyau. Le premier nombre de Betti de $M$ est
alors $b_1(N)+l$. En effet, on a vu que si une $1$-forme $\varphi=
\pi^*(\alpha)+\sum_{i=1}^k \pi^*(a_i)\cdot\omega_i$ est harmonique, cela
signifie, d'après (\ref{adap:dphi}) et (\ref{adap:deltaphi}) :
\begin{equation}
\Delta\alpha=0,\de a_i=0\text{ pour tout $i$, et }\sum_{i=1}^k a_i\cdot e_i=0.
\end{equation}
 Comme les fonctions $a_i$ sont constantes. L'ensemble des $a_i$ tels
que
$\sum_{i=1}^k a_i\cdot e_i=0$ est exactement le noyau de $e$.
L'espace des formes harmoniques de $M$ est donc l'espace
engendré par les relevés des formes harmoniques de $N$ et les
formes verticales induites par les éléments de noyau de $e$.

 On peut reprendre la démonstration précédente en prenant pour
$(\omega_i)_i$ une base de $\mathcal G^*$ telle que $\omega_{k-l+1},
\cdots,\omega_k$ soit une base de $\Ker e$ (le fait que la
forme $\varphi$ est orthogonale aux formes harmoniques se traduit
par le fait que $a_{k-l+1},\cdots,a_k=0$) et en étudiant $e^*e$
restreint à l'orthogonal de $\Ker e$. On obtient de la même façon
le résultat du fait \ref{adap:ft1}, à savoir que la première valeur 
propre du laplacien sur $M$ est minorée à une constante multiplicative 
près par la première valeur propre de $(e^*e)_{|(\Ker e)^\perp}$.

Pour minorer le spectre de $(e^*e)_{|(\Ker e)^\perp}$, on doit
être un peu plus attentif dans la manipulation des bases de 
$\mathcal G^*$.

 Soit $\mathcal B=(\omega_1,\cdots,\omega_k)$ une base orthonormée de 
$\mathcal G^*$ et $\mathcal B'=(\omega_1',\cdots,\omega_k')$ une
base du réseau des entiers de $\mathcal G^*$, telles que 
$(\omega_1,\cdots,\omega_l)$ et $(\omega_1',\cdots,\omega_l')$ soient
des bases de $\Ker e$ (comme l'image du réseau des entiers de $\mathcal G^*$
est contenue dans un réseau, le noyau de $e$ est effectivement engendré
par des éléments entiers). La matrice de passage de $\mathcal B$ à 
$\mathcal B'$ est de la forme $$P=\left(\begin{array}{cc}P_1 & P_2\\
0 & P_3\end{array}\right),$$ où $P_1$ est un bloc carré de taille $l$.
Si se donne une base orthonormée de $\Ima e$, la matrice de $e$ 
s'écrit sous la forme $(0,A)$ dans la base $\mathcal B$ et 
$(0,A')$ dans la base $\mathcal B'$, où $A$ et $A'$ sont des
blocs carrés de taille $k-l$ et vérifient $A'=AP_3$.

Le spectre de $(e^*e)_{|(\Ker e)^\perp}$ est celui de $A^*A$. On peut 
écrire, comme dans la démonstration du fait \ref{adap:ft2}~:
\begin{equation}
\lambda_1\geq\frac{\Det A^*A}{\|A^*A\|^{k-l-1}}
\geq \frac{(\Det A)^2}{\|A\|^{2(k-l-1)}},
\end{equation}
où $\lambda_1$ est la première valeur propre non nulle de $e^*e$.
De plus, on a $\Det A'=\Det A\cdot\Det P_2$ et donc
\begin{equation}
\Det A=\frac{\Det A'}{\Det P_3}=\Det A'\frac{\Det P_1}{\Det P}.
\end{equation} 
Le déterminant de $A'$ est, comme précédemment, minoré par le covolume
du réseau des formes entières dans $\Ima e$, et $\Det P$ 
s'interprète géométriquement comme l'inverse du volume de $T^k$. 

Il reste à minorer $\Det P_1$. Comme $\Ker e$ est engendré par des
éléments entiers de $\mathcal G^*$, l'orthogonal de $\Ker e$ pour la
dualité définit un sous-tore $T^{k-l}$ de $T^k$. De plus, le dual de 
de l'algèbre de Lie $\mathcal G(T^k/T^{k-l})$ du quotient $T^k/T^{k-l}$
est isomorphe à $\Ker e$. La matrice $P_1$ est donc la matrice de passage
d'une base orthonormée de $\mathcal G^*(T^k/T^{k-l})$ dans une base du
base orthonormée de $\mathcal G^*(T^k/T^{k-l})$ dans une base du
réseau des entiers de $\mathcal G^*(T^k/T^{k-l})$, et par conséquent
$\Det P_1$ est l'inverse du volume de $T^k/T^{k-l}$ pour la métrique
quotient. Le diamètre de $T^k/T^{k-l}$ est majoré par $\varepsilon$,
comme celui de $T^k$, et par conséquent son volume aussi.\hfill\carrenoir
\section{Petites valeurs propres et norme des $2$-formes harmoniques 
entières}
 Nous allons ici discuter de la valeur de la constante $c(n,a,d,(N,h))$
du théorème \ref{adap:th}, et en particulier de la manière dont elle
dépend de la métrique $h$ sur $N$.

 Dans la démonstration du théorème, la géométrie de $N$ intervient quatre
fois: on a besoin de contrôler sa courbure pour appliquer la formule
de O'Neill et majorer les $\|e_i\|$; en \ref{adap:ft1} apparaissent
les valeurs propre $\lambda_{1,1}(N,h)$ et $\lambda_{0,1}(N,h)$; enfin
on fait intervenir en \ref{adap:ft2} le minimum des normes des
2-formes harmoniques non nulles dont la classe de cohomologie est entière.

 On peut noter que dans la démonstration du fait \ref{adap:ft1}, on
a seulement besoin d'une minoration des deux valeurs propres
$\lambda_{1,1}(N,h)$ et $\lambda_{0,1}(N,h)$. L'idée est en fait
de s'assurer que le spectre de $(N,h)$ n'interfère pas dans la
recherche des petites valeurs propres de $M$. On sait par ailleurs
que $\lambda_{0,1}(N,h)$ peut être minoré en fonction du diamètre
et de la courbure de $N$ (cf. théorème \ref{intro:th2}), 
et $\lambda_{1,1}(N,h)$ en fonction du diamètre, de la courbure et
du rayon d'injectivité de $N$ (théorème \ref{intro:th3}). En outre,
une borne sur la courbure est aussi suffisante pour appliquer la formule
de O'Neill.

La démonstration du fait \ref{adap:ft2} introduit quand à elle
la constante
\begin{equation}
\rho(N,h)=\inf_{\stackrel{\alpha\in\mathcal H^2(N,h)\backslash\{0\}}
{[\alpha]\in H^2(N,\Z)}}\|\alpha\|_2
\end{equation} 
dans la minoration de la première valeur propre du laplacien.
Il est naturel de se demander si l'on peut contrôler $\rho(N,h)$
à l'aide des mêmes invariants géométriques que $\lambda_{0,1}(N,h)$
et $\lambda_{1,1}(N,h)$:
\begin{question}\label{adap:q1}
Existe-t-il une constante $c(n,a,d,r)>0$ telle que si $(N^n,h)$ est une
variété riemannienne vérifiant $\diam(N,h)\leq d$ et $|K(N,h)|\leq a$
et $\injrad(N,h)\geq r$, alors $\rho(N,h)\geq c$ ?
\end{question}

Un argument de compacité permet de montrer qu'avec
l'hypothèse de rayon d'injectivité minorée, on peut répondre
affirmativement à la question \ref{adap:q1}.
\begin{pr}
Pour tout réels $a,d,r>0$ et tout entier $n\in\N^*$, il
existe une constante $c(n,a,d,r)>0$ (non explicite) telle que si 
$(N^n,h)$ est une variété riemannienne vérifiant $\diam(N,h)\leq d$ 
et $|K(N,h)|\leq a$ et $\injrad(N,h)\geq r$, alors $\rho(N,h)\geq c$.
\end{pr}
\textbf{Démonstration :}
 On considère le tore $T=T^{b_2(N)}$, vu comme quotient de 
$\mathcal H^2(N,h)$ par le réseau des formes harmoniques entières.
La métrique $h$ sur $N$ induit une norme euclidienne sur $\mathcal H^2(N,h)$
qui passe au quotient sur $T$ en une métrique plate $\bar h$. Minorer
$\rho(N,h)$ revient à minorer le rayon d'injectivité de $(T,\bar h)$.

On sait (\cite{ach}) que l'espace des métriques $h$ sur $N$ telles
que $\diam(N,h)\leq d$ et $|K(N,h)|\leq a$ et $\injrad(N,h)\geq r$
est relativement compact pour la topologie $C^\alpha$. 
De plus, la métrique $\bar h$ dépend continument de $h$. En effet,
si on se donne un réel $\varepsilon>0$ et une métrique $h$ sur $N$,
on aura, pour toute métrique $h'$ suffisamment proche de $h$ et toute
$2$-forme $\alpha$ harmonique pour la métrique $h$, 
$\left|\|\alpha\|_2-\|\alpha\|'_2\right|\leq \varepsilon\|\alpha\|_\infty$, 
où $\|\cdot\|_p$ et $\|\cdot\|'_p$ désignent les normes pour les métriques 
$h$ et $h'$ respectivement. Comme à courbure et diamètre bornés, la norme
$L^\infty$ des formes harmoniques est contrôlée par leur norme
$L^2$ (cf. \cite{li} et inégalité (\ref{geom:li})), on peut écrire
$\left|\|\alpha\|_2-\|\alpha\|'_2\right|\leq\tau(n,a,d)
\varepsilon\|\alpha\|_2$. Si $\alpha'$ est le représentant harmonique
pour $h'$ de la classe de cohomologie de $\alpha$, on peut donc finalement
trouver un voisinage $\mathcal V$ de $h$ tel que pour toute métrique
$h'$ dans $\mathcal V$,
\begin{equation}
\|\alpha'\|'_2\leq\|\alpha\|'_2\leq(1+\varepsilon)\|\alpha\|_2.
\end{equation}
Quitte à restreindre le voisinage $\mathcal V$, on a réciproquement 
\begin{equation}
\|\alpha\|_2\leq(1+\varepsilon)\|\alpha'\|'_2,
\end{equation}
ce qui implique bien la continuité de $h\rightarrow \bar h$.
L'ensemble décrit par $\bar h$ quand $h$ varie est donc relativement 
compact dans l'espace des métriques plates du tore. Il existe par
conséquent une métrique
sur $T$ qui réalise la borne inférieure du rayon d'injectivité de $T$
quand $h$ varie. En particulier, cette borne inférieure est non nulle
\hfill\carrenoir

 On peut se demander si le résultat reste vrai avec des hypothèses
plus faibles :
\begin{question}\label{adap:q3}
Existe-t-il une constante $c(n,a,d)>0$ telle que si $(N^n,h)$ est une
variété riemannienne vérifiant $\diam(N,h)\leq d$ et $|K(N,h)|\leq a$
et alors $\rho(N,h)\geq c$ ?
\end{question}

Il faut noter par ailleurs qu'une minoration non explicite ne permet 
pas d'améliorer les minorations déjà connues de la première valeur propre 
du spectre. On a besoin d'estimations précises :
\begin{question}\label{adap:q2}
Si $\diam(N,h)\leq d$, $|K(N,h)|\leq a$ et $\injrad(N,h)\geq r$, 
peut-on minorer $\rho(N,h)$ par une constante explicite $c(n,a,d,r)>0$ ?
Plus précisément, peut-on trouver une constante $c$ de la forme
$c'(n,a,d)\cdot\Vol(N,h)^{\alpha(n)}$ ou 
$c'(n,a,d)\cdot\injrad(N,h)^{\alpha(n)}$ ?
\end{question}
 On peut noter qu'expliciter le rôle du volume de $N$ dans cette
minoration permet d'obtenir dans le théorème \ref{adap:th} une minoration
de $\lambda_{1,1}(M,g)$ en fonction du volume de $(M,g)$.

\renewcommand{\bibname}{Références}

\end{document}